\documentclass[11pt]{article}
\usepackage[top=1.1in, left=1in, right=1in, bottom=1.1in]{geometry}
\usepackage{amsmath,amssymb,amsthm,amsfonts,verbatim}
\usepackage{enumerate}
\usepackage{microtype}
\usepackage{url}
\usepackage{algpseudocode}
\usepackage{tikz}
\usetikzlibrary{arrows}
\usepackage{amsmath,amssymb}
\usepackage[justification=centering]{caption}
\usepackage[justification=centering]{subcaption}
\usepackage[english]{babel}
\usepackage{array}
\usepackage{psfrag}
\usepackage{graphicx}
\usepackage{color}
\usepackage{setspace}
\usepackage{hyperref}
\usepackage{algorithm2e}
\usepackage[utf8]{inputenc}
\usepackage{cleveref}
\usepackage{mathrsfs} 
\usepackage{authblk}
\usepackage{threeparttable}
\usepackage[titletoc,title]{appendix}
\usepackage{xcolor}
\usepackage{xfrac}
\usepackage{caption}
\SetKw{KwInitialization}{Initialization}
\SetKw{KwTheMainBody}{The main body}

\newtheorem{theorem}{Theorem}[section]

\newtheorem{definition}[theorem]{\bf Definition}

\newtheorem{example}{Example}[section]
\definecolor{red}{rgb}{1,0.2,0.2}

\def \mr {\mathbb{R}}
\def \U {\mathcal{U}}

\makeatletter
\newcommand*{\rom}[1]{\expandafter\@slowromancap\romannumeral #1@} 
\makeatother

\hypersetup{colorlinks=true,linkbordercolor=red,linkcolor=blue,pdfborderstyle={/S/U/W 1}}

\title{Computing the quasipotential for nongradient SDEs in 3D}%

\author[1]{Shuo Yang\thanks{shuoyang@math.umd.edu}}
\author[2]{Samuel F. Potter\thanks{sfp@umiacs.umd.edu}}
\author[1]{Maria K. Cameron\thanks{cameron@math.umd.edu} }
\affil[1]{Department of Mathematics, University of Maryland, College Park, MD 20742, USA}
\affil[2]{Department of Computer Science, University of Maryland, College Park, MD 20742, USA}                               
  
\begin{document}

\maketitle

\begin{abstract}
Nongradient SDEs with small white noise often arise when modeling biological and ecological time-irreversible processes.  
If the governing SDE were gradient, the maximum likelihood transition paths, transition rates, 
expected exit times, and the invariant probability distribution would be given in terms of its potential function. 
The quasipotential plays a similar role for nongradient SDEs.
Unfortunately, the quasipotential is the solution of a functional minimization problem 
that can be obtained analytically only in some special cases.
We propose a Dijkstra-like solver for computing the quasipotential on regular rectangular meshes in 3D. 
This solver results from a promotion and an upgrade of the previously introduced 
ordered line integral method with the midpoint quadrature rule for 2D SDEs.
The key innovations that have allowed us to keep the CPU times reasonable while maintaining good accuracy are 
$(i)$ a new hierarchical update strategy, 
$(ii)$ the use of Karush-Kuhn-Tucker theory for rejecting unnecessary simplex updates, and
$(iii)$ pruning the number of admissible simplexes and a fast search for them.  
An extensive numerical study is conducted on a series of linear and nonlinear examples where the quasipotential
is analytically available or can be found at transition states by other methods. 
In particular, the proposed solver is applied to  Tao's examples  where the transition states are 
hyperbolic periodic orbits, and to a genetic switch model by Lv et al. (2014). 
The C source code implementing the proposed algorithm is available at M. Cameron's web page.
\end{abstract}


\section{Introduction}
Nongradient stochastic differential equations (SDEs) with small white noise accounting for random environmental factors
 often arise in modeling biophysical \cite{lv} and ecological  \cite{nolting,savanna}  processes. 
It is argued by Q. Nie and collaborators that noise plays a fundamental role in 
biological processes such as cell sorting and boundary formation between different gene expression domains \cite{NieNoise}.
The development of computational tools for the study of the stability of attractors of SDEs
and noise-induced transitions between them are important  for understanding the dynamics of such systems.
In this work, we present a numerical algorithm for computing the quasipotential for nongradient SDEs in 3D,
a key function of large deviation theory (LDT) \cite{FW},
giving asymptotic estimates for the dynamics of the system in the limit as the noise term
tends to zero. This algorithm is the result of an enhancement 
of the 2D ordered line integral method (OLIM) with the midpoint quadrature rule  \cite{olim} to 3D. We will refer to it as {\tt olim3D}.

\subsection{Definition and significance of the quasipotential}
We consider a system evolving according to the SDE
\begin{equation}
\label{sde}
dx = b(x)dt + \sqrt{\epsilon}dw, \quad x\in \mr^d,
\end{equation}
where $b(x)$ is a continuously differentiable vector field, $w$ is the standard Brownian motion, and $\epsilon>0$ is a small parameter.
We assume that $b$ is such that the corresponding deterministic system $\dot{x} = b(x)$ has a finite number of attractors, any trajectory of $\dot{x}=b(x)$ 
remains in a bounded region as $t\rightarrow\infty$, and
almost all trajectories approach some attractor $A$ as $t\rightarrow \infty$. 

The dynamics according to SDE \eqref{sde}, no matter how small $\epsilon$ is, are qualitatively different from those of $\dot{x}=b(x)$.
The system will escape from any bounded neighborhood  of any attractor $A$ of $\dot{x}=b(x)$ 
with probability close to one if you wait long enough.
As $\epsilon$ approaches $0$, possible escapes 
from the basin of attraction of any attractor $A$ 
can be quantified by the function called the quasipotential \cite{FW}.
Once the quasipotential is computed, one can readily obtain a number of useful 
asymptotic estimates.
The maximum likelihood escape path can readily be found by numerical integration. The
expected escape time from the basin of $A$ and the invariant probability 
distribution near $A$ are determined by the quasipotential up to exponential order \cite{FW}.
Moreover, a sharp estimate for the expected escape time from a basin of attraction was developed for a common special case \cite{bouchet}.

The quasipotential with respect to an attractor $A$ is the function $U_A:\mr^n\rightarrow [0,\infty)$ given by
\begin{equation}
\label{Qpot}
U_A(x) = \inf_{\psi,L}\{S(\psi)~|~\psi(0)\in A,~\psi(L) = x\},
\end{equation}
where the infimum is taken over all  absolutely continuous paths $\psi$ starting in $A$ and ending at $x$, $L$ is the length of $\psi$, and
$S(\psi)$ is the geometric action defined by
\begin{equation}
\label{Gaction}
S(\psi) = \int_0^L(\|b(\psi)\|\|\psi'\| - b(\psi)\cdot\psi')ds.
\end{equation}
The geometric action is 
obtained from the Freidlin-Wentzell action
\begin{equation}
\label{FWaction}
S_{T_0,T_1}(\phi) = \frac{1}{2}\int_{T_0}^{T_1}\|\dot{\phi} -  b(\phi) \|^2dt
\end{equation}
 by analytic optimization with respect to time $t$ and 
reparametrization of the path $\phi$ by its  arclength \cite{hey1,hey2,FW}.
From now on, we will assume that the attractor $A$ is fixed and omit the subscript in the notation for the quasipotential: $U(x)\equiv U_A(x)$.

\subsection{Significance of computing the quasipotential in the entire region}
Numerical methods for finding transition paths, transition rates, and estimating the invariant probability measure in the limit of $\epsilon\rightarrow 0$
for processes evolving according to SDE \eqref{sde} have developed in two directions.

The first direction advanced the family of path-based methods that aim at finding the maximum likelihood transition paths 
(a.k.a. the minimum action paths (MAPs), or instantons)
by numerical minimization of the Freidlin-Wentzell action \eqref{FWaction}
 \cite{mam,zhou1,zhou2} or the geometric action \cite{hey1,hey2}. 
 Path-based methods are applicable in arbitrary dimensions, even for discretized PDEs (see e.g. \cite{hey1,NLSE}).
 They tend to work well in any dimension when 
 the length of the MAP exceeds the distance between its endpoints by a modest factor.  
However, various numerical issues arise when the true MAP exhibits complex behavior.
Path evolution can stall when the true MAP spirals near an attractor or a transition state, or if the MAP has a kink.
Improved path-based methods (e. g. \cite{tao,sunzhou}) 
have been proposed to tackle these issues in some common special cases.

While the slow convergence of path-based methods is in the process of being resolved, 
there is another issue that cannot be completely eliminated.
The output of path-based methods is always biased to the initial guess. 
The path can converge to a local minimum or to a stationary path in the path space.
Even worse, it might be unable to converge to the true MAP in principle. For example, imagine the case where the phase space is a cylinder.
For any two distinct points on a cylinder, 
there exist countably many paths indexed by the number of revolutions that cannot be continuously deformed one into another. 

Furthermore,  path-based methods do not give any information about the 
quasipotential beyond the path. Hence, even if the MAP is found, the information
about the width and the geometry of the transition channel is missing. 
Of course, one can introduce a mesh and find the quasipotential on it 
by integrating  the geometric action \eqref{Gaction} or the  Freidlin-Wentzell action \eqref{FWaction}
 along the MAPs connecting each mesh point with the attractor $A$, 
but this is computationally expensive and wasteful.

The second direction started in \cite{quasi} is associated with the development of 
methods for computing the quasipotential on meshes.
 In this work, we introduce
the first 3D quasipotential solver, while previous works featured 2D solvers \cite{quasi,Rjournal,olim,DC2}.
Although these methods are limited to low dimensions, they have important advantages over path-based methods. 

First, knowledge of the quasipotential in a region provides important visual 
information about the dynamics of the system.
The quasipotential found within a level set completely lying in the basin of attraction of $A$
gives an estimate for the invariant probability measure within this level set.
Besides finding the maximum likelihood exit path from the basin of $A$
and estimating the expected exit time (see Section \ref{sec:background}), 
one can infer the width and the geometry of the exit channel.
For example, our computation of the quasipotential for the Lv et al. genetic switch model \cite{lv}
presented in Section \ref{sec:gene} has shown that the dynamics of the system are virtually 
limited to a small neighborhood of a 2D manifold, while the transition channel is relatively wide within it.
Visualization of the quasipotential allows one to capture the most important information 
about the dynamics of the system in the vanishing noise limit in a single glance.
A comprehensive review of using the quasipotential for analysis of ecological models is given in \cite{nolting}.

Second, the maximum likelihood escape location from the basin of $A$ might be not known a priori.
The computation of the quasipotential readily gives it, because as soon as the boundary of the basin of 
$A$ is reached, the quasipotential remains constant along trajectories going to other attractors.
For example, the vector field $b$ can be very complicated, as it is in the Sneppen-Aurell genetic switch model in Lambda Phage 
\cite{Aurell}, so that the location of the transition state is hard to 
determine by setting $b$ to zero, but becomes apparent from the
computed quasipotential \cite{DC2}. An interesting example of a maximum likelihood escape path not associated with 
any special-trajectory-type transition state such as saddle or
unstable limit cycle was found in the FitzHugh-Nagumo system 
using the quasipotential computed in an entire region \cite{FitzHugh-Nagumo}.

Third,  the MAP found by numerical integration using the 
computed quasipotential is guaranteed to be the global minimizer of the
geometric action \eqref{Gaction}. 
This is important because only global minimizers
are the maximum likelihood transition paths in the vanishing noise limit.

The computation of the quasipotential in 2D and 3D can be extended to that on 2D and 3D 
manifolds embedded into higher dimensions. This will be useful when the dynamics of a system 
evolving according to SDE \eqref{sde} are virtually restricted to a neighborhood of a low-dimensional manifold. 
This methodology is currently under development and will be presented in the future.

In summary, while the computation of the quasipotential in the entire region can be done 
only in low dimensions, it is worthwhile as it provides more complete 
 information about the dynamics
of the system, along with a useful means of  visualizing it in a region.

\subsection{An overview of Dijkstra-like Hamilton-Jacobi and quasipotential solvers}
\label{sec:briefo}
The proposed 3D quasipotential solver {\tt olim3D} 
is an extension and an upgrade of the OLIM with the midpoint quadrature rule 
introduced in \cite{olim}.

The OLIMs \cite{olim} can be viewed as a further development of the ordered upwind method (OUM) \cite{OUM2001,OUM2003} for solving the 
Hamilton-Jacobi equation of the form 
\begin{equation}
\label{hj0}
F\left(x,\tfrac{\nabla u(x)}{\|\nabla u(x)\|}\right)\| \nabla u(x)\| = 1,~~{\rm where}~~0< F_{\min}\le F\le F_{\max}<\infty 
\end{equation} 
is the speed of front propagation in the normal direction. 
The OUMs are Dijkstra-like solvers that inherit their general structure from the
fast marching method (FMM) \cite{SethPNAS,SethSIAM,SethBook} for solving the eikonal equation, $F(x)\|\nabla u(x)\| = 1$, 
where the speed function $F$
is isotropic. The anisotropy of the speed function $F$ in \eqref{hj0} rendered the FMM inapplicable.
This issue was resolved in the OUM \cite{OUM2001,OUM2003} by increasing the radius of neighborhoods for updating mesh points from $h$ to
$\Upsilon h$, where $h$ is the mesh step size and $\Upsilon : =F_{\max}/F_{\min}$ is the anisotropy ratio.
While, in theory, the OUM can be implemented in any finite dimension, in practice, unfortunately,
it  has remained limited to 2D because of exceedingly large CPU times faced already in 3D. 
The same is true of the first quasipotential solver based on the OUM \cite{quasi}. 
The key issue in adjusting the OUM for computing the quasipotential was associated with unbounded $\Upsilon$.
A finite \emph{update radius} of $Kh$ where $K$ 
is an
appropriately chosen integer, 
was introduced \cite{quasi}, and it was shown that the  additional error decays
quadratically with mesh refinement.

The OLIMs \cite{olim} are algorithmically similar to the OUM-based quasipotential solver, 
but involve two important innovations that make them much more accurate
 and 
up to four times faster.
The numerical errors of the OLIMs, mainly error constants, 
were reduced
by two to three orders of magnitude by the direct local solution of the minimization problem \eqref{Qpot} using quadrature rules such as
the midpoint, trapezoid and Simpson's  rules instead of the upwind finite difference  scheme used by the FMM \cite{SethPNAS,SethSIAM,SethBook}
and the OUM \cite{OUM2001,OUM2003,quasi}.
The CPU times were decreased by reducing the number of  computationally expensive triangle updates that are 
very unlikely to give the sought minimizer. 
The optimal choice of  quadrature rule in terms of the relationship between the computational error 
and the CPU time was  the midpoint rule.
 
\subsection{A brief summary of main results}
\label{sec:briefsum}
 The main result of the present work is the first 3D quasipotential solver {\tt olim3D}
 implementing the OLIM with the midpoint rule for 
solving minimization problem \eqref{Qpot} on a regular rectangular mesh.
The C source code  is available at M. Cameron's web page \cite{mariakc}.

A straightforward extension of the 2D OLIM-Midpoint  \cite{olim} to 3D would still lead to large 
CPU times and render the solver unappealing. 
While keeping high accuracy, we have managed to reduce CPU times  dramatically: e.g., from a few days to a few hours on a $513^3$ mesh.
This was achieved by introducing the following technical innovations:
\begin{itemize}
\item First, we have promoted  and upgraded the hierarchical update strategy proposed in \cite{olim}. We have made it more radical.
In \cite{olim}, it was applied only to new {\sf Considered} points, while now it is used to update all {\sf Considered} points (Section \ref{sec:hierarchy}).
\item Second, we have applied the Karush-Kuhn-Tucker constrained optimization theory 
\cite{nocedal} to reject unnecessary simplex updates that are unlikely to result in an
interior point solution (Section \ref{sec:KKT}).
\item Third, we have pruned the number of simplex updates by reducing the set of admissible simplexes 
and implementing a fast search for them (Section \ref{sec:admissible}).
\end{itemize}

We have conducted an experimental study of {\tt olim3D} on linear and nonlinear examples with ratios of the 
magnitudes of the rotational and potential  components ranging from one to ten.  The quasipotential in these examples is 
available analytically. 
A recommendation for choosing  the update factor $K$, an important parameter of the OLIMs, is given. 
Surprisingly at  first glance, we find that  larger ratios of the magnitude of the rotational and the potential components of
the vector field can lead to faster convergence rates. We give an explanation of this phenomenon on a 2D model in Section \ref{sec:why}.
Effects of local factoring near asymptotically stable equilibria \cite{QV}  adjusted for OLIMs are studied in Section \ref{sec:ELF}.
It is found that local factoring reduces computational errors for linear SDEs but may or may not be beneficial for nonlinear ones.

We have applied {\tt olim3D} to two of  Tao's examples \cite{tao} where the value of the quasipotential is known analytically at
hyperbolic periodic orbits serving as transition states between two attractors. Finally, we have applied {\tt olim3D} 
to a genetic switch model  (Lv et al. \cite{lv}). Our computation shows that the dynamics of this 3D
system is limited to a neighborhood of a 2D manifold.

The rest of the paper is organized as follows. 
Background on LDT \cite{FW} and the quasipotential is given in Section \ref{sec:background}.
The  proposed solver is described in details in Section \ref{sec:olim3D}.
The numerical study of {\tt olim3D} and the aforementioned applications are presented in Section \ref{sec:tests}.
The results are summarized and perspectives are given in Section \ref{sec:conclusion}.

\section{Background on the quasipotential: useful facts and estimates}
\label{sec:background}
In this section, we set up notations and terminology and recap useful facts about the quasipotential  which we will refer to throughout the rest of the paper.
We also state some important formulas involving the quaspotential that can be evaluated using the output of {\tt olim3D}.

One can show \cite{FW,quasi} that the quasipotential is a viscosity solution \cite{visc,ishii} of the boundary value problem for the Hamilton-Jacobi equation
\begin{equation}
\label{HJ}
\|\nabla U(x)\|^2 + 2 b(x)\cdot\nabla U(x)= 0,\quad \left.U\right\vert_A = 0.
\end{equation}
Equation \eqref{HJ} implies that $\nabla U(x)$ is orthogonal to 
\begin{equation}
\label{rot}
l(x): = \frac{1}{2}\nabla U(x) + b(x),
\end{equation}
the \emph{rotational component} of $b(x)$. It is easy to see that $b(x)$  decomposes into 
$b(x) = -\tfrac{1}{2}\nabla U(x) + l(x)$. The term $ -\tfrac{1}{2}\nabla U(x)$ is called the \emph{potential component} of $b$.

The invariant probability measure within any level set of the quasipotential completely lying in the basin of attraction of $A$
 is approximated by \cite{FW}
 \begin{equation}
 \label{inv}
 \mu(x)\asymp \exp( -U(x)/\epsilon),\quad {\rm i.e.}\quad \lim_{\epsilon\rightarrow 0}\epsilon\log\mu(x) = -U(x).
 \end{equation}
 The invariant probability measure is of the form $Z^{-1}e^{-U(x)/\epsilon}$, where $Z$ is a normalization constant,
 if and only if $\nabla \cdot l(x) = 0$ (see, e.g., \cite{quasi}).

The characteristics of  \eqref{HJ} are the MAPs. 
A MAP from $A$ to $x$ can be readily found by 
numerical integration \cite{quasi} of
\begin{equation}
\label{ShootMAP}
\psi' = -\frac{b(\psi) + \nabla U(\psi)}{\|b(\psi) + \nabla U(\psi)\|},\quad \psi(0) = x.
\end{equation}
However, the quasipotential defined by  \eqref{Qpot} is not the unique solution of  \eqref{HJ} 
with the homogeneous boundary condition on an attractor \cite{ishii}.
For example,  \eqref{HJ} written for a linear SDE where $b(x) = Bx$, 
where $B$ is a square matrix with all eigenvalues having negative real parts, 
has as many 
solutions as there are invariant subspaces of the linear transformation associated with $B$. 

The expected escape time from the basin of attraction  
$\mathcal{B}(A)$ of  $A$ can be estimated up to  exponential order \cite{FW}:
\begin{equation}
\mathbb{E}[\tau_{\mathcal{B}(A)}]\asymp \exp\{\min_{y\in\partial\mathcal{B}(A)}U(y)/\epsilon\},~~{\rm i.e.}~~
\lim_{\epsilon\rightarrow 0}\epsilon\log\mathbb{E}[\tau_{\mathcal{B}(A)}] = \min_{y\in\partial\mathcal{B}(A)}U(y).
\end{equation}
One can find the maximum likelihood escape path 
by integrating  \eqref{ShootMAP} starting from 
$$
x_s= \arg\min_{y\in\partial\mathcal{B}(A)}U(y).
$$
If $A$ is an asymptotically stable equilibrium  $x^{\ast}$, if $x_s$ is a Morse index one saddle, and if
the quasipotential twice continuously differentiable 
near $x^{\ast}$ and $x_s$, then one can obtain a sharp estimate for the expected 
escape time from $\mathcal{B}(x^{\ast})$ using the Bouchet-Reygner formula \cite{bouchet} written for SDE \eqref{sde}:
\begin{equation}
\label{BR}
\mathbb{E}[\tau_{\mathcal{B}(x^{\ast})}] \approx 
\frac{2\pi}{{\lambda}_{+}} \sqrt{\frac{|\det ~H(x_s)|}{\det~ H(x^{\ast}) }} 
\exp\left(\int_{0}^{L} l(\psi(s)) ds\right) \exp\left(\frac{U(x_s)}{\epsilon}\right),
\end{equation}
where, ${\lambda}_{+}$ is the positive eigenvalue of the Jacobian of $b(x)$
at the saddle point $x_s$, 
$\epsilon$ is a small parameter, $H(x^{\ast})$ and $H(x_s)$ 
are the Hessian matrices of the quasipotential at the equilibrium $x^{\ast}$ and the saddle $x_s$ respectively,
and $l(x)$ is the rotational component of $b(x)$ given by  \eqref{rot}.


\section{Description of the algorithm}
\label{sec:olim3D}
As  pointed out in Section \ref{sec:briefsum}, the OLIMs 
are Dijkstra-like solvers that can be viewed as a development of the line of thought originating with the 
fast marching method (FMM) \cite{SethPNAS,SethSIAM,SethBook}, and followed by the ordered upwind method  (OUM) \cite{OUM2001,OUM2003}. 
These methods are contrasted in Table \ref{table:compare}.
\begin{table}[htp]
\caption{A comparison of Dijkstra-like solvers.}
\begin{center}
\begin{tabular}{|c|c|c|c|c|}
\hline
Year&Method&Problem&Update radius&Update rule\\
\hline
1996&FMM&$F(x)\|\nabla u\| = 1$&$h$& Upwind finite difference \\
\hline
2001 &OUM &$F(x,\tfrac{\nabla u}{\|\nabla u\|})\|\nabla u\| = 1$& $\frac{ F_{\max}}{F_{\min}}h$& Upwind finite difference \\
\hline 
2017 & OLIM & $U(x)=\inf_{\psi,L}S(\psi)$ & $Kh$, $1\ll K\ll N$ & Quadrature rules for $S(\psi)$\\
\hline
\end{tabular}
\end{center}
\label{table:compare}
\end{table}%

These methods belong to the family of 
label-setting algorithms, a comprehensive overview of which is given in \cite{CV}. 
The terminology for labels of mesh points used by the OLIMs is borrowed from the OUM: {\sf Unknown}, {\sf Considered}, {\sf Accepted Front}, and {\sf Accepted}.
\begin{itemize}
\item {\sf Unknown points:}  points where the solution $u$ has not been computed yet, 
and none of its nearest neighbors are {\sf Accepted} or {\sf Accepted Front}.
\item {\sf Considered:} points { that} have {\sf Accepted Front} nearest neighbors and that have 
tentative values of $u$ that might change as the algorithm proceeds.
\item {\sf Accepted Front:}  points at which $u$ has been computed and will no longer be updated. 
These points have at least one {\sf Considered} nearest neighbor, and 
are used to update {\sf Considered} points.
\item {\sf Accepted:} points whose $u$ has been computed and fixed, 
having only {\sf Accepted} or {\sf Accepted Front} nearest neighbors.
\end{itemize}
Breaking with the FMM and OUM, the OLIMs completely abandon the use of 
finite difference schemes. Instead, local functional minimization problems are solved at each step. 
The main contribution to the numerical error
in the OUM comes from the first order upwind finite difference scheme applied on obtuse triangles with two of their sides exceeding the mesh step size by a 
significant factor, bounded by the anisotropy ratio $\Upsilon$ in \cite{OUM2001,OUM2003}, or the update factor $K$ in \cite{quasi}. 
The OLIMs also compute updates using the same type of triangles. However, the use of the
higher order  (midpoint, trapezoid, or Simpson's) quadrature rules for approximating the functional significantly reduces the numerical error. 
The numerical experiments in \cite{olim} demonstrated error reduction by about three orders of magnitude 
and  superlinear error decay for practical mesh sizes ($N\times N$ where $N=2^p$, $p\le 12$). 

The proposed algorithm solves minimization problem \eqref{Qpot} with boundary condition $\left.U\right\vert_A = 0$. 
The initialization near asymptotically stable equilibria is discussed in Section \ref{sec:init}.
The initialization near stable limit cycles can be done as proposed in \cite{olim}. 
In this work, we focus on the computation of the quasipotential
with respect to  asymptotically stable equilibria as this is the most important case for applications. 
At each step, a {\sf Considered} mesh point $x_0$ with the smallest tentative value of the quasipotential changes its status to {\sf Accepted Front}.
Then, the following two series of updates are performed:
\begin{enumerate}
\item
All current {\sf Considered} points lying within a distance of $Kh$ from $x_0$ are updated using $x_0$ and 
its {\sf Accepted Front} nearest neighbors.
\item
The status of each {\sf Unknown} nearest neighbor $x$ of $x_0$ is changed to {\sf Considered} 
and $x$ is updated using all {\sf Accepted Front} points
lying within a distance of $Kh$ from $x$.
\end{enumerate}
Here, $h$ is the mesh step, and $K$ is the update factor. 
Since the anisotropy ratio \cite{OUM2001,OUM2003} is unbounded for  minimization problem \eqref{Qpot}, 
the update radius has to be forcefully truncated.
In the OUM-based quasipotential solver \cite{quasi} and the OLIMs \cite{olim} it is set to $Kh$ where $K$ is a user-chosen 
positive integer. A rule of thumb for choosing a good value of $K$  in 2D was proposed in \cite{olim}. We refer an interested reader to
the accompanying discussions in \cite{olim,quasi}. In this work, we  will tune the rule of thumb for the 3D case in Section \ref{sec:tests}. 

The update procedures are elaborated in Sections \ref{sec:updates} and \ref{sec:simplex}.
The computation terminates as soon as  a boundary mesh point becomes {\sf Accepted Front}.
The rationale for this is that the computation approximately follows the minimizers of geometric action \eqref{Gaction} which can return 
to the computational domain after escaping it.

The selection of admissible simplexes and a fast search for them will be discussed in Section \ref{sec:admissible}.
In Section \ref{sec:hierarchy}, the 3D hierarchical update strategy will be presented. 
The use of the Karush-Kuhn-Tucker theory (Chapter 12 in \cite{nocedal}) for skipping unnecessary simplex updates
will be explained in Section \ref{sec:KKT}.
Our implementation of local factoring will be described in Section \ref{sec:locfac}.

\subsection{Initialization near asymptotically stable equilibria}
\label{sec:init}
Let $x^{\ast}$ be an asymptotically stable equilibrium of $\dot{x}=b(x)$.
Let $J$ be the Jacobian matrix of $b(x)$ evaluated at $x=x^{\ast}$. 
Then the dynamics according to SDE \eqref{sde}
near $x^{\ast}$ are approximated by those of the linear SDE for the new variable $y: = x - x^{\ast}$:
\begin{equation}
\label{linsde}
dy = Jy + \sqrt{\epsilon}dw.
\end{equation}
Since $x^{\ast}$ is asymptotically stable, all eigenvalues of $J$ have negative real parts. 
The quasipotential for  \eqref{linsde} is a quadratic form $U(y) = y^\top Qy$ where $Q$ is a symmetric positive definite matrix
that needs to be found. We will call $Q$ the \emph{quasipotential matrix}. Plugging the gradient $\nabla U(y) = 2Qy$ to the
 Hamilton-Jacobi equation \eqref{HJ} and canceling the factor of 4 we obtain
\begin{equation}
\label{HJL}
y^\top Q\left(Q + J\right)y = 0\quad\text{for all}\quad y\in\mr^n.
\end{equation}
Hence the matrix $Q(Q+J)$ must be antisymmetric, i.e.
\begin{equation}
\label{anti}
Q(Q+J) + (Q+J^\top)Q = QJ + J^\top Q + 2Q^2 = 0.
\end{equation}
The solution of  \eqref{anti} is given by the Chen-Freidlin formula \cite{CF,chen}
\begin{equation}
\label{chen}
Q =  \left(\int_0^{\infty}e^{Jt}e^{J^\top t} dt\right)^{-1}.
\end{equation}
The direct evaluation of $Q$ from  \eqref{chen} is inconvenient.
Instead, one can observe that, multiplying  \eqref{anti} by $Q^{-1}$ 
on the left and on the right, one obtains the Sylvester equation (a.k.a. the Lyapunov equation) with respect to $Q^{-1}$:
\begin{equation}
\label{sylvester}
JQ^{-1} + Q^{-1} J^\top = -2I.
\end{equation}
 \eqref{sylvester} is solved by the Bartels-Stewart algorithm \cite{BS} that is implemented in MATLAB in the command {\tt sylvester}.
Therefore, the quasipotential matrix $Q$ can be calculated in MATLAB using the following command\footnotemark[1]:
\begin{verbatim}
Q = inv(sylvester(J,J',-2*eye(size(J))))
\end{verbatim}
\footnotetext[1]{We thank Prof. Daniel Szyld for pointing out this simple way to find the quasipotential decomposition for linear SDEs using MATLAB.}

It is convenient to set up the mesh so that $x^{\ast}$ is a mesh point. Then, the nearest 26 mesh points $x$ can be initialized by
$U(x) = (x - x^{\ast})^\top Q(x - x^{\ast})$.


\subsection{The three update types: one-point, triangle, and simplex}
\label{sec:updates}
Our implementation in {\tt olim3D} involves three types of updates: 
the one-point update and the triangle update  are similar to the ones in 2D \cite{olim}, 
and a new simplex update is added. 
Let $x$ be a {\sf Considered} point that is to be updated.
Let us imagine for a moment a surface $\Sigma$ surrounding the attractor $A$ and
passing through all {\sf  Accepted Front}  points. 
Then true value at $x$ would be 
\begin{equation}
\label{true}
U(x) = \inf_{y,\psi}\{U(y) + \int_0^L(\|b(\psi)\|\|\psi'\| - b(\psi)\cdot\psi')ds~|~\psi(0) = y\in\Sigma,~\psi(L) = x\}.
\end{equation}
The infimum in  \eqref{true} is achieved at $(x^{\star},\psi^{\star})$, where $\psi^{\star}$ is the segment of the MAP arriving at $x$ 
from the attractor $A$ cut off by the surface $\Sigma$, and $x^{\star} = \psi^{\star}(0)$ is the point where the MAP crosses the surface $\Sigma$. 
The OLIMs approximate 
the solution of  \eqref{true} by the minimum
among all straight line segments $[y,x]$ where $y\in \Sigma$ and $\|x - y\|\lesssim Kh$.

{\bf One-point update.}
The proposed value at $x$ from the {\sf Accepted Front} point $x_0$ is given by
\begin{equation}
\label{1u}
\mathsf{Q}_{1}(x_0,x) = U(x_0) + \mathcal{Q}(x_0,x),
\end{equation}
where $\mathcal{Q}(x_0,x)$ is the midpoint quadrature rule applied to the integral in  \eqref{true}:
\begin{equation}
\label{mid1}
\mathcal{Q}(x_0,x) = \|b_m\|\|x-x_0\| - b_m\cdot (x-x_0),\quad {\rm where}~~  b_m\equiv b\left( \frac{x+x_0}{2}\right).
\end{equation}

{\bf Triangle update.}
The proposed value at $x$ from the {\sf Accepted Front} points $x_0$ and $x_1$ is the solution of the following constrained minimization problem:
\begin{align}
\mathsf{Q}_{2}(x_0,x_1,x) &= \min_{\lambda\in[0,1]}\{U_{\lambda} + \mathcal{Q}(x_\lambda,x)\},~~{\rm where}\label{2u} \\
 U_{\lambda} &= U(x_0) + \lambda(U(x_1) - U(x_0)),~~x_{\lambda} = x_0 + \lambda(x_1- x_0), ~~{\rm and} \notag \\
 \mathcal{Q}(x_\lambda,x) &=  \|b_{m\lambda}\|\|x - x_{\lambda}\| - b_{m\lambda}\cdot (x -x_{\lambda}),\quad {\rm with}\notag\\
b_{m\lambda}&\equiv b_{m0} + \lambda (b_{m1}-b_{m0}),~~
b_{m0} \equiv b\left(\frac{x_0 + x}{2}\right), ~~ b_{m1} = b\left(\frac{x_1 + x}{2}\right).\notag
\end{align}

\begin{figure}[htbp]
\begin{center}
\includegraphics[width = 0.5\textwidth]{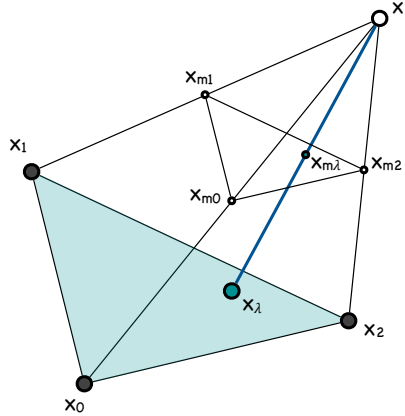}
\caption{An illustration for the simplex update corresponding to  minimization problem \eqref{3u}. 
The points $x_0$ , $x_1$, and $x_2$ are {\sf Accepted Front}. The point $x$ is to be updated. The points $x_{m0}$, $x_{m1}$,
and $x_{m2}$ are the midpoints of the edges of the simplex $(x_0,x_1,x_2,x)$ ending at $x$.
}
\label{fig:simplex}
\end{center}
\end{figure}

{\bf Simplex update.}
The proposed value at $x$ from the {\sf Accepted Front} points $x_0$, $x_1$, and $x_2$ 
is the solution of the following constrained minimization problem illustrated in Fig. \ref{fig:simplex}:
\begin{align}
\mathsf{Q}_{3}(x_0,x_1,x_2,x)& = \min\{U_{\lambda} + \mathcal{Q}(x_\lambda,x)\},~~{\rm where}\label{3u} \\
U_{\lambda} &= U(x_0) + \lambda_1(U(x_1) - U(x_0)) + \lambda_2 (U(x_1) - U(x_0)),\notag \\
x_{\lambda} &= x_0 + \lambda_1(x_1- x_0) + \lambda_2(x_2- x_0),\notag \\
 \mathcal{Q}(x_\lambda,x) & =  \|b_{m\lambda}\|\|x - x_{\lambda}\| - b_{m\lambda}\cdot (x -x_{\lambda}),\quad {\rm with}\notag\\
b_{m\lambda} & \equiv b_{m0} + \lambda_1 (b_{m1}-b_{m0}) + \lambda_2(b_{m2}-b_{m0}),~~
b_{mi}\equiv b\left(\frac{x_i + x}{2}\right),~~ i= 0,1,2, \notag \\
&\text{subject to} \quad
\lambda_1\ge 0,\quad
\lambda_2\ge 0,\quad
\lambda_1 + \lambda_2 \le 1.\notag
\end{align}

In all three updates, the proposed value $\mathsf{ Q}_i$, $i = 1,2,3$, 
replaces the current value $U(x)$ if and only if $\mathsf{Q}_i < U(x)$.

The solution of the  constrained minimization problem \eqref{2u}  for the triangle update in 3D is done in the same way as it is done in 2D. 
The function $f(\lambda) : = U_{\lambda} + \mathcal{Q}(x_\lambda,x)$ minimized in  \eqref{2u} is differentiated with respect to $\lambda$.  
If $f'(0)f'(1) < 0$  then
Wilkinson's hybrid method \cite{wilkinson} (a combination of the secant and bisection methods) is applied to find the root.
Otherwise, the minimum in \eqref{2u} is achieved at one of the end-points and must already have been 
 found by the routinely used one-point update.  Over 67\%  of triangle updates are 
rejected by this rule. 
The details of the triangle update are worked out in \cite{olim}.


\subsection{The simplex update} 
\label{sec:simplex}
The solution of the  constrained minimization problem for the simplex update  \eqref{3u}
is found using Newton's method. In {\tt olim3D}, whenever the simplex update is called,
there is a local minimum $y$ on the boundary of the simplex base
found by a triangle update. This point $y$ is used as a warm start for the iteration. 
Then, the Karush-Kuhn-Tucker (KKT) optimality conditions \cite{nocedal} (Chapter 12) 
are checked for $y$. If they are satisfied, a solution to  \eqref{3u} has already been found. 
Our experiments show that 50\% to 70\%  of 
the simplex updates are rejected by this KKT criterion. 
Otherwise, we try to find the solution of  \eqref{3u} using Newton's method.  
We will elaborate on the application of the KKT conditions in Section \ref{sec:KKT}.

The function to be minimized in  \eqref{3u} is
\begin{align}
f(\lambda) & : = U_{\lambda}  + \|x-x_{\lambda}\| \|b_{m\lambda}\| - (x-x_{\lambda})\cdot b_{m\lambda},~~{\rm where} \label{fmin}\\
b_{m\lambda}& :=b(x_{m0}) + \lambda_1\left[b(x_{m1})-b(x_{m0})\right] + \lambda_2\left[b(x_{m2}) - b(x_{m0})\right].\notag
\end{align}
The gradient $\nabla f$ and Hessian $H$ are given by
\begin{align}
\nabla f &= \delta U + \frac{\|x-x_{\lambda}\|}{ \|b_{m\lambda}\|} B^\top b_{m\lambda} + \frac{ \|b_{m\lambda}\|}{\|x-x_{\lambda}\|}X^\top (x-x_{\lambda}) 
- B^\top (x-x_{\lambda}) - X^\top b_{m\lambda}, \label{grad3}\\
H & =\frac{X^\top (x-x_{\lambda} ) [B^\top b_{m\lambda}]^\top +
\left[X^\top(x-x_{\lambda} )[B^\top b_{m\lambda}]^\top\right]^\top} {\|x-x_{\lambda}\| \|b_{m\lambda}\|}\notag\\
&+\frac{\|x-x_{\lambda}\|}{ \|b_{m\lambda}\|}B^\top B + \frac{ \|b_{m\lambda}\|}{\|x-x_{\lambda}\|}X^\top X\notag\\
&-\frac{\|x-x_{\lambda}\|}{ \|b_{m\lambda}\|^3} B^\top b_{m\lambda} \left[B^\top b_{m\lambda} \right]^\top
-\frac{ \|b_{m\lambda}\|}{\|x-x_{\lambda}\|^3}X^\top (x-x_{\lambda})\left[X^\top (x-x_{\lambda})\right]^\top\notag \\
&-\left(B^\top X + [B^\top X]^\top\right), \label{hess3}
\end{align}
where
\begin{align}
B &: = \left[b(x_{m1}) - b (x_{m0}),b(x_{m2}) - b (x_{m0})\right]~~({\rm a}~~3\times 2 ~~{\rm matrix}),\label{Bmatr}\\
X & : = \left[x_0 -x_1,x_0-x_2\right]~~({\rm a}~~3\times 2 ~~{\rm matrix}),\label{Xmatr}\\
\delta U &: = \left[\begin{array}{c}U(x_1) - U(x_0)\\U(x_2) - U(x_0)\end{array}\right].\notag
\end{align}
Some details of the calculation of $\nabla f$ and $H$ are worked out in Appendix A.


\subsubsection{Is the Hessian positive definite?}
Let us check if $H$ is positive definite.
First, note that we have not used the fact that $x_{m0}$, $x_{m1}$, and $x_{m2}$ are the midpoints of the simplex edges emanating from $x$.
Hence, the  gradient \eqref{grad3} and  Hessian \eqref{hess3} 
will have the same form for any choice of points $y_0\in[x_0,x]$, $y_1\in[x_1,x]$ and $y_2\in[x_2,x]$
to play the roles of 
$x_{m0}$, $x_{m1}$, and $x_{m2}$, respectively. 
In particular, if $y_0\equiv y_1\equiv y_2\equiv x$ as it is in the OLIM with the righthand quadrature rule \cite{olim}, 
$B$ is the zero matrix. Then the Hessian becomes
\begin{align}
H_R & =  \frac{ \|b(x)\|}{\|x-x_{\lambda}\|}X^\top X -
\frac{ \|b(x)\|}{\|x-x_{\lambda}\|^3}X^\top (x-x_{\lambda})\left[X^\top (x-x_{\lambda})\right]^\top\notag \\
& = \beta X^\top\left(I - pp^\top\right)X, \quad {\rm where}\quad 
\beta: = \frac{ \|b(x)\|}{\|x-x_{\lambda}\|},\quad
p:=\frac{x-x_{\lambda}}{\|x-x_{\lambda}\|}.
\label{HR}
\end{align}
The symmetric matrix $I-pp^\top$ in  \eqref{HR} is the projection matrix onto the plane normal to the unit vector $p$.
It has one zero eigenvalue corresponding to the eigenvector $p$ and two eigenvalues equal to one corresponding to the two-dimensional eigenspace
orthogonal to $p$. Hence $H_R$ is symmetric nonnegative definite. Furthermore, since $p\notin{\rm span}(X)$, 
there is no vector $\lambda\in\mr^2$ such that $X\lambda \parallel p$. Hence, $H_R$ 
is positive definite. 

For any other choice of $y_0$, $y_1$, and $y_2$, the situation is more involved.
In general, the Hessian is not positive definite. An example is presented in Appendix B.
However, as the mesh step size $h$ tends to zero, and if the update radius $K$ is proportional to $h^{-\alpha}$ where $0<\alpha<1$,
the Hessian becomes positive definite away from neighborhoods of equilibria of $b(x)$ (see Appendix B).

\subsection{Admissible triangles and simplexes and a fast search for them}
\label{sec:admissible}
The cost of a Dijkstra-like quasipotential solver in 3D is at least $O(N^3\log N)$, and the CPU times tend to be large 
for $N\sim 500$. Therefore, it is important to design time-efficient 3D codes that avoid unnecessary floating point operations.
One way to reduce the number of triangle and simplex updates is to limit the number of types of admissible simplexes while preserving
the directional coverage of characteristics. 

\begin{figure}[htbp]
\begin{center}
\includegraphics[width = 0.6\textwidth]{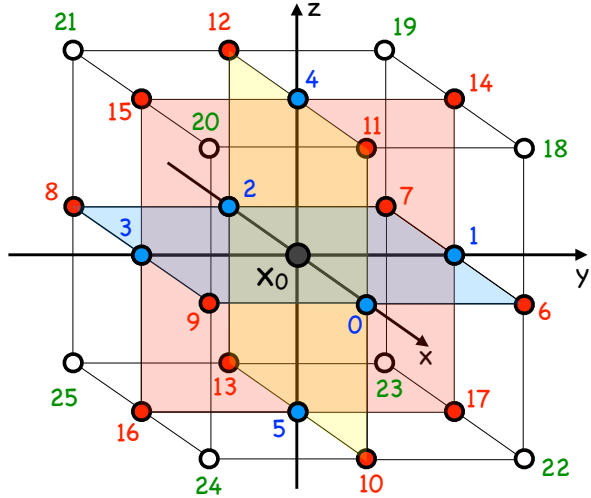}
\caption{An illustration for Definition \ref{def:neib}. The neighborhoods 
$\mathcal{N}_1(x_0)$, $\mathcal{N}_2(x_0)$ and $\mathcal{N}_3(x_0)$ consist of blue, red, and clear points respectively.}
\label{fig:neib}
\end{center}
\end{figure}

We allow for different mesh step sizes along $x$-, $y$-, and $z$-axes. 
Despite this, the definitions of the nearest neighborhoods and the far neighborhood below are given in terms of  indices of mesh points rather 
than  Euclidean distances between the mesh points. 
So, let us map our 3D mesh into the lattice $\mathbb{Z}^3$.
\begin{definition}
\label{def:neib}
Let  $x_0$ and $x$ be lattice points with indices  $(i_0,j_0,k_0)\in \mathbb{Z}^3$ and $(i,j,k)\in \mathbb{Z}^3$ respectively.
The $l_1$ and $l_{\infty}$ distances between $x_0$ and $x$ are  
$$
\|x - x_0\|_1 = |i-i_0| + |j - j_0| + |k - k_0| \quad{\rm and}\quad \|x - x_0\|_{\infty} = \max\{ |i-i_0|, |j - j_0| ,|k - k_0|\},
$$ 
respectively. Then, the nearest neighborhoods $\mathcal{N}_1(x_0)$, $\mathcal{N}_2(x_0)$, $\mathcal{N}_3(x_0)$, and $\mathcal{N}(x_0)$
are defined by
\begin{align}
\mathcal{N}_1(x_0) & := \{ x = (i,j,k)\in \mathbb{Z}^3~|~\|x - x_0\|_1 = 1\}\label{N1}\\
\mathcal{N}_2(x_0) & := \{ x = (i,j,k)\in \mathbb{Z}^3~|~\|x - x_0\|_1 = 2 ~{\rm and}~\|x - x_0\|_{\infty} = 1\}\label{N2}\\
\mathcal{N}_3(x_0) & := \{ x = (i,j,k)\in \mathbb{Z}^3~|~\|x - x_0\|_1 = 3 ~{\rm and}~\|x - x_0\|_{\infty} = 1\}\label{N3}\\
\mathcal{N}(x_0) &: =\mathcal{N}_1(x_0)\cup\mathcal{N}_2(x_0)\cup\mathcal{N}_3(x_0)\label{NN}.
\end{align}
The far neighborhood $\mathcal{N}_{{\rm far}}^K(x_0)$ of $x_0$ for update factor $K$ consists of all
 lattice points $x = (i,j,k)\in \mathbb{Z}^3$ excluding  $x_0$ such that
 \begin{align}
 |i - i_0| & \le K \label{far_i}\\
 |j - j_0| & \le {\rm ceil}\left(\sqrt{K^2 - |i - i_0|^2}\right), \label{far_j}\\
 |k - k_0| & \le {\rm ceil}\left(\sqrt{K^2 - \min\{|i - i_0|^2 +|j - j_0|^2,K^2\} }\right). \label{far_k}
\end{align}
\end{definition}
The far neighborhood of $x_0$ is slightly larger than the Euclidean ball of radius $K$.
Definition \ref{def:neib} is illustrated in Fig. \ref{fig:neib}. 
The blue, red, and clear lattice points constitute the neighborhoods 
$\mathcal{N}_1(x_0)$, $\mathcal{N}_2(x_0)$, and $\mathcal{N}_3(x_0)$ respectively. The sizes of these neighborhoods are 
$|\mathcal{N}_1(x_0)| = 6$, $|\mathcal{N}_2(x_0)| = 12$ and $|\mathcal{N}_3(x_0)| = 8$.
Altogether, they comprise the 26-point nearest neighborhood $\mathcal{N}(x_0)$. 

The triangle and simplex update rules given by Eqs. \eqref{2u} and \eqref{3u}, respectively, use linear interpolation 
to approximate the values of $U$ at the bases of  triangles and simplexes, respectively.
This incurs an interpolation error. It is clear that larger bases tend to lead to larger interpolation errors. 
To develop a quantitative insight, let us consider the linear interpolation of a strictly convex quadratic function
$$
u(x) = u_0 + g^\top x + \tfrac{1}{2}x^\top Hx,\quad x\in \mr^d,\quad d = 1,2.
$$
In 1D, the error of linear interpolation of $u(x)$ on an interval of length $h$ is at most $\tfrac{1}{8}h^2H$.
In 2D, the error of linear interpolation of $u(x)$
in a triangle $(x_0,x_1,x_2)$
depends on the mutual orientation of the level sets of $u(x)$ and the triangle,  and the ratio of the eigenvalues of $H$.
Nevertheless, the trend remains the same: larger triangles lead to larger errors. For example, if $H$ is the identity matrix,
and the vertex of the quadratic function $u$ lies at the center of the circumcircle of  $(x_0,x_1,x_2)$,
the interpolation error is at most $R^2/2$ where $R$ is the radius of the circumcircle. 
For a right isosceles triangle with legs of length $h$, the radius of the circumcircle is
$R = h/\sqrt{2}$, hence the error is bounded by $h^2/4$.
For an equilateral triangle with side length $h\sqrt{2}$, $R = h\sqrt{2/3}$ and the error is at most $h^2/3$.
Motivated by this, we limit the set of admissible simplexes to those whose bases are right isosceles triangles 
such that the vertices at the acute angles belong to the $\mathcal{N}_1$ neighborhood of the vertex at the right angle.
Respectively, the sides of the bases of admissible simplexes are the bases of admissible triangles; i.e.,
the bases of admissible triangles are pairs of $\mathcal{N}_1$ or $\mathcal{N}_2$ nearest neighbors. 

\begin{figure}[htbp]
\begin{center}
\includegraphics[width = 0.9\textwidth]{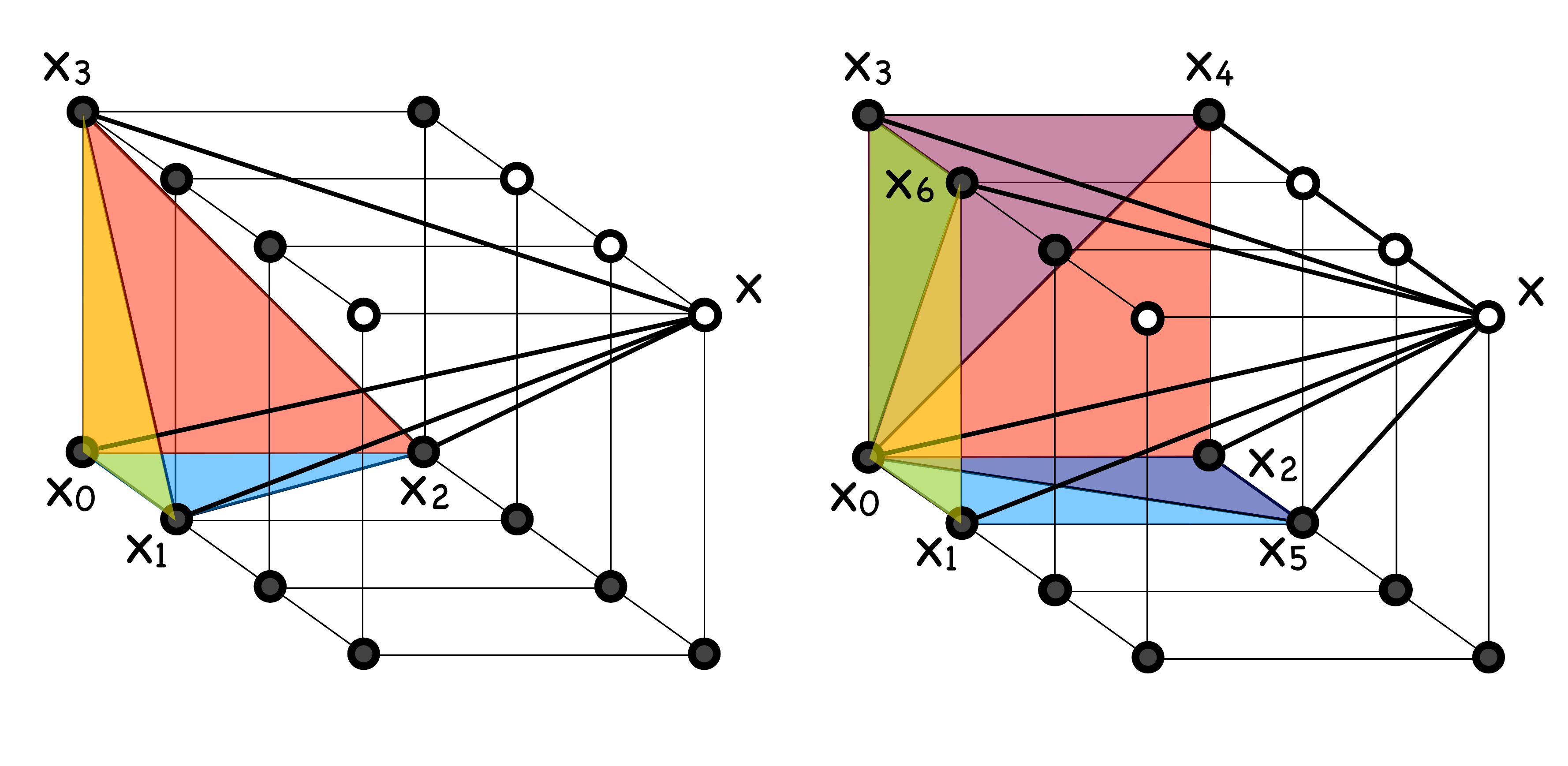}
\caption{The set of admissible simplexes with bases containing a given {\sf Accepted Front} point $x_0$
and lying in a chosen octant; $x$ is a {\sf Considered} point for which a simplex  update is being performed. 
Left: Three admissible simplexes with $x_0$ at the right angle of their bases:
$(x_0,x_1,x_2,x)$ (blue base), $(x_0,x_1,x_3,x)$ (yellow base),  and $(x_0,x_2,x_3,x)$ (red base).
Right: Six admissible simplexes with $x_0$ at an acute angle of their bases:
$(x_0,x_1,x_5,x)$ (blue base), $(x_0,x_2,x_5,x)$ (violet base), $(x_0,x_1,x_6,x)$ (yellow base), 
$(x_0,x_3,x_6,x)$ (green base), $(x_0,x_2,x_4,x)$ (red base),  and $(x_0,x_3,x_4,x)$ (magenta base). 
}
\label{fig:simplexes}
\end{center}
\end{figure}

Let a {\sf Considered} point $x$ be up for a simplex update. 
Figure \ref{fig:simplexes} shows all bases of admissible simplexes  containing a given {\sf Accepted Front} point $x_0$ and lying in a chosen octant.
There are  three admissible simplexes (Figure \ref{fig:simplexes}, Left) such that $x_0$ is at the right angle of their bases, 
and six admissible simplexes (Figure \ref{fig:simplexes}, Right) such that  $x_0$ is at an acute angle. 
Of course, any of these simplexes is admissible provided that all vertices of its base are {\sf Accepted Front}.

A fast search for admissible simplexes in {\tt olim3D} is organized as follows. 
The  $\mathcal{N}_1\cup\mathcal{N}_2$ neighbors of each mesh point are indexed  from $0$ to $17$ as shown in Fig. \ref{fig:neib}.
For each neighbor index $0\le i\le 17$, the set of indexes $0\le j\le 17$
such that $(x_0,i,j)$ make a base of an admissible simplex is listed. Then in order to find a base of an admissible simplex given its two
vertices $x_0$ and its  $\mathcal{N}_1\cup\mathcal{N}_2$ neighbor with index $i$, the routine checks the list of possible third indexes $j$. 
This procedure is organized so that repetitions of simplex bases are avoided.

\subsection{The hierarchical update strategy}
\label{sec:hierarchy}
In this work, we upgrade the hierarchical update strategy introduced in \cite{olim} to a new level of efficiency.
The new version of the hierarchical update is outlined in Algorithm \ref{alg:HU} below.
For comparison, we first sketch the straightforward brute-force approach for
solving minimization problem \eqref{Gaction} on a 3D mesh (see Algorithm \ref{alg:bruteforce}).

\noindent\makebox[\linewidth]{\rule{\textwidth}{0.4pt}}
 \begin{algorithm}[H]
 \label{alg:bruteforce}
\KwInitialization{
Start with all mesh points being {\sf Unknown}. Let $x^{\ast}$ be an asymptotically stable equilibrium  located at a mesh point.
Compute tentative values { of $U$} at the 26 nearest neighbors of $x^{\ast}$ and mark them {\sf Considered}.
}

\KwTheMainBody\\
\While { the boundary of the mesh has not been reached {\bf and} the set of {\sf Considered} points is not empty}{
 {\bf 1:} Change the status of the  {\sf Considered} point $x_0$ with the smallest tentative value of $U$ to {\sf Accepted Front}. \\
{\bf 2:} Set all  {\sf Accepted Front} points in $\mathcal{N}(x_0)$ that no longer have {\sf Considered} nearest neighbors to  {\sf Accepted}.\\
{\bf 3:} \For{ all  {\sf Considered} points $x\in\mathcal{N}_{{\rm far}}^K(x_0)$}{
Form all admissible simplexes $(x_0,x_1,x_2,x)$ where $x_1,x_2\in\mathcal{N}_1(x_0)\cup\mathcal{N}_2(x_0)$ are {\sf Accepted Front}, 
find the solution $\mathsf{Q}_3(x_0,x_1,x_2,x)$ to the constrained minimization problem \eqref{3u}, and set
$U(x) = \min\{U(x),\mathsf{Q}_3(x_0,x_1,x_2,x)\}$.
}
{\bf 4:} \For{ all {\sf Unknown} points $x\in\mathcal{N}(x_0)$}{
{\bf 4a:} Change the status of $x$ to {\sf Considered}. \\
{\bf 4b:} \For{ all {\sf Accepted Front} points $y\in\mathcal{N}_{{\rm far}}^K(x)$}{
Form all admissible simplexes $(y,x_1,x_2,x)$ where $x_1,x_2\in\mathcal{N}_1(y)\cup\mathcal{N}_2(y)$ are {\sf Accepted Front}, 
find the solution $\mathsf{Q}_3(y,x_1,x_2,x)$ to the constrained minimization problem \eqref{3u}, and set
$U(x) = \min\{U(x),\mathsf{Q}_3(y,x_1,x_2,x)\}$.
}
}
}
\caption{A brute-force version of the OLIM for computing the quasipotential in 3D with respect to an
asymptotically stable equilibrium $x^{\ast}.$}
\end{algorithm}
\noindent\makebox[\linewidth]{\rule{\textwidth}{0.4pt}}

Unfortunately, the brute-force approach would lead to excessively large CPU times 
due to a huge number of simplex updates. Most of these updates would give a solution to the
constrained minimization problem \eqref{3u} lying on the boundary of the corresponding simplex base. 
We did not try to implement the
brute-force  version of the OLIM  for finding the quasipotential in 3D. 
However, we did implement it for solving the eikonal equation $\|u\| = s(x)$ where $s(x)> 0$ is the given slowness function. 
Essentially, the brute-force OLIM-based  eikonal solver is obtained from Algorithm \ref{alg:bruteforce} 
by setting $K=1$ and replacing the integrand in  \eqref{Gaction} with merely $s(\psi)$.
Applied to the eikonal equation, our implementation using the hierarchical update strategy reduced  CPU times by a
factor exceeding 5 while  only increasing the numerical errors by less than 1\%. 
The application of the OLIM to the eikonal equation will be reported in \cite{PC}.

The idea of the hierarchical update strategy is to limit  triangle and simplex updates  with arguments
$(x_0,x_1,x)$ and $(x_0,x_1,x_2,x)$ to those  where the point $x_0$ is such that
\begin{equation}
\label{q1min}
\mathsf{Q}_1(x_0,x) = \min\{\mathsf{Q}_1(y,x)~|~y\in \mathcal{N}_{{\rm far}}^K(x),~y~\text{is {\sf Accepted Front}}\}.
\end{equation}
In other words, the triangle update with arguments $(x_0,x_1,x)$  or the simplex update with arguments $(x_0,x_1,x_2,x)$
are attempted only if the point $x_0$ is the minimizer of the one-point update for $x$.
The rationale is that the update function
\begin{equation}
\label{ufun}
f(x_{\lambda}): = Iu(x_{\lambda}) + \left\|b\left(\frac{x + x_{\lambda}}{2}\right)\right\|\|x - x_{\lambda}\| - b\left(\frac{x + x_{\lambda}}{2}\right)\cdot (x - x_{\lambda})
\end{equation}
 continuously depends on the basis of the characteristic $x_{\lambda}$ arriving at $x$.
Here, $Iu(x_{\lambda})$ is the linear interpolation of $U$ at the point $x_{\lambda}$ obtained from the
three {\sf Accepted Front} vertices of the base of the admissible simplex containing $x_{\lambda}$.
In the typical case (though not always), $f(x_{\lambda})$ has a unique local minimizer. 
The hierarchical update is validated by numerical evidence of convergence of the numerical solution
obtained by {\tt olim3D} to the true solution with mesh refinement (see Section \ref{sec:tests}).

One more useful consideration allows us to further reduce the number of unnecessary simplex updates.
Suppose that  the Hessian \eqref{hess3} is positive definite.
Then, if there is an interior point solution to the simplex update \eqref{3u}, 
there will be also interior point solutions to the three triangle
updates whose bases form the boundary of the simplex update in question. 
This leads to the following update strategy.
Let $x_0$ be the minimizer of the one-point update. Then triangle updates are attempted for every {\sf Accepted Front}
point $x_1\in\mathcal{N}_1(x_0)\cup\mathcal{N}_2(x_0)$. Whenever the triangle update $\mathsf{Q}_2(x_0,x_1,x)$
gives an interior point solution, 
a simplex update is attempted for every {\sf Accepted Front} point $x_2$ forming an admissible simplex base $(x_0,x_1,x_2)$.

Algorithm \ref{alg:HU} outlines the hierarchical update strategy implemented in {\tt olim3D}.
We store the minimal one-point update values $\mathsf{Q}_1(x_0,x)$ and the indices of their minimizers $x_0$.

\noindent\makebox[\linewidth]{\rule{\textwidth}{0.4pt}}
 \begin{algorithm}[H]
 \label{alg:HU}
\KwInitialization{
Start with all mesh points being {\sf Unknown}. Let $x^{\ast}$ be an asymptotically stable equilibrium  located at a mesh point.
Compute tentative values of $U$ at the 26 nearest neighbors of $x^{\ast}$ and change their status to {\sf Considered}.
}

\KwTheMainBody\\
\While { the boundary of the mesh has not been reached {\bf and} the set of {\sf Considered} points is not empty}{
 {\bf 1:} Set the  {\sf Considered} point $x_{\rm new}$ with the smallest tentative value of $U$  to {\sf Accepted Front}. \\
{\bf 2:} Change all  {\sf Accepted Front} points in $\mathcal{N}(x_{\rm new})$ that no longer have {\sf Considered} nearest neighbors to  {\sf Accepted}.\\
{\bf 3:} Form the set $\mathcal{AN}(x_{\rm new})$ of {\sf Accepted Front} points $y\in\mathcal{N}_1(x_{\rm new})\cup \mathcal{N}_2(x_{\rm new})$.\\
{\bf 4:} Update all  {\sf Considered} points $x\in\mathcal{N}_{{\rm far}}^K(x_{\rm new})$.\\
{\bf 5:} Change the status of all {\sf Unknown} points $x\in\mathcal{N}(x_{\rm new})$ to {\sf Considered} and update them.
}
\caption{The hierarchical update strategy implemented in {\tt olim3D} for computing the quasipotential with respect to an
asymptotically stable equilibrium $x^{\ast}.$}
\end{algorithm}
\noindent\makebox[\linewidth]{\rule{\textwidth}{0.4pt}}

Steps {\bf 4} and {\bf 5} of Algorithm \ref{alg:HU} are elaborated below (see Algorithms \ref{alg:HU1} and \ref{alg:HU2} respectively).

\noindent\makebox[\linewidth]{\rule{\textwidth}{0.4pt}}
 \begin{algorithm}[H]
 \label{alg:HU1}
{\bf 4:} \For{ all  {\sf Considered} points $x\in\mathcal{N}_{{\rm far}}^K(x_{\rm new})$}{
{\bf 4a:} Compute the one-point update $\mathsf{Q}_1(x_{\rm new},x)$. \\
{\bf 4b:} \If{ $\mathsf{Q}_1(x_{\rm new},x)$ is the minimal one-point update for $x$} {
{\bf 4c:} \For{ all $y\in \mathcal{AN}(x_{\rm new})$ } { 
{\bf 4d:} Do the triangle update  $\mathsf{Q}_2(x_{\rm new},y,x)$.\\
{\bf 4e:} \If{ an interior point solution is found}{
{\bf 4f:} \For{ all $z\in \mathcal{AN}(x_{\rm new})$ such that $(x_{\rm new},y,z)$ forms an admissible simplex base
which has not been attempted}{
{\bf 4g:} Do the simplex update  $\mathsf{Q}_3(x_{\rm new},y,z,x)$.
}
}
}
}
{\bf 4h:} \Else{
{\bf 4i:} Let $x_0$ be the base of the minimal one-point update for $x$.\\
{\bf 4j:} \If{ $x_0\in\mathcal{AN}(x_{\rm new})$ }{
{\bf 4k:} Do the triangle update $\mathsf{Q}_2(x_0,x_{\rm new},x)$.\\
{\bf 4l:} \If{  an interior point solution is found}{
{\bf 4m:} \For{ all $z\in \mathcal{AN}(x_{\rm new})$ such that $(x_0,x_{\rm new},z)$ forms an admissible 
simplex base}{
{\bf 4n:} Do the simplex update $\mathsf{Q}_3(x_0,x_{\rm new},z,x)$.
}
}
}
}
}
\caption{Substeps of Step {\bf 4} of Algorithm \ref{alg:HU}.}
\end{algorithm}
\noindent\makebox[\linewidth]{\rule{\textwidth}{0.4pt}}

\noindent\makebox[\linewidth]{\rule{\textwidth}{0.4pt}}
 \begin{algorithm}[H]
 \label{alg:HU2}
{\bf 5:} \For{ all {\sf Unknown} points $x\in\mathcal{N}(x_{\rm new})$}{
{\bf 5a:} Change the status of $x$ to {\sf Considered}. \\
{\bf 5b:} \For{ all {\sf Accepted Front} points $y\in\mathcal{N}_{{\rm far}}^K(x)$}{
{\bf 5c:} Do the one-point update $\mathsf{Q}_1(y,x)$;
}
{\bf 5d:} Let $x_0$ be the base of the minimal one-point update;\\
{\bf 5e:} \For{ all $y\in \mathcal{AN}(x_0)$}{
{\bf 5f:} Do the triangle update $\mathsf{Q}_2(x_0,y,x)$;\\
{\bf 5g:} \If{ an interior point solution is found}{
{\bf 5h:} \For{ all $z\in \mathcal{AN}(x_0)$ such that $(x_0,y,z)$ forms an admissible 
simplex base which has not been attempted}{
{\bf 5i:} Do the simplex update $\mathsf{Q}_3(x_0,y,z,x)$.
}
}
}
}
\caption{Substeps of Step {\bf 5} of Algorithm \ref{alg:HU}.}
\end{algorithm}
\noindent\makebox[\linewidth]{\rule{\textwidth}{0.4pt}}

\subsection{Skipping simplex update by using the KKT conditions}
\label{sec:KKT}
During a triangle update, the first step is to check whether the derivative of the function to be minimized in  \eqref{2u}
has opposite signs at the endpoints $\lambda = 0$ and $\lambda = 1$. 
This criterion rejects 50\% to 70\% of triangle update attempts.
Since each simplex update is called only if
an interior point solution is found for a triangle update, 
and since this solution is used as the initial iterate, 
the Karush-Kuhn-Tucker (KKT) theorem (see Chapter 12 of \cite{nocedal})
gives a simple criterion to reject the simplex update right away. 
As mentioned in Section \ref{sec:simplex}, 50\% to 70\% of simplex updates are rejected as a result.
Let us elaborate. 
Suppose a triangle update gave the solution $\mathsf{Q}_2(x_0,x_1,x)$ with the corresponding minimizer $\lambda^{\ast}\in(0,1)$,
and now the simplex update with arguments $(x_0,x_1,x_2,x)$ is to be done.
The corresponding constrained minimization problem  in canonical form is:
\begin{align}
\min f(\lambda) &~~{\rm where}~~f(\lambda) = U_{\lambda} +\|b_{m\lambda}\||x-x{_\lambda}\| - b_{m\lambda}\cdot(x-x_{\lambda}) \label{fmin1}\\
&\text{subject to}\notag\\
&\lambda_1\ge 0,\label{i1}\\
&\lambda_2\ge 0,\label{i2}\\
&1-\lambda_1-\lambda_2\ge 0,\label{i3}
\end{align}
where $\lambda =( \lambda_1,\lambda_2)$, 
$U_{\lambda}$, $x_{\lambda}$ and $b_{m\lambda}$ are the linear interpolants defined in Eq. \eqref{3u}.
The Lagrangian for Eq.  \eqref{fmin1} is given by
\begin{equation}
\label{Lag}
L(\lambda,\mu) = f(\lambda) - \mu_1\lambda_1 -\mu_2\lambda_2-\mu_3(1-\lambda_1-\lambda_2).
\end{equation}
The KKT optimality conditions applied to  \eqref{fmin1}--\eqref{i3} state that if $\lambda$ is a local solution of Eqs. \eqref{fmin1}--\eqref{i3}
then there exist Lagrange multipliers $\mu_1$, $\mu_2$, and $\mu_3$ such that the following conditions are satisfied:
\begin{align}
\nabla _{\lambda}L(\lambda,\mu) & = \nabla f -\mu_1\left[\begin{array}{c}1\\0\end{array}\right]  -\mu_2\left[\begin{array}{c}0\\1\end{array}\right] 
-\mu_3\left[\begin{array}{c}-1\\-1\end{array}\right]  = \left[\begin{array}{c}0\\0\end{array}\right] \label{kkt1}\\
\mu_1&\ge 0,\quad \mu_2\ge 0,\quad \mu_3\ge 0,\label{kkt2}\\
\lambda_1&\ge0,\quad\lambda_2\ge 0,\quad 1-\lambda_1-\lambda_2\ge 0,\label{kkt3}\\
\lambda_1\mu_1& = 0,\quad \lambda_2\mu_2 = 0,\quad (1-\lambda_1-\lambda_2)\mu_3 = 0.\label{kkt4}
\end{align}
Let us check whether the initial guess $\lambda = (\lambda^{\ast},0)$ where $ \lambda^{\ast}\in(0,1)$
 is a solution 
of \eqref{fmin1}--\eqref{i3}, i.e., whether we can find $\mu_1$, $\mu_2$ and $\mu_3$ such that
the KKT optimality conditions \eqref{kkt1}--\eqref{kkt4} are satisfied.
Condition \eqref{kkt4} forces $\mu_1$ and $\mu_3$ to be zero.
The first component of $\nabla f$ must be zero
at $\lambda = (\lambda^{\ast},0)$ since $\lambda^{\ast}\in(0,1)$ is a solution of \eqref{2u}. 
Then $\mu_2$ found from  \eqref{kkt1} satisfies
\begin{equation}
\label{kkt5}
\frac{\partial f}{\partial\lambda_2} - \mu_2 = 0.
\end{equation}
If $\mu_2 = \tfrac{\partial f}{\partial\lambda_2}\ge 0$, then the KKT conditions \eqref{kkt1}--\eqref{kkt4}
hold and hence $\lambda = (\lambda^{\ast},0)$ is a local solution. In this case, we skip the simplex update.
Otherwise, we use Newton's method to solve the constrained minimization problem \eqref{fmin1}--\eqref{i3}.

\subsection{Remarks about local factoring}
\label{sec:locfac}
It was shown in a series of works
that the accuracy of numerical solutions to the eikonal equation initialized near point sources
 can be significantly enhanced by factoring.
Originally, a multiplicative factoring  
was introduced in \cite{FLZ}  for the fast sweeping eikonal solver 
as a tool to decrease error constants and increase order of convergence in the case of point sources. 
An additive factoring serving the same purpose was proposed in \cite{LQ}. 
The idea of factoring was adapted for the fast marching method (FMM) in \cite{QV}. 
Furthermore, the fact that the FMM propagates the solution throughout the domain from smaller values to  larger values without 
iteration allows for factoring locally; i.e., the eikonal equation only needs to be factored around point sources and rarefaction fans. 

In this work, we combine local factoring with the OLIM quasipotential solver.
Our results show that the effect of local factoring depends 
on the ratio of magnitudes of the rotational and potential components of the vector field $b$.
\begin{itemize}
\item Local factoring tends to reduce numerical errors by about 15\% to 30\% in cases where the rotational component of $b$ is
of the same order of magnitude as the potential one.
\item Local factoring tends to have no effect and may even increase numerical errors 
in cases where the rotational component is significantly  larger than the potential one or
$b$ is not differentiable at the asymptotically stable equilibrium (see Example 6 in Section \ref{sec:tao}).
\end{itemize}
Since local factoring enhances accuracy in some cases, and since mesh refinement 
is limited in 3D, it is worth trying it if it is known that the rotational component is not large in comparison with the potential component.
Below, we describe the implementation of local factoring for the OLIMs and explain why its effect depends on the relationship between the
rotational and potential components.

The implementation of local factoring for the OLIMs is simple. 
Let $x^{\ast}$ be an asymptotically stable equilibrium. Let $J$ be the Jacobian matrix of the vector field $b(x)$ evaluated at $x^{\ast}$.
Then, one can find the quasipotential matrix $Q$ for $J$ as explained in Section \ref{sec:init}. 
Thus, in the neighborhood of $x^{\ast}$, the quasipotential is approximated by 
\begin{equation}
\label{ustar}
U^{\ast}(x): = (x-x^{\ast})^\top Q(x-x^{\ast})
\end{equation}
 We choose a radius $R_{{\rm fac}}$ and decompose the quasipotential so that
 \begin{equation}
\label{efac1}
U(x) = U^{\ast}(x) + u(x)\quad\text{for all}~~x~~\text{such that}~~\|x - x^{\ast}\|\le R_{\rm fac}.
\end{equation} 
Here $U^{\ast}(x)$ is given by  \eqref{ustar}, and $u(x)$, the correction to the linear approximation, is the new unknown function.
Ansatz \eqref{efac1} changes the term $U_{\lambda}$ in Eqs. \eqref{2u} and \eqref{3u} to
\begin{equation}
\label{efac2}
U_{\lambda} \longrightarrow u_{\lambda} + U^{\ast}(x_{\lambda}),
\end{equation}
where  
$$
u_{\lambda} : = u_0 + \delta u^\top \lambda~~{\rm with}~~ u_i: = U(x_i) -U^{\ast}(x_i).
$$ 
For the triangle update, $\delta u = u_1 - u_0$, while
for the simplex update, $\delta u^\top = [u_1 - u_0,u_2 - u_0]$.
The modifications to the gradient \eqref{grad3} and the Hessian \eqref{hess3}
are readily found from Eqs. \eqref{ustar} -- \eqref{efac2}. 

The effect of local factoring in the case of the quasipotential solver is less dramatic than that in the case of the eikonal equation \cite{QV,PC}.
Numerical errors committed by the OLIMs in computing the quasipotential 
near asymptotically stable equilibria are mainly due to $(i)$ the  approximation of $U$ in the 
triangle and simplex bases using linear interpolation, and $(ii)$ the use of line segments to approximate segments of the MAPs.
Issue $(i)$ is largely resolved by local factoring, while issue $(ii)$ is not. 
Furthermore, unlike the solution to the eikonal equation, the quasipotential has no singularity 
near asymptotically stable equilibria since the vector field $b$ vanishes at these points.
Hence, the quasipotential is $O(h^2)$ and absolute numerical errors are at most $O(h^2)$
near the asymptotically stable equilibria.

\begin{figure}[htbp]
\begin{center}
\centerline{
(a)\includegraphics[height = 0.35\textwidth]{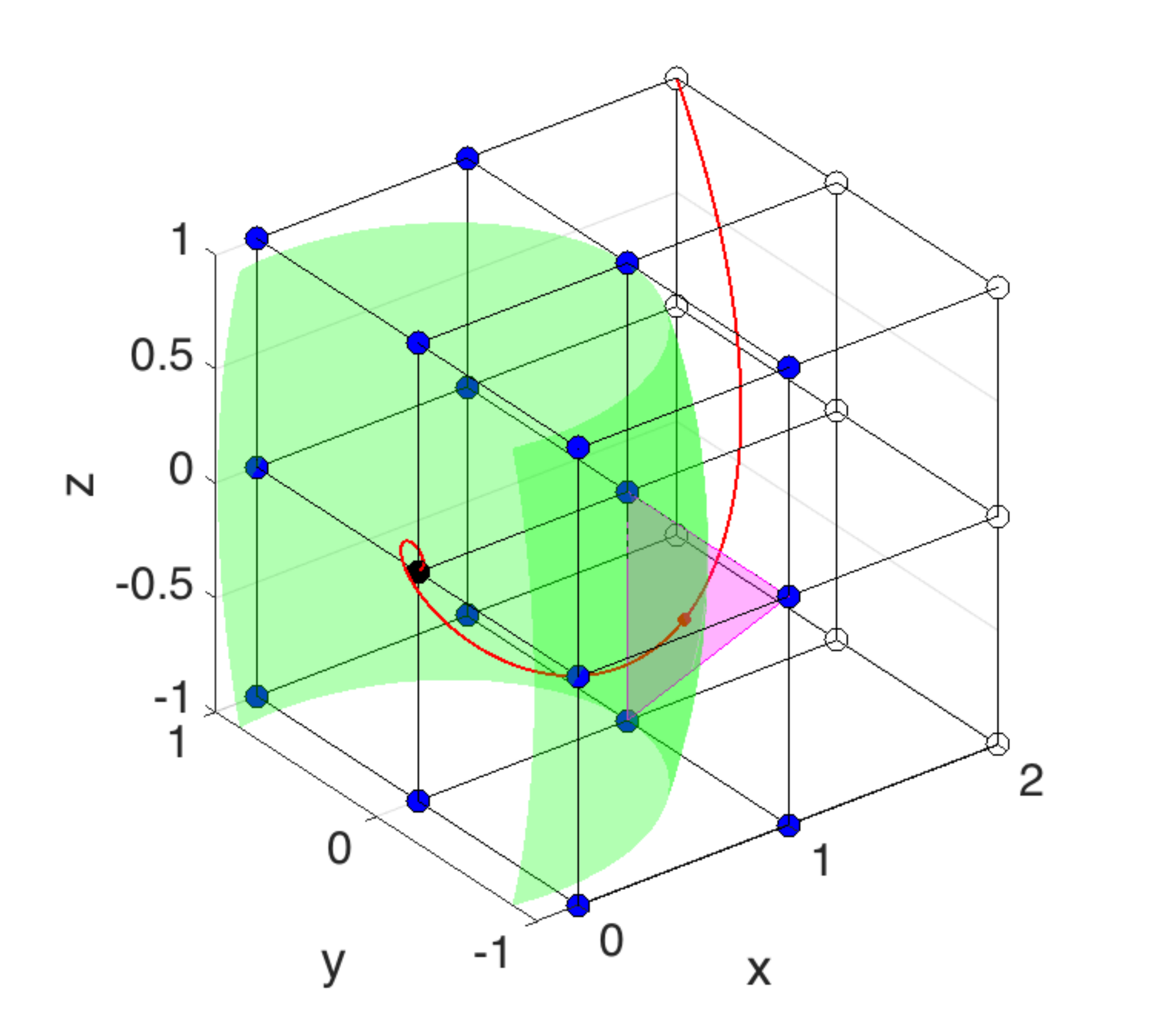}
(b)\includegraphics[height = 0.35\textwidth]{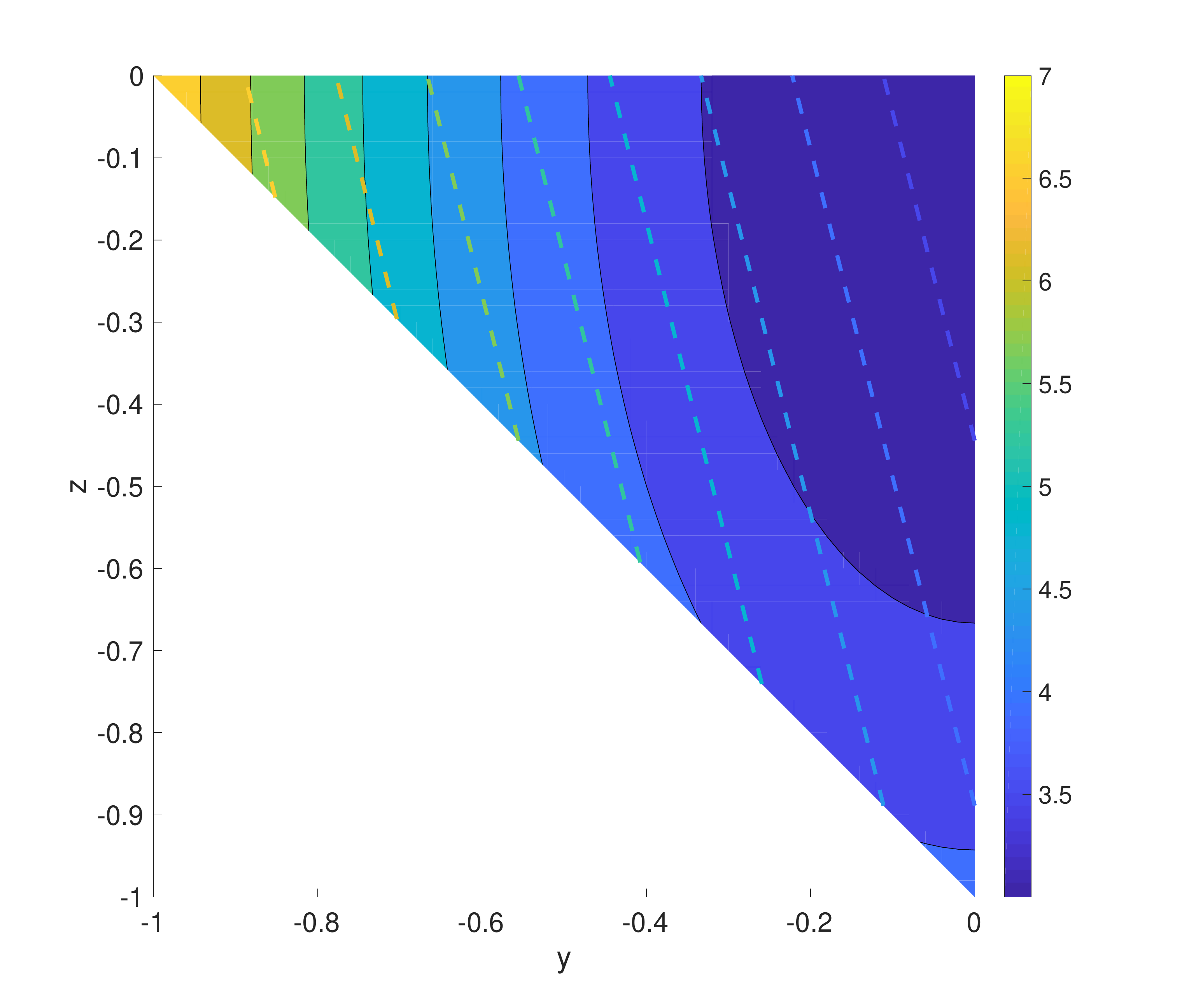}
}
\caption{An illustration for local factoring. 
(a): The red curve is the MAP arriving at the point $(2,1,1)$. 
The red dot is the point of intersection of the MAP with the {\sf Accepted Front} which lies in the triangle $\Delta$ (the magenta patch).
The green surface is the level set of the quasipotential passing through the point $(1,0,-1)$.
(b): The contour plot of the true quasipotential in the triangle $\Delta$ superimposed with the dashed lines representing the contour plot of the linear interpolant
of the true quasipotential in $\Delta$.}
\label{fig:locfac}
\end{center}
\end{figure}

Let us illustrate local factoring on the following example.
Consider the linear vector field 
\begin{equation}
\label{bfield_a}
b(x,y,z) = \left[\begin{array}{rrr} -3&-4&-1\\3&-4&-1\\3&4&-1\end{array}\right]\left[\begin{array}{c}x\\y\\z\end{array}\right].
\end{equation}
The quasipotential matrix for this field is diagonal: $Q = {\rm diag}\{3,4,1\}$.
Let us choose the mesh step size $h=1$.
Suppose the 26 nearest neighbors surrounding the origin, a stable spiral point, are initialized using the exact solution.
The MAP arriving at the mesh point $(2,1,1)$ (the red curve in Fig. \ref{fig:locfac}(a)) intersects the {\sf Accepted Front} as a point
lying in the triangle  $\Delta: = \{(1,0,0),(1,-1,0),(1,0,-1)\}$ which forms a base of an admissible simplex.
A contour plot of $U^{\ast}$ in $\Delta$ is shown in Fig. \ref{fig:locfac}(b). For comparison, a contour plot with dashed lines of the linear interpolant of $U^{\ast}$
is superimposed. The difference between these two is notable. Local factoring largely eliminates this source of error.


\section{Numerical tests}
\label{sec:tests}
We have tested {\tt olim3D} on a collection of examples that includes: 
\begin{itemize}
\item
Three linear examples (Examples 1--3). 

\item Two nonlinear examples (Examples 4--5) 
where the quasipotential is available analytically.
\item
Two of Tao's examples \cite{tao} (Examples 6--7) with two point attractors and hyperbolic periodic orbits 
serving as transition states between them.
In these examples, the quasipotential 
can be found at the transition states using certain orthogonal decompositions (Section \ref{sec:tao}).
\item
A genetic switch model \cite{lv} with two asymptotically stable equilibria separated by a saddle.
The quasipotential is not known analytically, but can be computed along the MAPs (Section \ref{sec:gene}).
\end{itemize}
The results of a large series of measurements performed in Examples 1--5 are presented in Section \ref{sec:meas}.
Effects of local factoring will be exposed in Section \ref{sec:ELF}.

\subsection{Measurements and least squares fits}
\label{sec:meas}
In this Section, we present the results of our study of the dependence of numerical errors and the CPU times  
on the mesh size and the update factor.
The  data are obtained  for the version of {\tt olim3D} without local factoring on a 2017 iMac desktop\footnotemark[1].
\footnotetext[1]{To be precise, here are the specifications. Processor: 4.2GHz Quad- core Intel Core i7, Turbo Boost up to 4.5GHz.
Memory:  64GB 2400MHz DDR4 SODIMM SDRAM - 4$\times$16GB}
The computations have been performed on meshes of sizes $N^3$ where $N = 2^p + 1$, $p = 5,6,7,8,9$, i.e., $N = 33$, 65, 129, 257, and 513.
The values of the update factor $K$  range from 1 to 30.

We will often mention the ratio of the magnitudes of the rotational and potential components of the vector field $b(x)$.
So, we introduce a notation for it:
\begin{equation}
\label{xi}
\xi(x): = \frac{\|l(x)\|}{\|\tfrac{1}{2}\nabla U(x)\|}.
\end{equation}
The first three examples are linear, and the computational domain is the cube $[-1,1]^3$ in all of them.

\textbf{Example 1.}
We consider a linear SDE with $b(x)$ is given by  \eqref{bfield_a}. The exact quasipotential is the quadratic form $U(x) = x^\top Qx$
with $Q = {\rm diag}\{3,4,1\}$. The eigenvectors of $Q$ are aligned with the coordinate axes. 
The ratio $\xi(x)$ varies from 0 to $\sqrt{3}$.


Examples 2 and 3 feature linear systems where the eigenvectors of their quasipotential matrix $Q$ are not aligned with the coordinate axes. 
These examples are constructed so that the ratio $\xi$ can be chosen as desired. 
We pick a matrix $J$ of the form
\begin{equation}
\label{Jmatr}
J= \left[ {\begin{array}{ccc}
   -1 & 0 & 0 \\
   0 & -\tfrac{1}{2} & -\rho \\
   0 & \rho & -\tfrac{1}{2} \\
  \end{array} } \right].
\end{equation}
The corresponding quasipotential matrix and the rotational component are
\begin{equation}
\label{rot23}
Q = \left[\begin{array}{ccc}1&0&0\\0&\tfrac{1}{2}&0\\0&0&\tfrac{1}{2}\end{array}\right]\quad{\rm and}\quad 
l(x) = \left[\begin{array}{c}0\\-\rho x_3\\\rho x_2\end{array}\right].
\end{equation}
The ratio $\xi$ varies from 0 to $2\rho$.
Next, we define a rotation matrix $R$
that rotates  by angle  $\tfrac{\pi}{5}$ about the $x_3$-axis, then by the angle $\tfrac{\pi}{8}$ about $x_2$-axis, 
and then by the angle $\tfrac{2\pi}{3}$ about the $x_1$-axis and  set up the linear SDE
\begin{equation}
\label{linsde23}
dx = R^{^\top} J Rx + \sqrt{\epsilon}dw.
\end{equation}
It is easy to check by plugging $R^{\top}JR$ and $R^{\top}QR$ into  \eqref{sylvester}
that $R^{\top}QR$ is the exact quasipotential matrix for SDE \eqref{linsde23}. 
The rotational component is $l(x) = R^{\top}(J + Q)Rx$. The ratio $\xi$ varies from 0 to $2\rho$ as before.

\textbf{Example 2.}
We pick $\rho = \tfrac{\sqrt{3}}{2}$ in  \eqref{Jmatr}. Hence, $\xi(x)$ ranges from 0 to $\sqrt{3}$.

\textbf{Example 3.}
We pick $\rho = 5$ in  \eqref{Jmatr} which causes $\xi(x)$ to vary from 0 to $10$.

Examples 4 and 5 are constructed from a double-well potential so that
the ratio $\xi$ can be easily prescribed:
\begin{equation}
\label{bfield45}
b(x) = \left[\begin{array}{c} 
- 2(x_1^3 -x_1) - \rho(x_2 + x_3)\\
-x_2 + 2\rho(x_1^3 - x_1)\\
-x_3 + 2\rho(x_1^3 -x_1)\end{array}\right].
\end{equation}
The potential and the rotational components are, respectively,
\begin{equation}
\label{potrot}
- \frac{\nabla U(x)}{2} = \left[\begin{array}{c} 
- 2(x_1^3 -x_1) \\
-x_2 \\
-x_3\end{array}\right]~~{\rm and}~~
l(x)  = \rho\left[\begin{array}{r} 
 - (x_2 + x_3)\\
2(x_1^3 - x_1)\\
2(x_1^3 -x_1)\end{array}\right].
\end{equation}
The ratio $\xi(x)$ equals $\rho$ everywhere except for the equilibria at $(-1,0,0)$, $(0,0,0)$ and $(1,0,0)$ where $\xi(x)$  is undefined.
The equilibria $(\pm 1,0,0)$ are asymptotically stable, while the origin is a Morse index 1 saddle. 
We compute the quasipotential with respect to the attractor $(-1,0,0)$. The computational domain is the cube 
$[-2,0]\times[-1,1]^2$. 
The exact quasipotential within  the level set passing through the origin is given by
\begin{equation}
\label{Qpot45}
U(x) = x_1^4 - 2x_1^2 +x_2^2 + x_3^2 + 1.
\end{equation}

{\bf Example 4.}
The vector field $b(x)$ is given by  \eqref{bfield45} with $\rho = 1$.

{\bf Example 5.}
The vector field $b(x)$ is given by  \eqref{bfield45} with $\rho = 10$.

\begin{figure}[htbp]
\begin{center}
(a)\includegraphics[width = 0.3\textwidth]{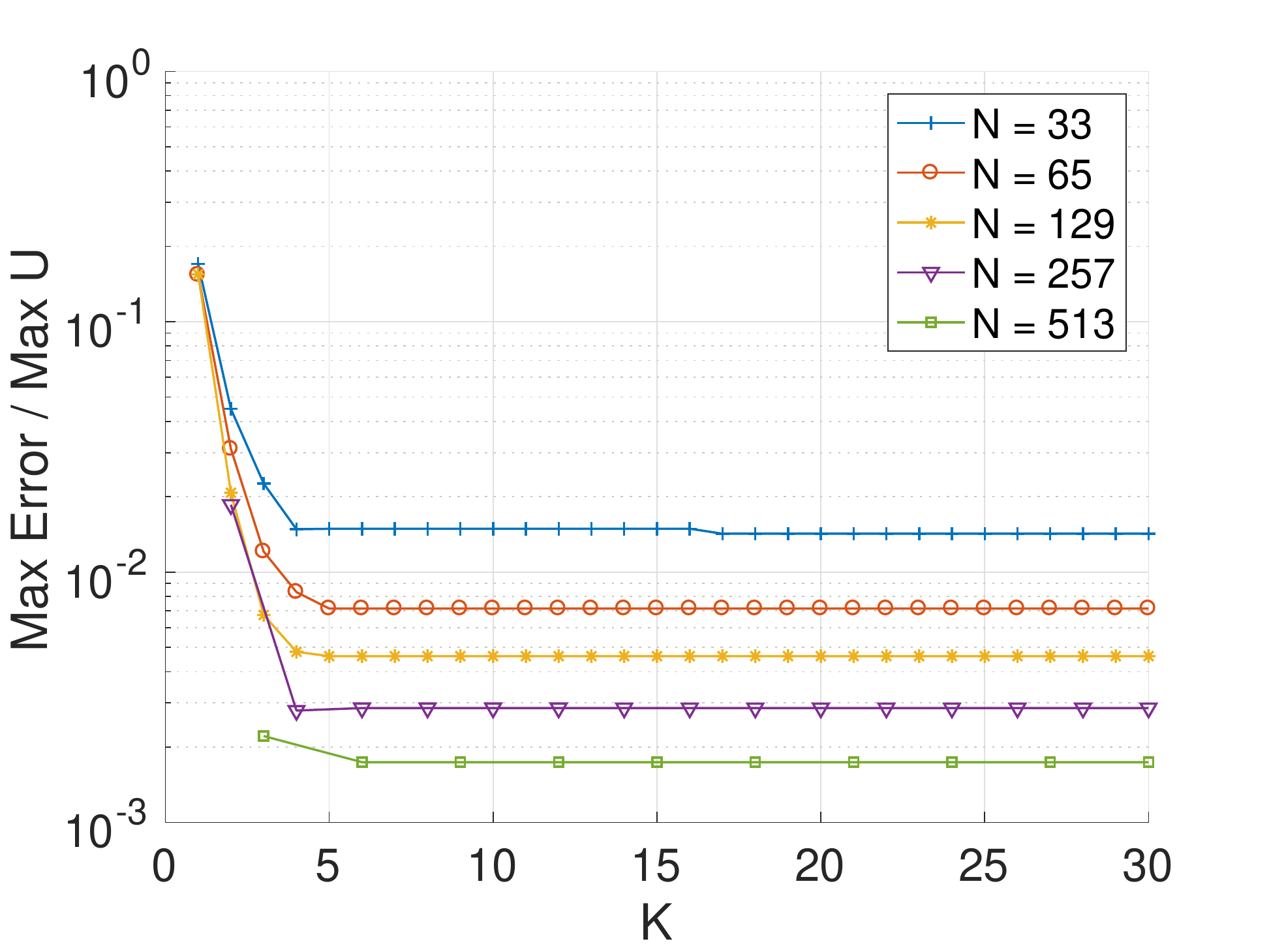}
(b)\includegraphics[width = 0.3\textwidth]{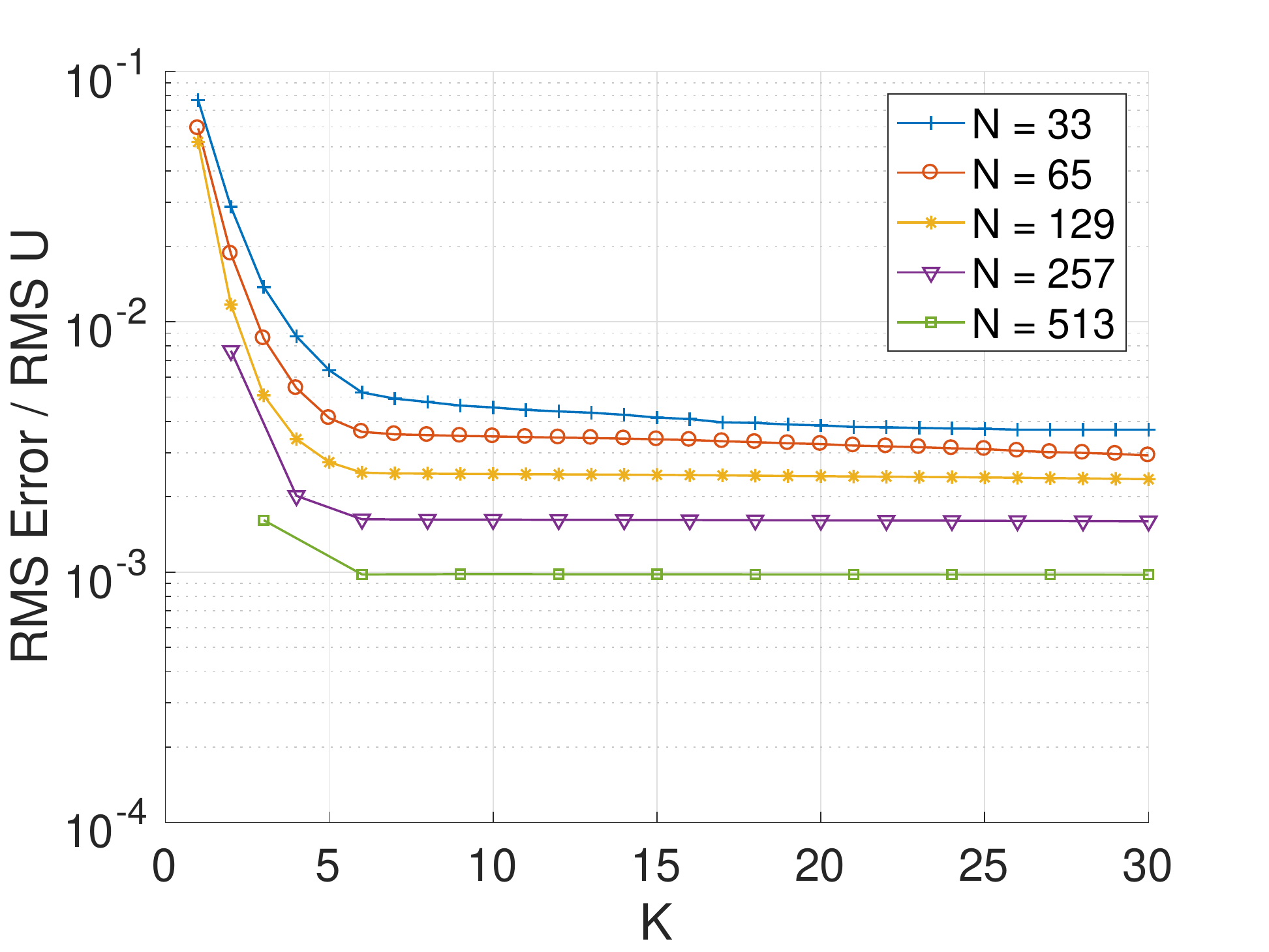}\\
(c)\includegraphics[width = 0.3\textwidth]{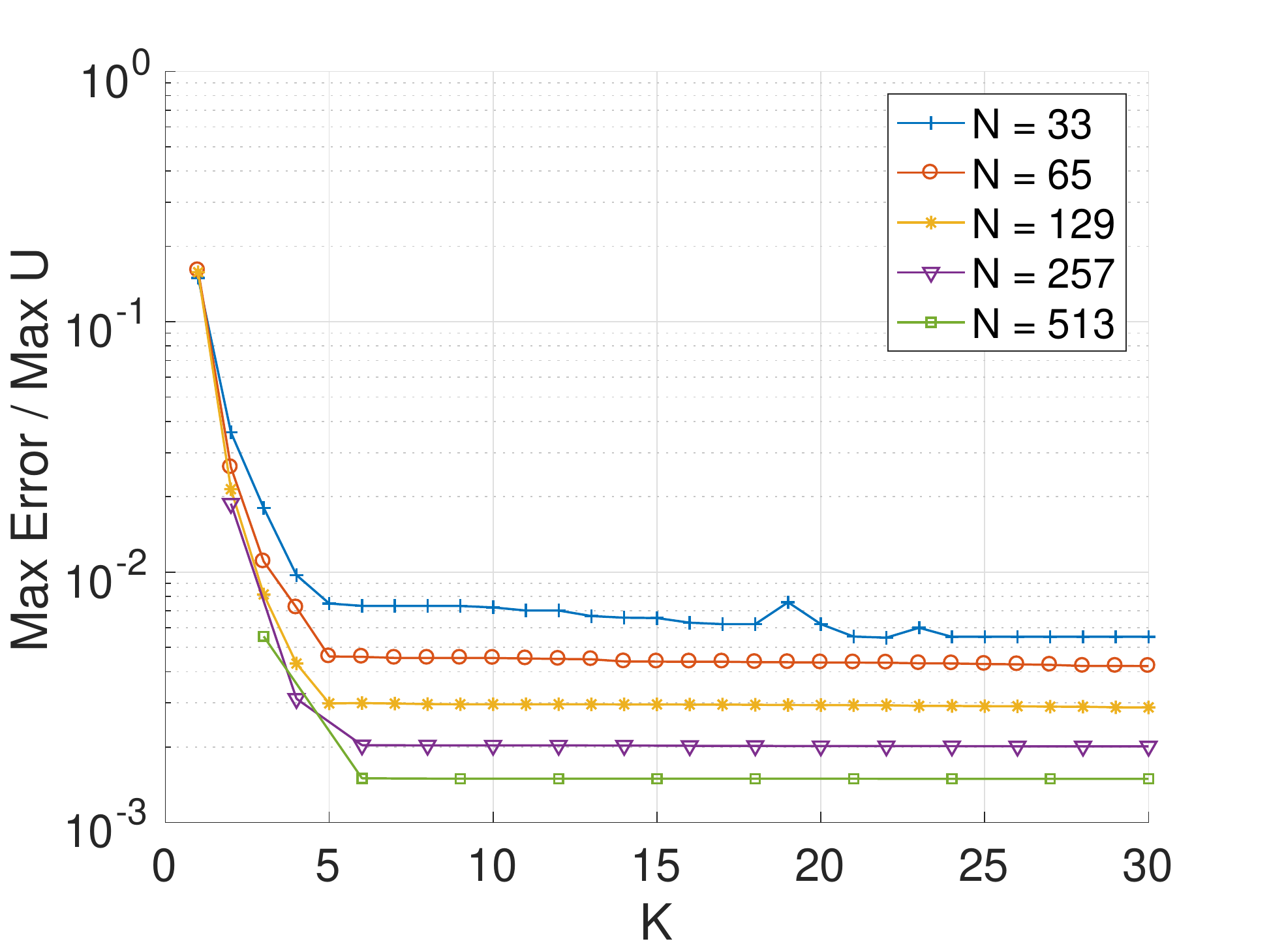}
(d)\includegraphics[width = 0.3\textwidth]{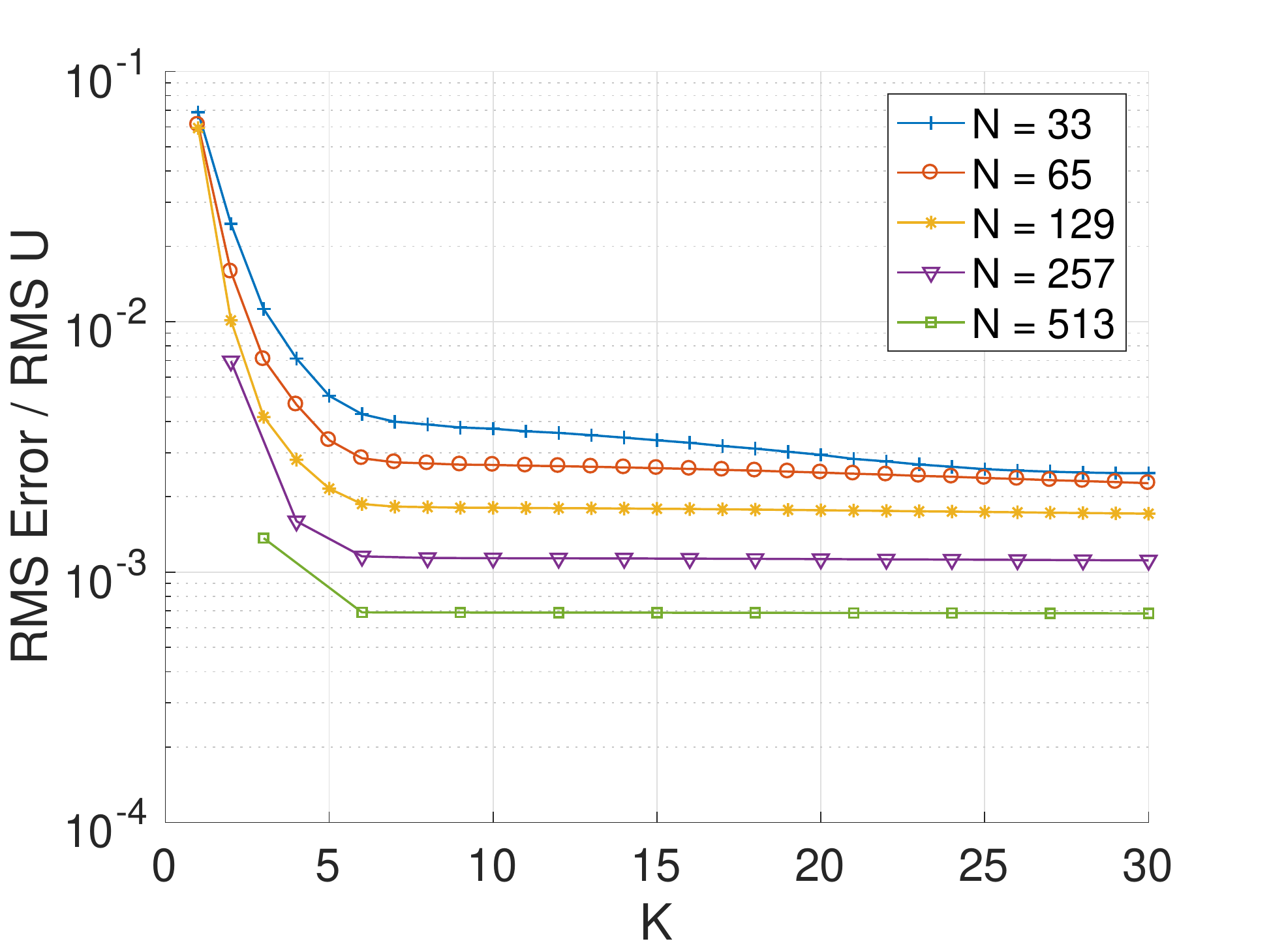}\\
(e)\includegraphics[width = 0.3\textwidth]{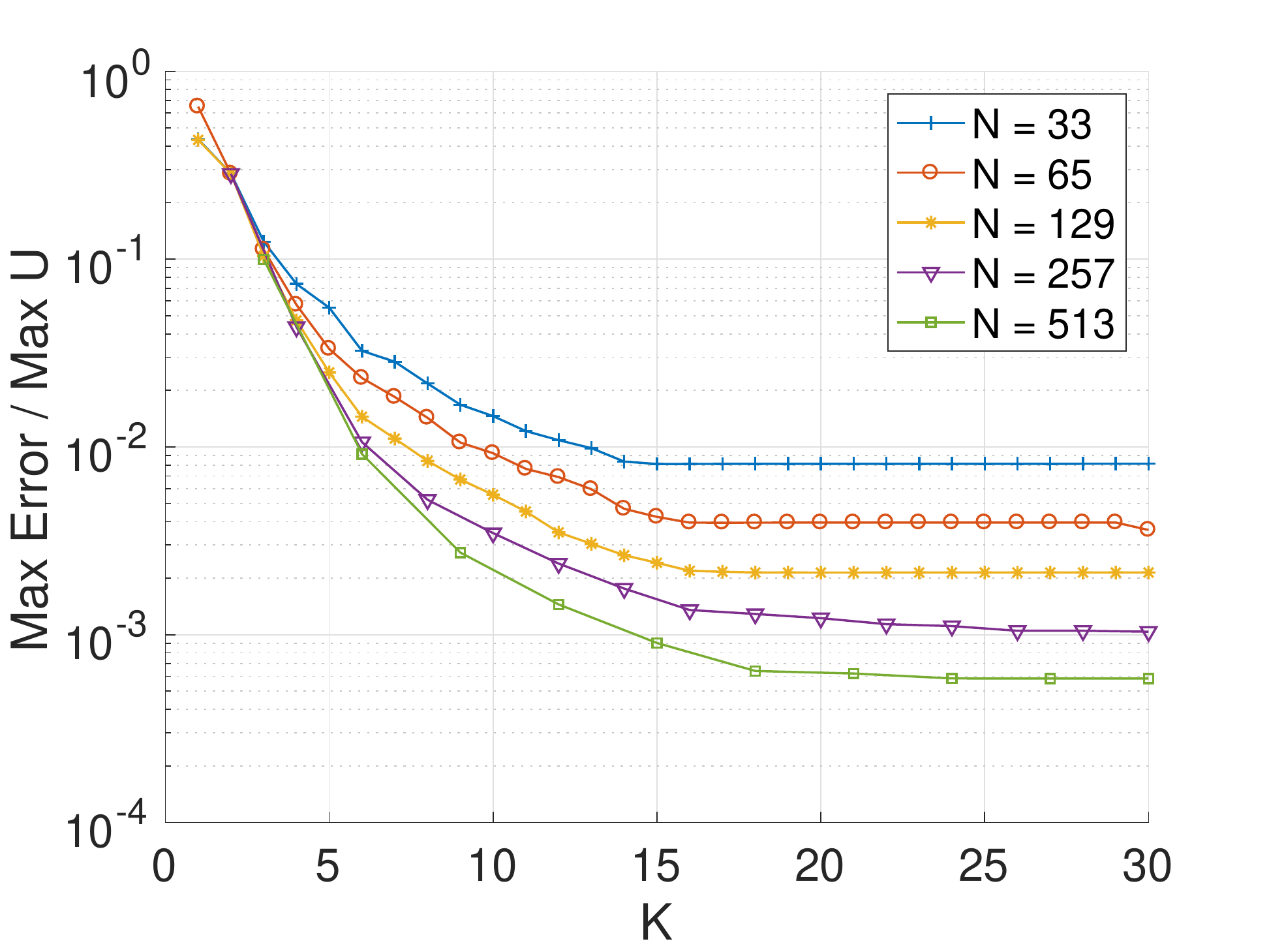}
(f)\includegraphics[width = 0.3\textwidth]{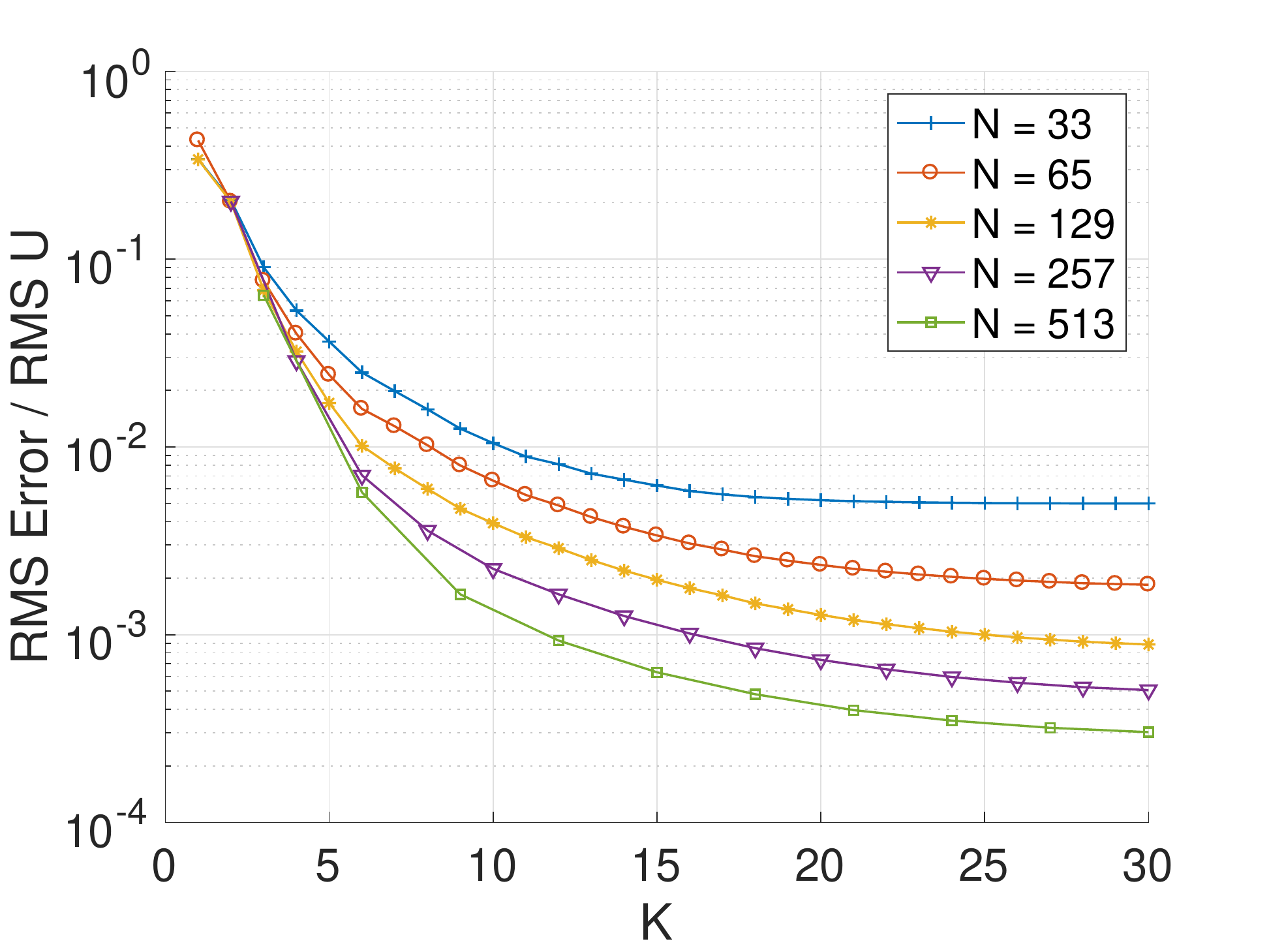}\\
(g)\includegraphics[width = 0.3\textwidth]{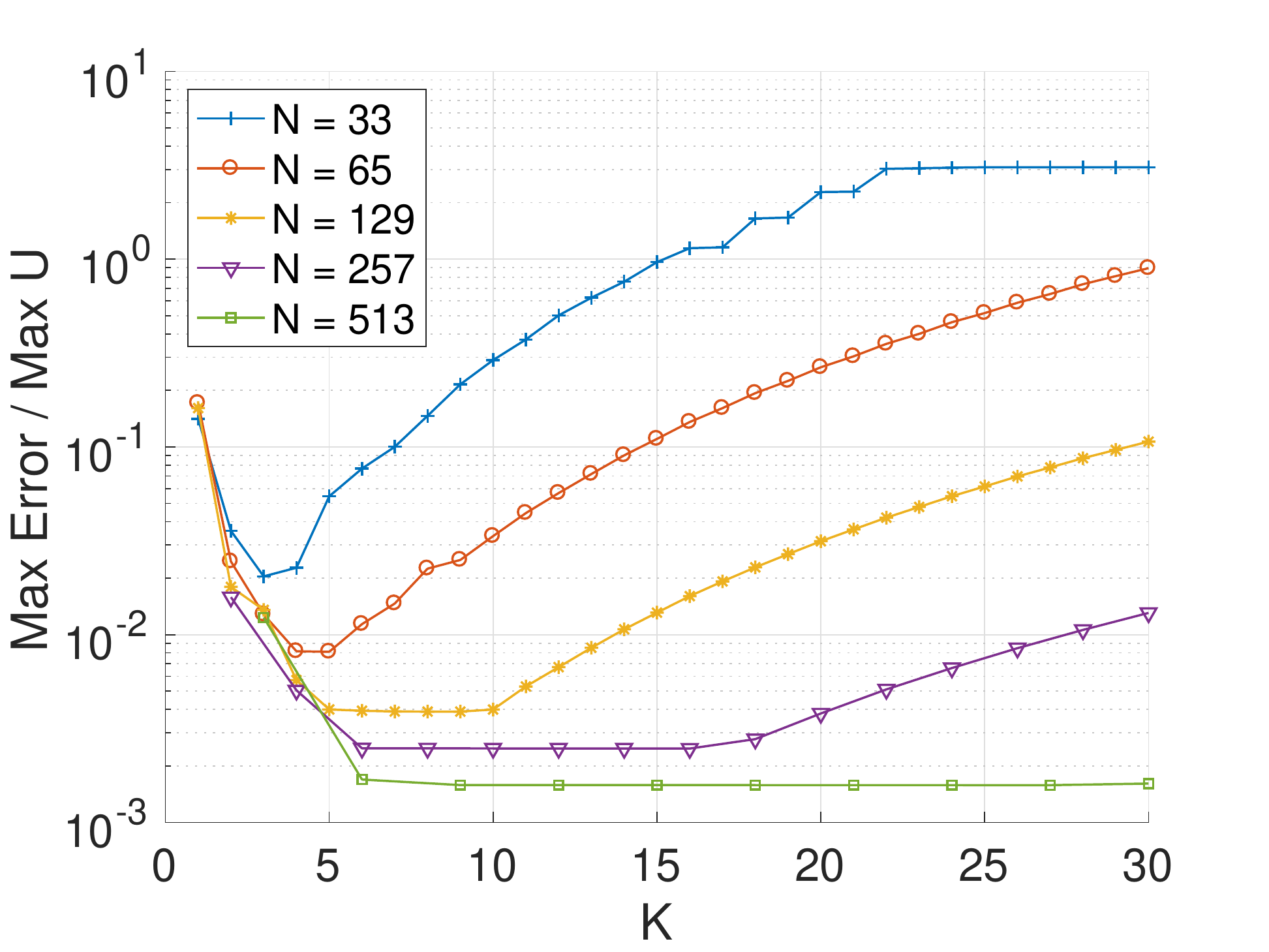}
(h)\includegraphics[width = 0.3\textwidth]{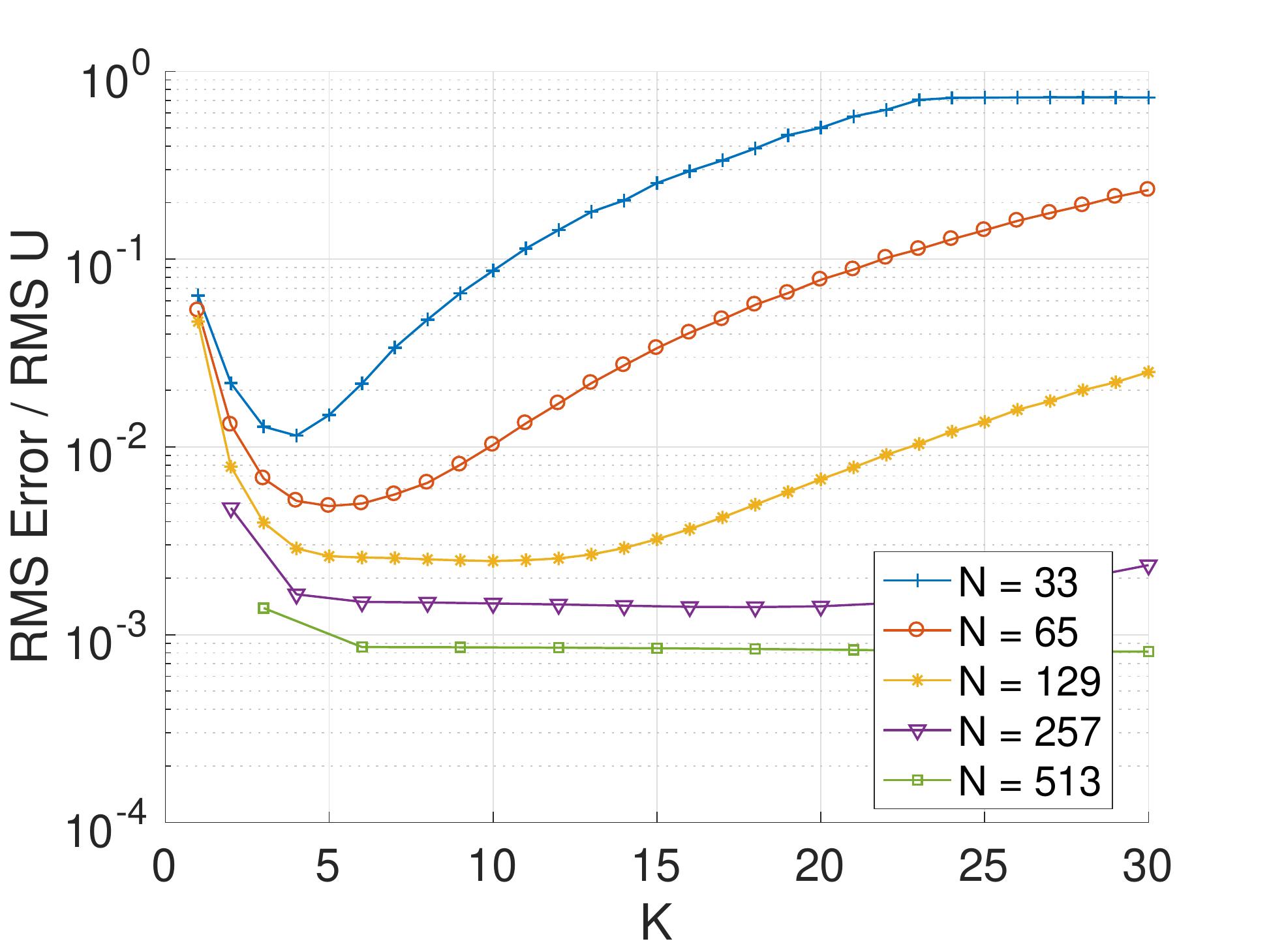}\\
(i)\includegraphics[width = 0.3\textwidth]{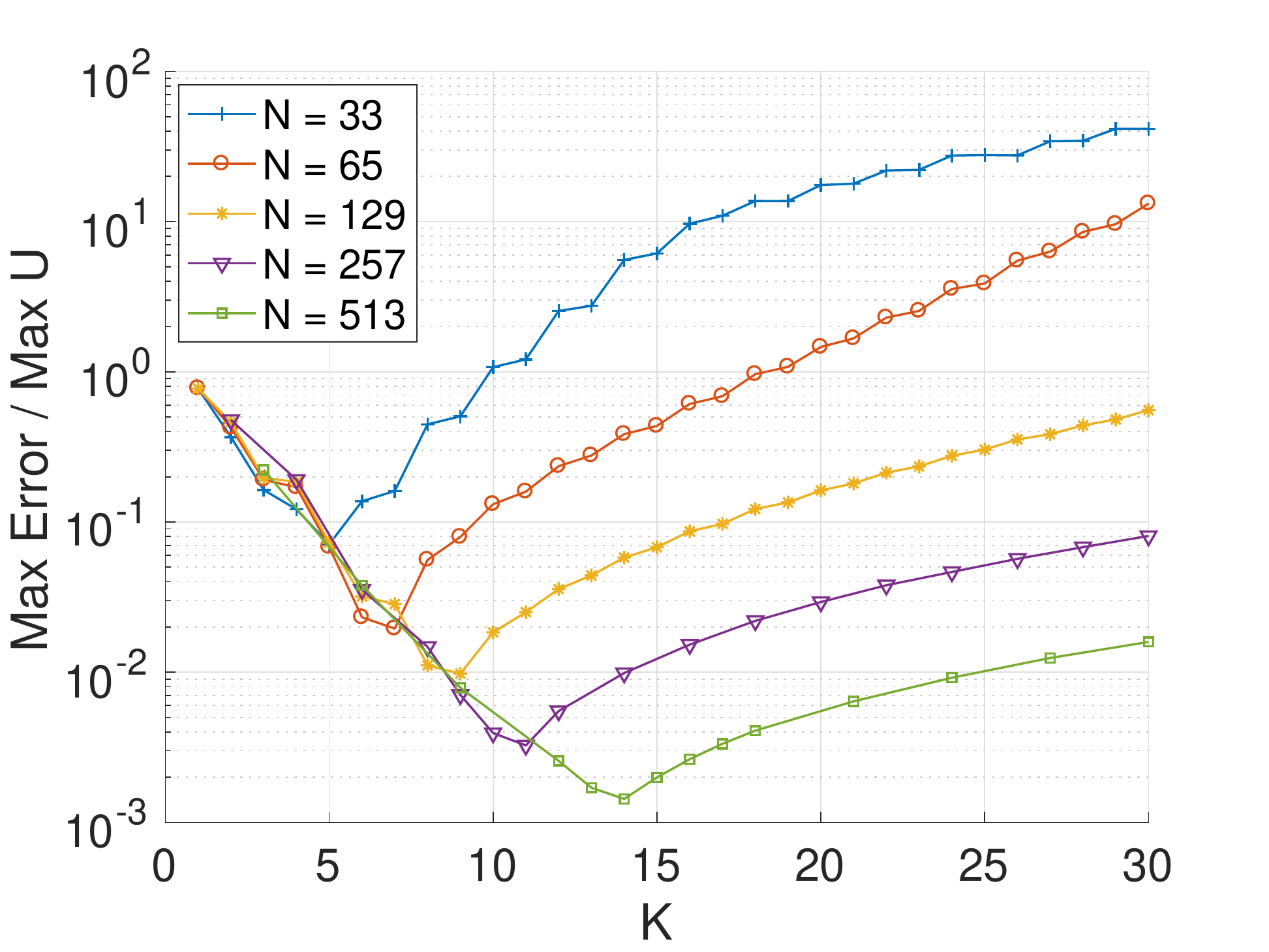}
(j)\includegraphics[width = 0.3\textwidth]{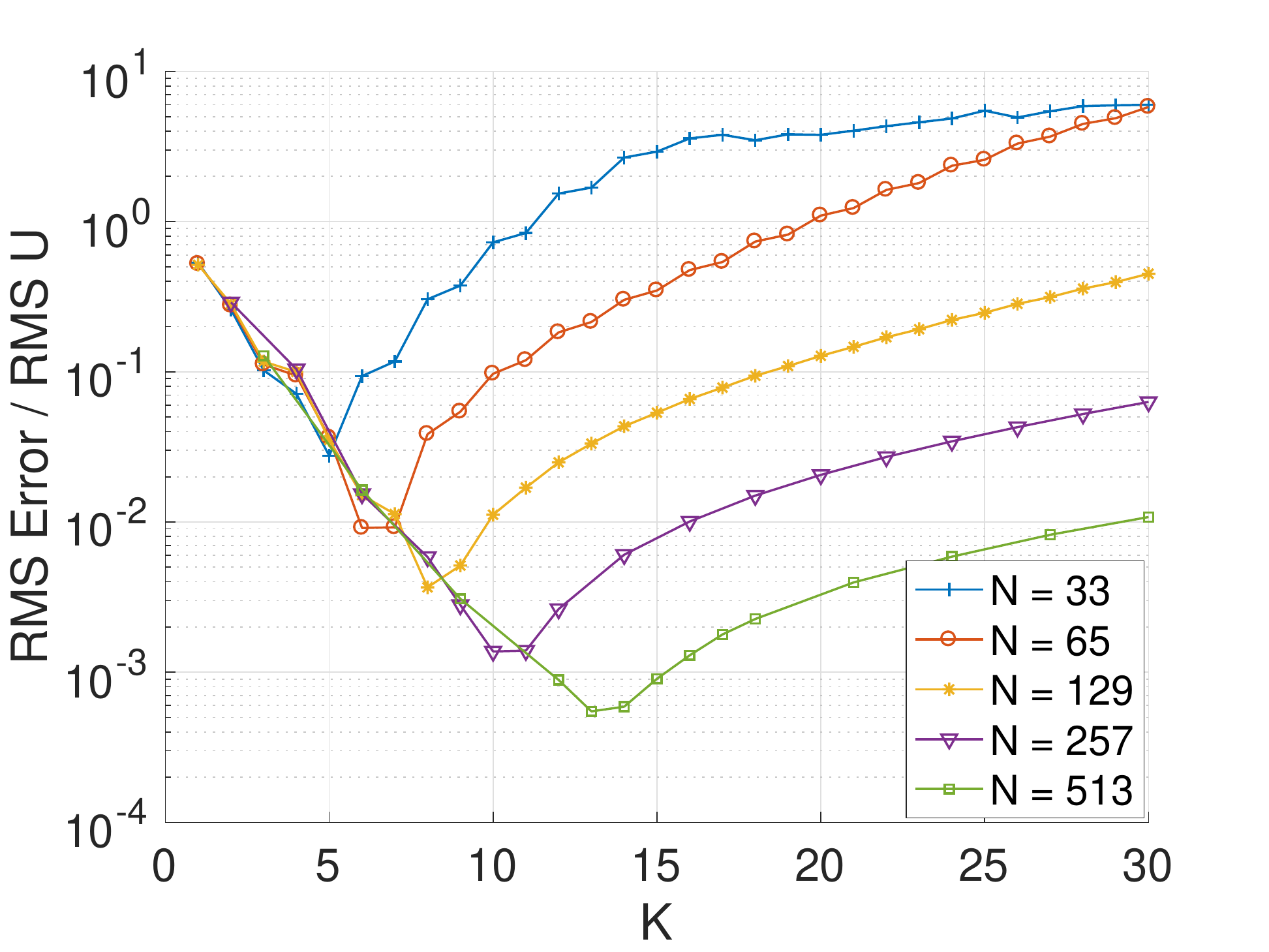}
\caption{(a--b): Example 1. (c--d): Example 2. (e--f): Example 3. (g--h): Example 4. (i--j): Example 5.
(a, c, e, g, i): The normalized maximum absolute error versus the update parameter $K$.
(b, d, f, h, j): The normalized RMS error versus $K$.
}
\label{fig1:ex}
\end{center}
\end{figure}

The graphs of the normalized maximal absolute error and the normalized RMS error versus the update factor 
$K$ for mesh sizes $N^3 = (2^p+1)^3$, $p = 5,6,7,8,9$, for Examples 1--5 are displayed in Fig. \ref{fig1:ex}.
A careful  choice of $K$ is important for achieving optimal accuracy. 
These graphs allow us to make recommendations regarding a reasonable choice of $K$.
First, we make a few observations.
\begin{itemize}
\item 
In Examples 1--3, where the SDEs are linear, the computational errors are approximately monotone functions of $K$.
Choosing $K$ larger than enough will simply lead to an increase 
in CPU time but will not diminish the accuracy.
 In Examples 1 and 2, where $\xi\le\sqrt{3}$, 
 it suffices to pick $K$ smaller than in Example 3 where $\xi\le 10$.
\item
In Examples 4 and 5, where the SDEs are nonlinear, the dependences of numerical errors on $K$ are not monotone.
The errors will significantly exceed the minimal possible ones if $K$ is either too small or too large.
The good news is that if $\xi(x)$  is rather small like in Example 4, there are
 large ranges of values of $K$ for which errors are nearly minimal. 
 Therefore, if an upper bound for $\xi$ is unknown a priori, 
 it is safer to choose $K$ assuming that $\xi\sim10$.
\item
It might seem surprising that the errors in Examples 3 and 5 where $\xi$ is large reach smaller values than the errors in 
Examples 2 and 4 respectively where $\xi$ is small.
\end{itemize}
Based on Fig. \ref{fig1:ex}, a guideline for choosing $K$ is given in Table \ref{table:K}.

\begin{table}[htp]
\caption{A guideline for choosing $K$ for mesh sizes $N^3$. }
\begin{center}
\begin{tabular}{|c|c|c|c|c|c|}
\hline
$N$& 33 & 65 & 129 & 257 & 513 \\
\hline
$K$& 4 & 6 & 8 & 10 & 14\\
\hline
\end{tabular}
\end{center}
\label{table:K}
\end{table}%

\begin{figure}[htbp]
\begin{center}
(a)\includegraphics[width = 0.45\textwidth]{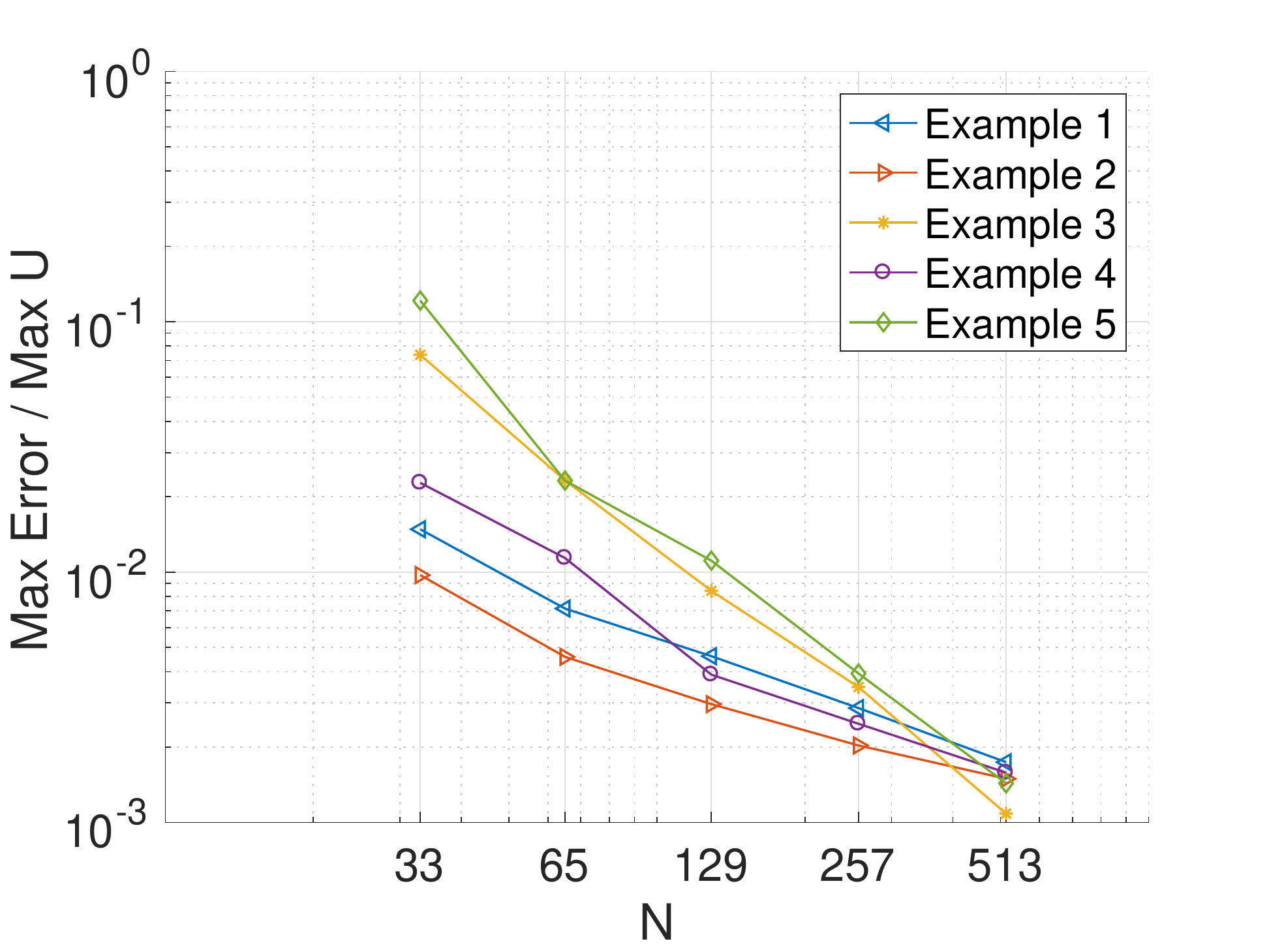}
(b)\includegraphics[width = 0.45\textwidth]{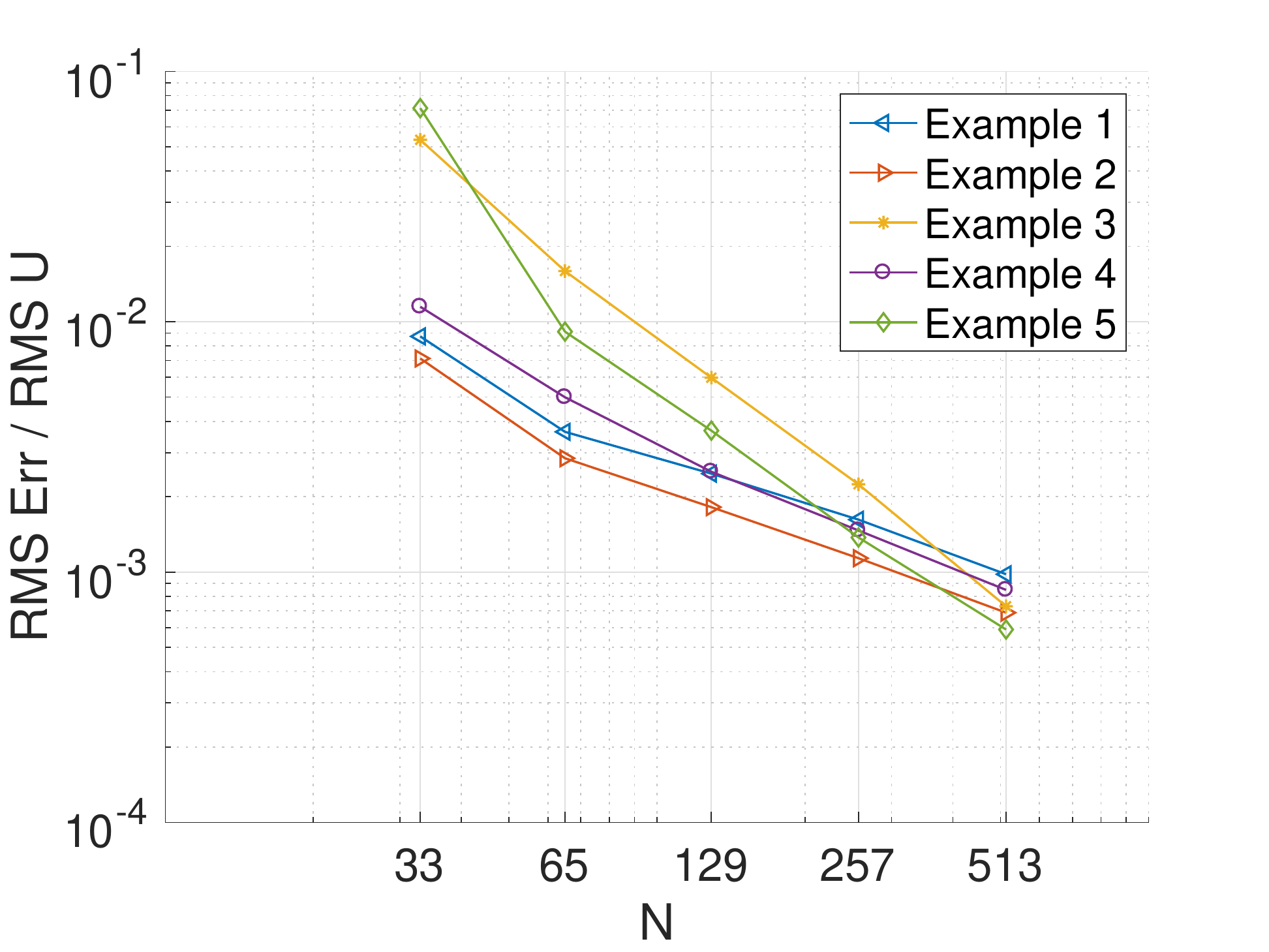}\\
(c)\includegraphics[width = 0.45\textwidth]{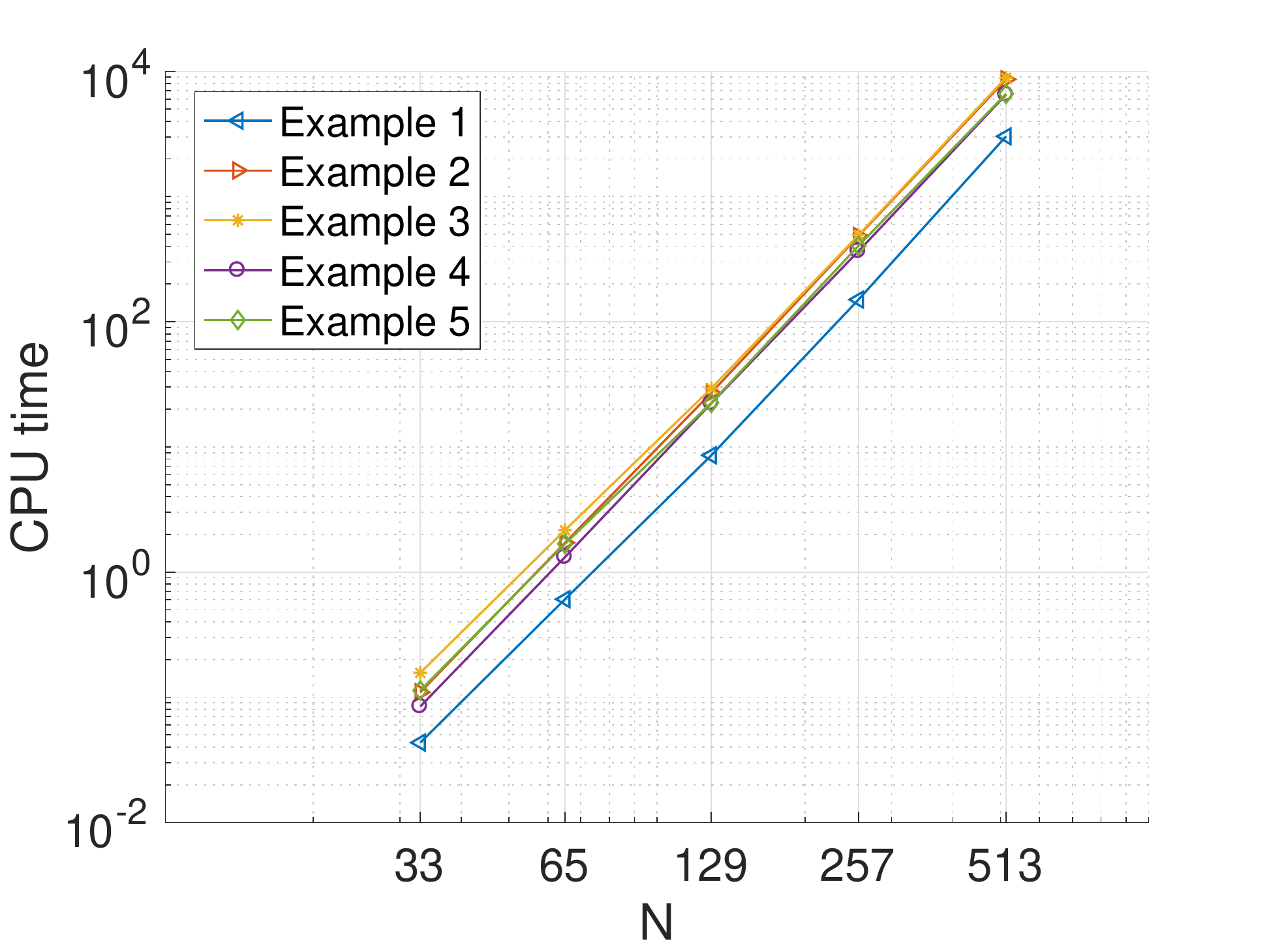}
(d)\includegraphics[width = 0.45\textwidth]{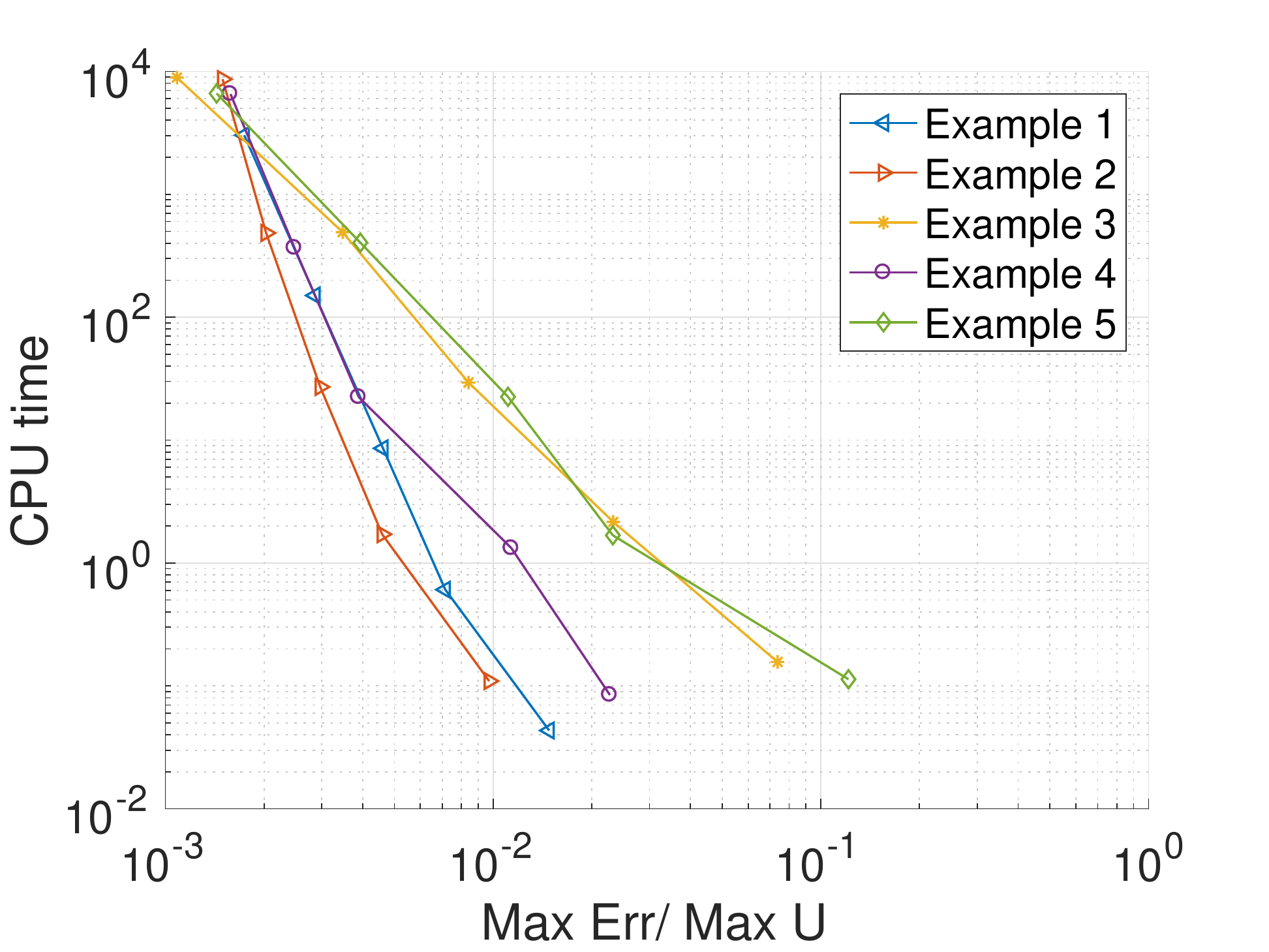}
\caption{(a): The normalized maximum absolute error versus $N$, where the mesh size is $N^3$.
(b): The normalized RMS error versus $N$.
(c): The CPU time versus $N$.
(d): The CPU time versus the normalized maximum absolute error.
}
\label{fig2:ex}
\end{center}
\end{figure}

The next series of graphs shown in Fig. \ref{fig2:ex} is done for the values of $K$ chosen according to Table \ref{table:K}.
 Fig. \ref{fig2:ex} (c) displays the dependence of the CPU time on $N$. 
Table \ref{table:tests} shows that the CPU time  $T$ grows a bit slower than quadratically with $K$ for a fixed $N$,
and approximately as the fourth power of $N$ for an optimal choice of $K$ for each $N$.

Table \ref{table:tests} also shows that the observed rate of convergence is superlinear for Examples 3 and 5 where the ratio $\xi$ is large, 
while they are sublinear for Examples 1, 2, and 4 where $\xi$ is relatively small.
We will give an explanation for this phenomenon in Section \ref{sec:why}.

\begin{table}[htp]
\caption{Least squares fits for Examples 1--5. $T$ is the CPU time, $N^3$ is the mesh size, $h$ is the mesh step, $E_{NM}$ 
is the normalized maximal absolute error:
$E_{NM} : = E_{\max}/U_{\max}$. $E_{NR}$ is the normalized maximal RMS error: $E_{NR}: = E_{RMS}/U_{RMS}$.}
\begin{center}
\begin{tabular}{|p{3.5cm}|c|c|c|c|c|}
\hline
Example \#& 1&2&3&4&5 \\
\hline
$T$ vs $K$ for $N = 513$, seconds & $31.4K^{1.79}$ &  $ 8.57K^{1.80}$ & $ 93.8K^{1.78}$ & $ 66.3K^{1.80}$ & $ 1.94K^{1.80}$ \\
\hline
$E_{NM}$ vs $h$ for optimal $K$ & $0.181h^{0.750}$ & $0.0808h^{0.658}$ & $ 12.3h^{1.49}$ & $0.635h^{0.989}$ & $19.4h^{1.54}$ \\
\hline
$E_{NR}$ vs $h$ for optimal $K$ & $0.0987h^{0.748}$ & $0.0983h^{0.807}$ & $ 9.74h^{1.52}$ & $0.255h^{0.928}$ & $14.0h^{1.66}$ \\
\hline
$T$ vs $N$ for optimal $K$, nanoseconds & $27.4N^{4.05}$ & $61.0N^{4.11}$ & $131.0N^{3.98}$ & $49.2N^{4.10}$ & $92.6N^{4.00}$ \\
\hline
\end{tabular}
\end{center}
\label{table:tests}
\end{table}%


\subsection{Effects of local factoring}
\label{sec:ELF}
\begin{figure}[htbp]
\begin{center}
(a)\includegraphics[width = 0.4\textwidth]{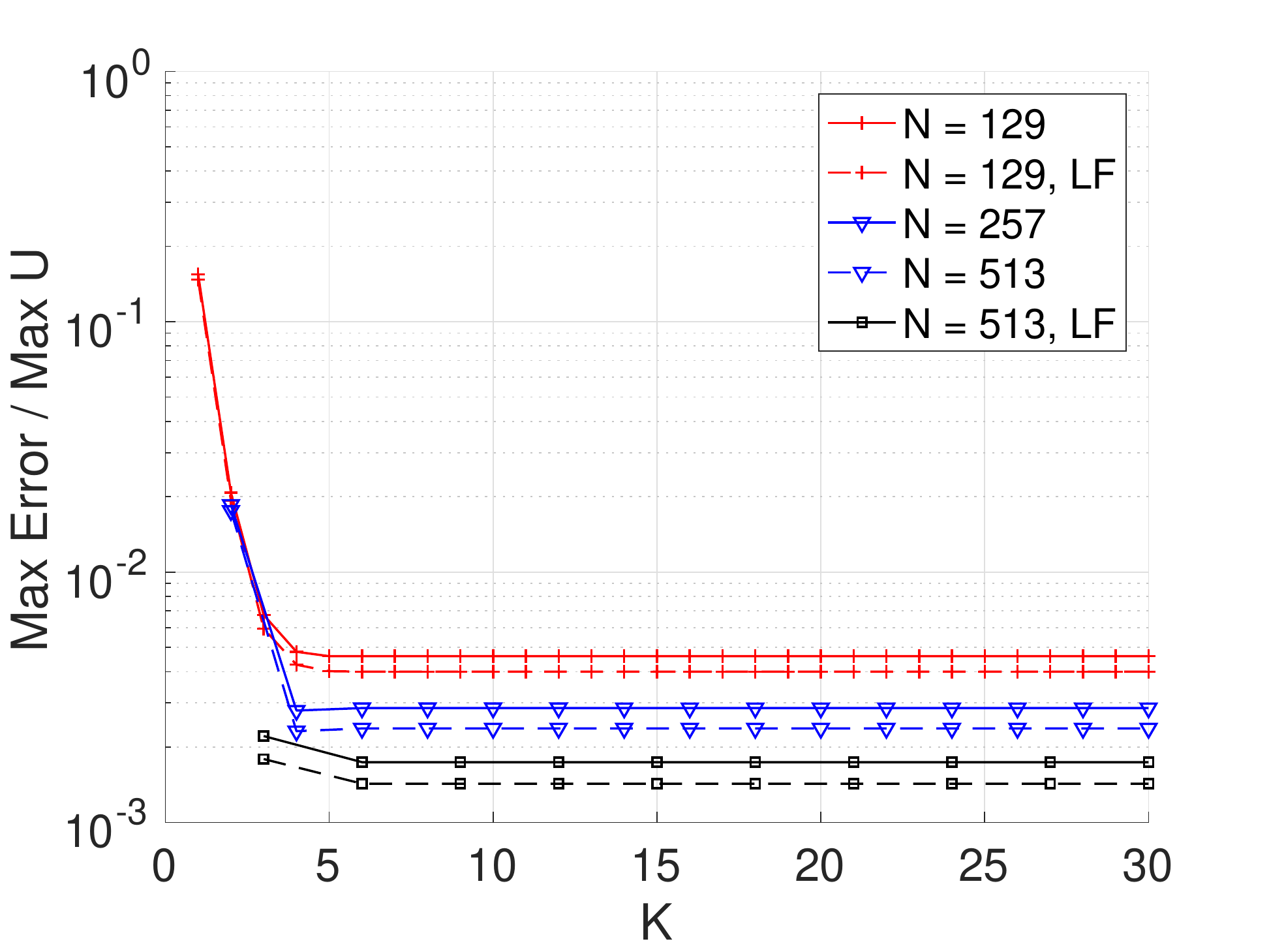}
(b)\includegraphics[width = 0.4\textwidth]{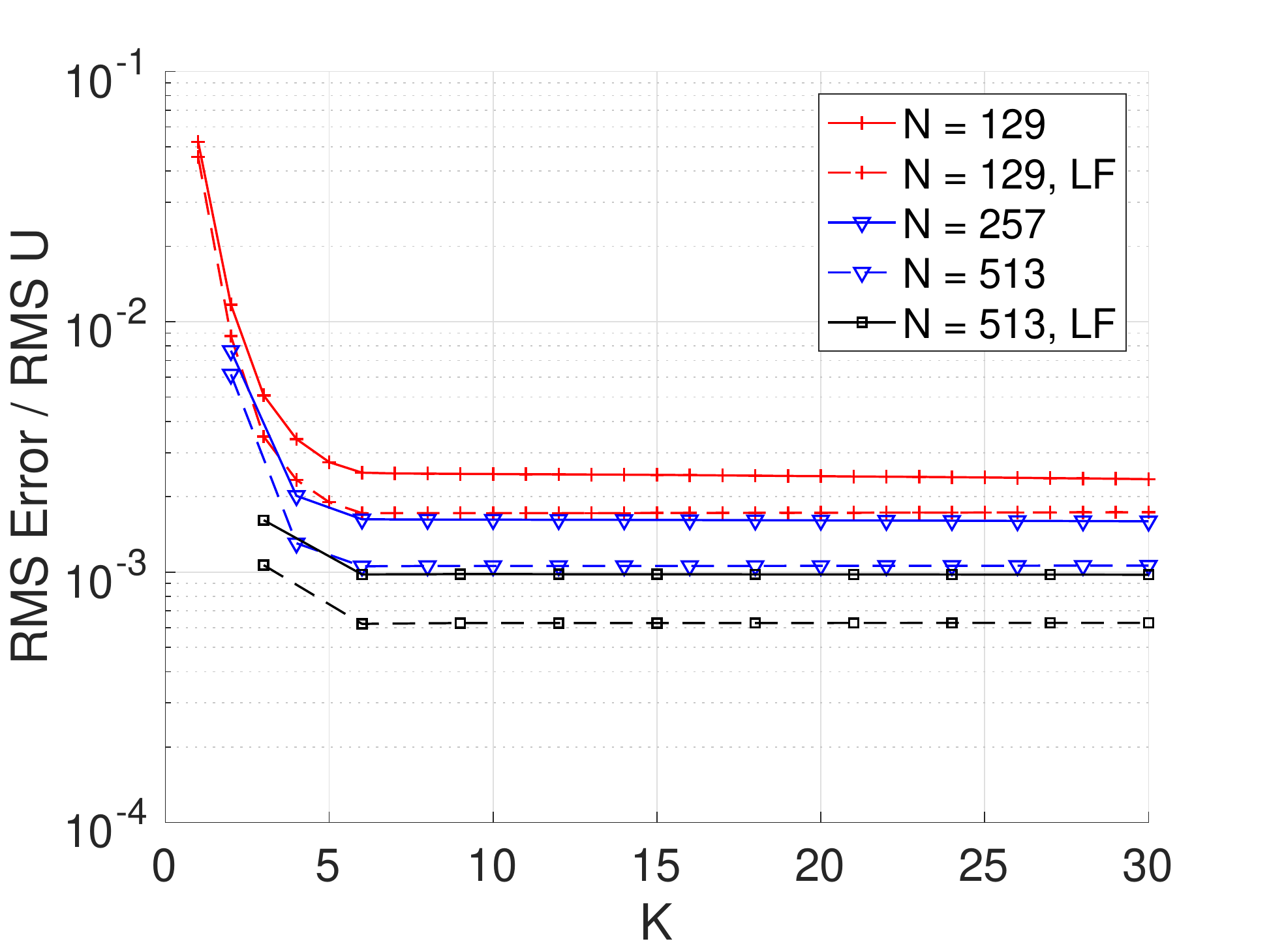}\\
(c)\includegraphics[width = 0.4\textwidth]{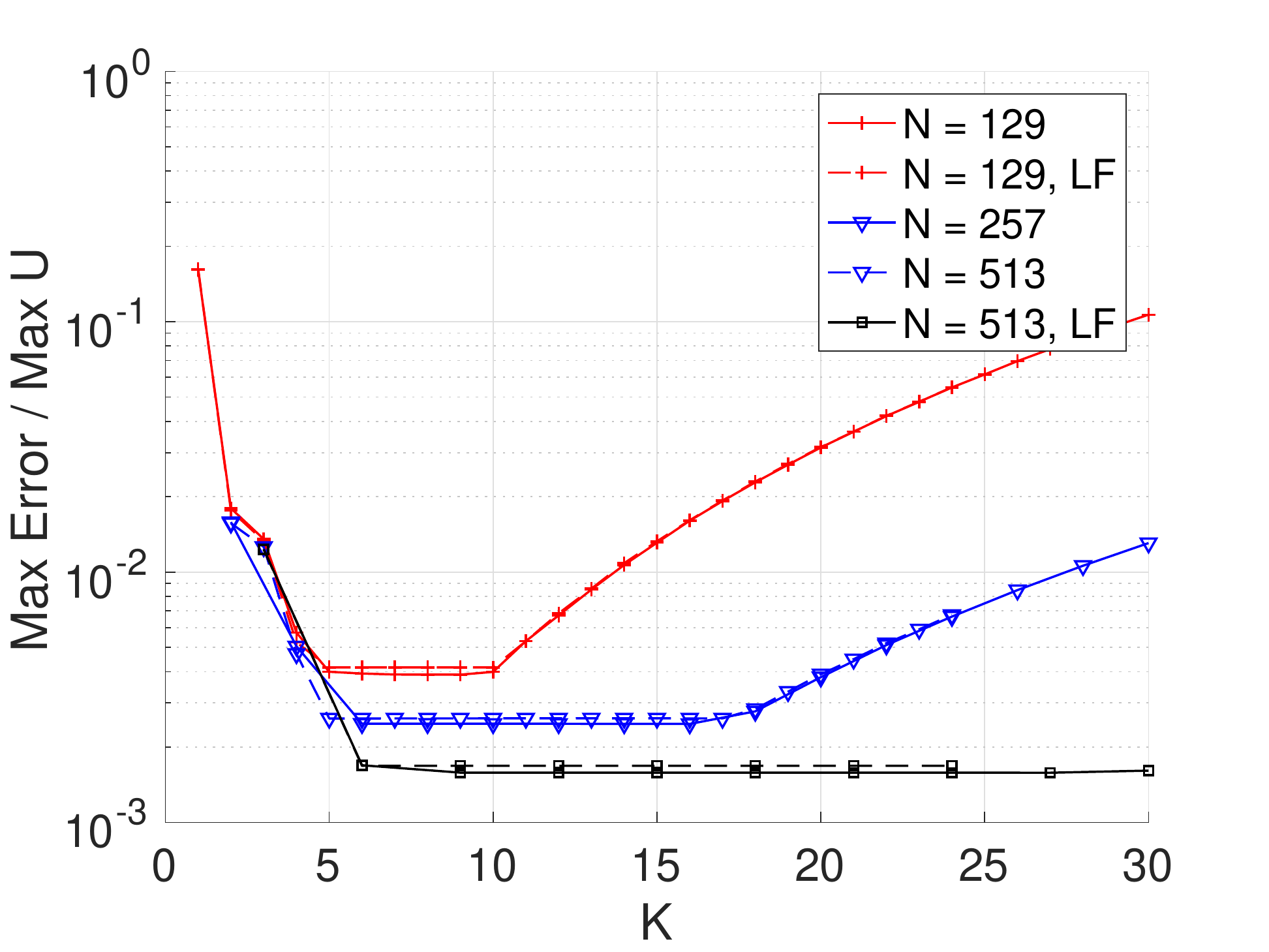}
(d)\includegraphics[width = 0.4\textwidth]{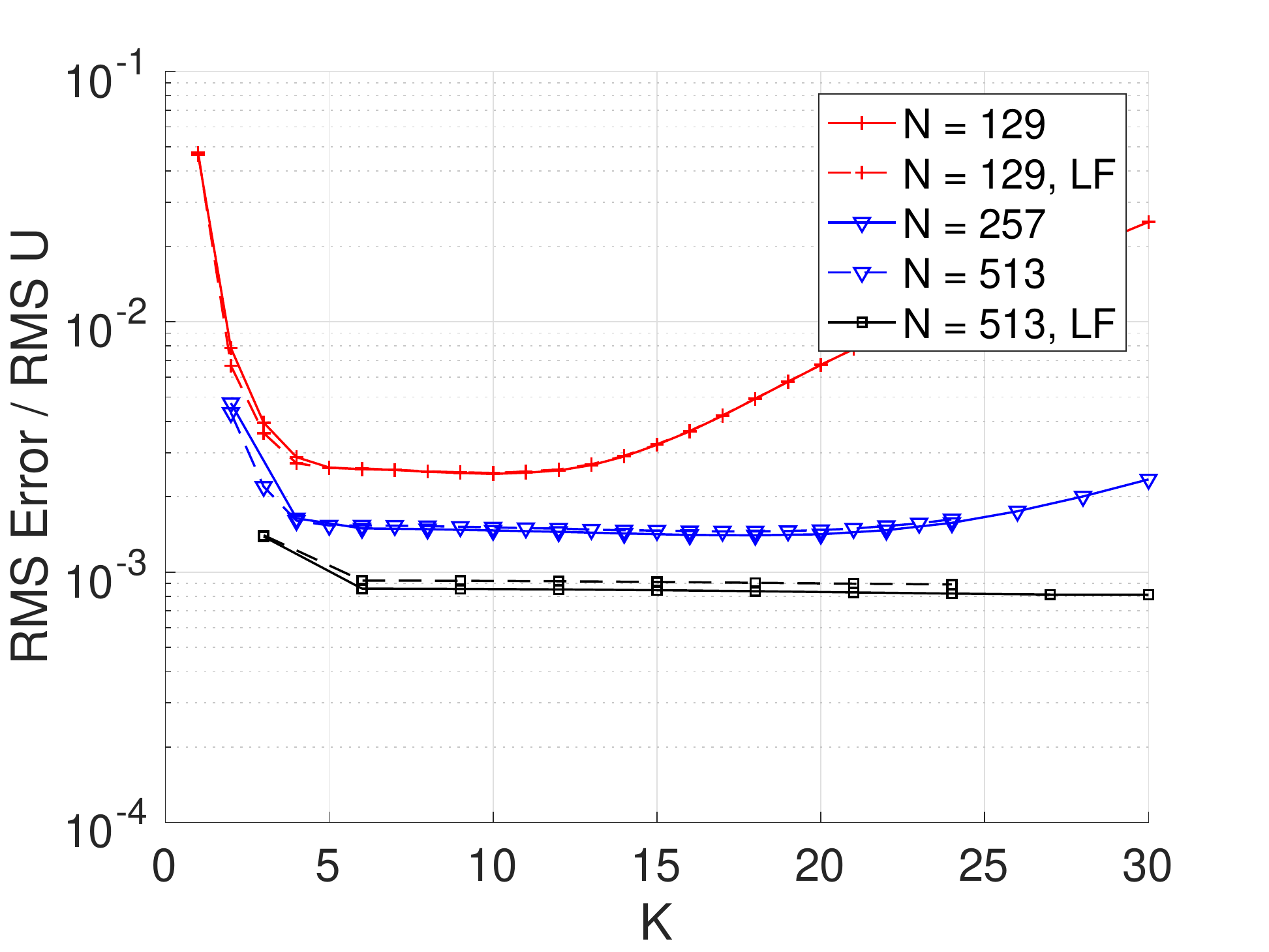}\\
(e)\includegraphics[width = 0.4\textwidth]{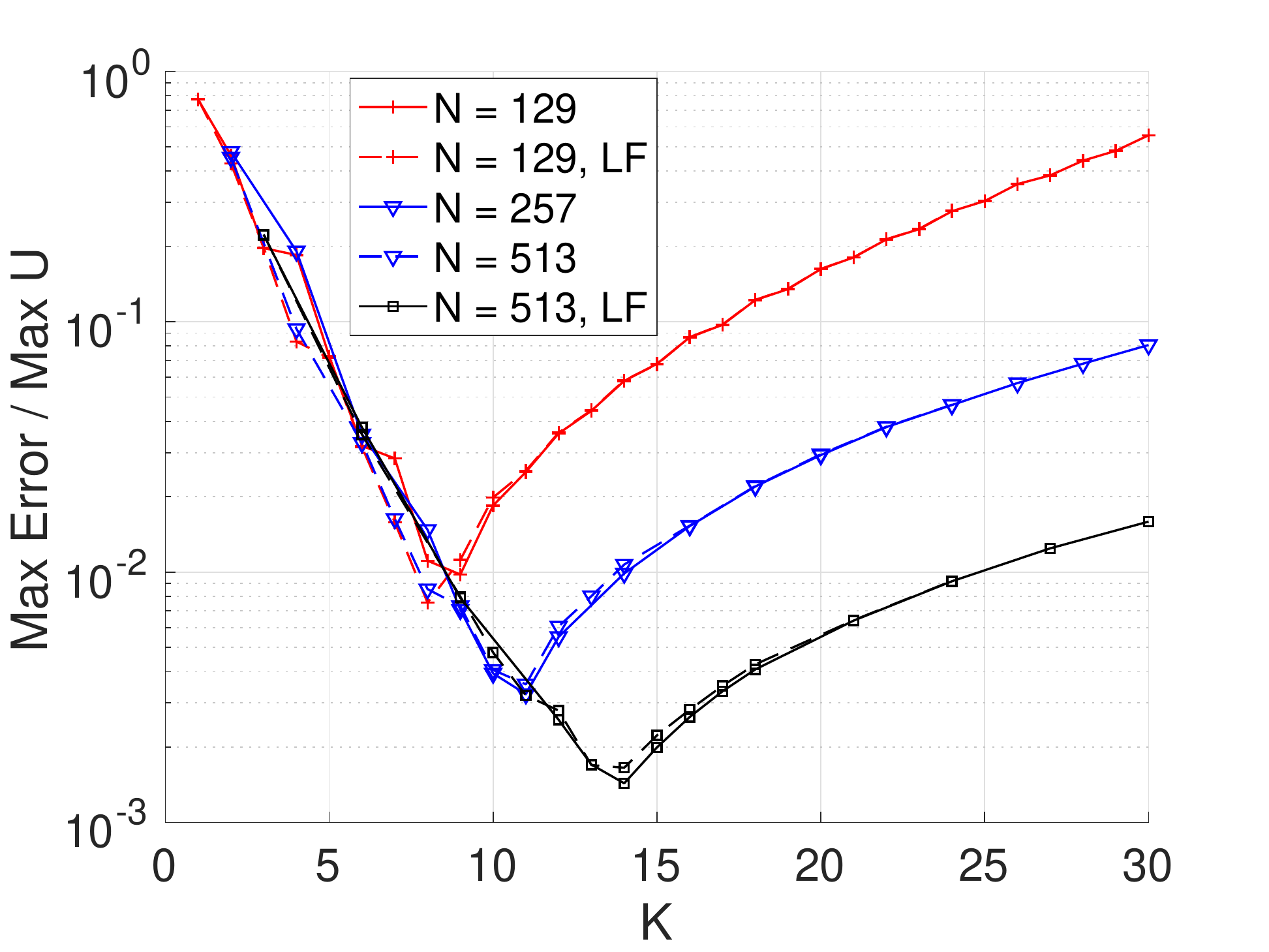}
(f)\includegraphics[width = 0.4\textwidth]{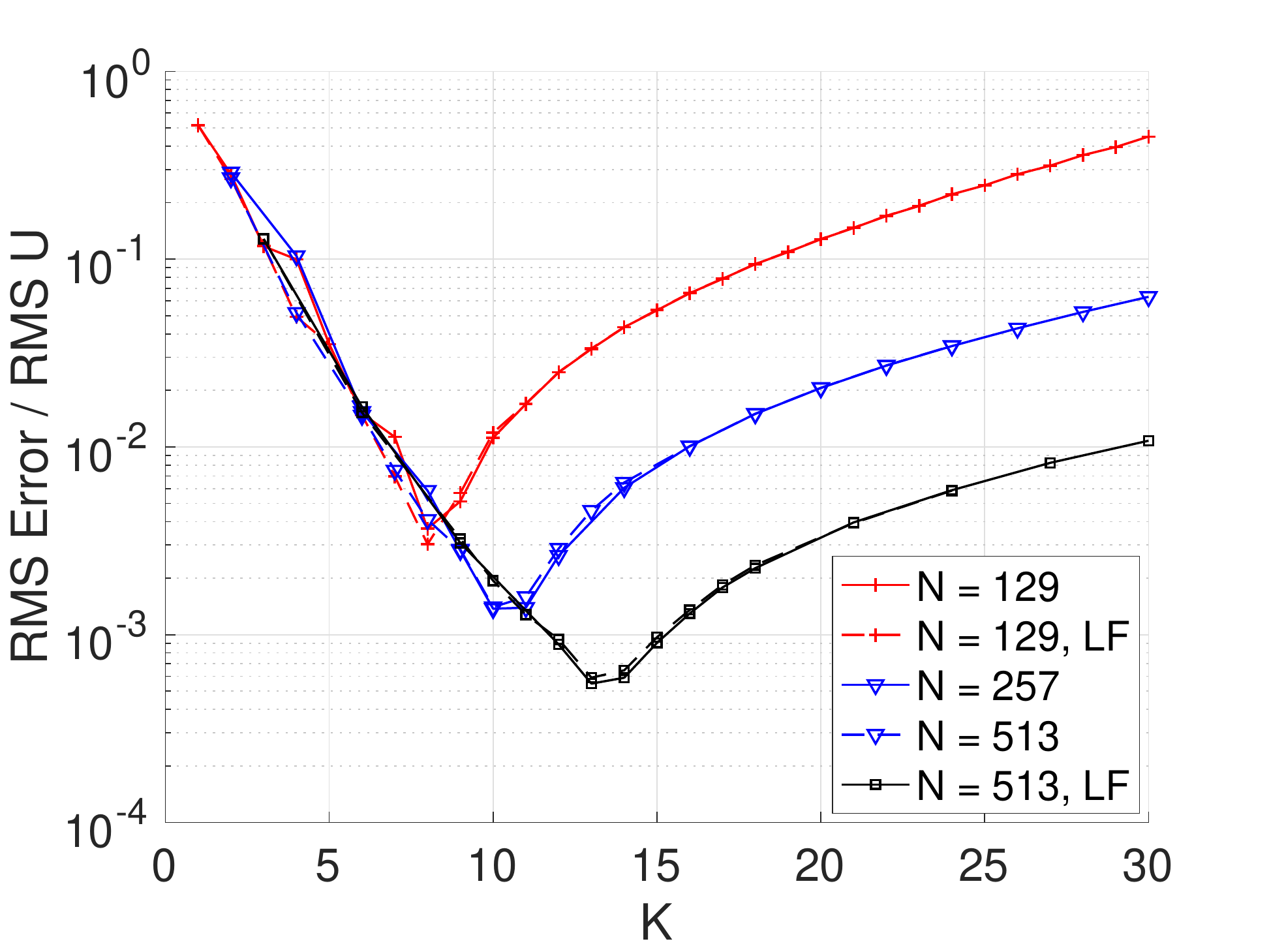}
\caption{Effects of local factoring. Dashed and solid curves depict errors committed by {\tt olim3D} 
with and without local factoring respectively.
(a--b): Example 1. (c--d): Example 4. (e--f): Example 5.
(a,c,e): The normalized maximum absolute error versus the update parameter $K$.
(b,d,f): The normalized RMS error versus $K$.
}
\label{fig:ELF}
\end{center}
\end{figure}

We have implemented two versions of {\tt olim3D}, with and without local factoring, as described in Section \ref{sec:locfac}.
The radius of the neighborhood where local factoring is done is $R_{\rm fac} = 0.1$. 
Our measurements are performed on Examples 1, 4, and 5.
Our results are presented in Fig. \ref{fig:ELF}.
The increase of CPU times due to local factoring is within natural variations of CPU times.
Our measurements show that the computational errors noticeably decrease due to local factoring for the linear SDE in Example 1,
while  they are almost unaffected--even they may slightly increase--for the nonlinear SDEs in Examples 4 and 5. 
As a result, we conclude that local factoring is not worth implementing routinely in quasipotential solvers, but may be helpful in some cases.
One such case is the genetic switch model presented in Section \ref{sec:gene}.

\subsection{Systems with hyperbolic periodic orbits}
\label{sec:tao}
An interesting class of SDEs was considered in \cite{tao}.
The corresponding ODE $\dot{x} = b(x)$ 
has two asymptotically stable equilibria whose basins of attraction are separated by a manifold $\Sigma$ 
containing a periodic orbit $C$. This orbit is the only attractor of $\dot{x}=b(x)$ restricted to $\Sigma$. 
The orbit $C$ is the transition state for the corresponding SDE. Additionally, $b(x)$
admits an orthogonal decomposition 
\begin{equation}
\label{ortdec}
b(x) = -\frac{\nabla v(x)}{2} + r(x),\quad {\rm where}\quad \nabla v(x) \cdot r(x) = 0.
\end{equation}
The function $v(x)$ is not necessarily the quasipotential, however, it is a solution of the Hamilton-Jacobi  \eqref{HJ} which is easy to check.
Let $x^{\ast}$ be one of the asymptotically stable equilibria of $\dot{x}=b(x)$, and let $y$ be a point for which there exists a path $\zeta(s)$ 
from $x^{\ast}$ to $y$ such that 
\begin{equation}
\label{zetapath}
\zeta'(s) = \frac{b(s) + \nabla v(s)}{\|b(s) + \nabla v(s)\|}\equiv \frac{\frac{\nabla v(s)}{2} + r(s)}{\left\|\frac{\nabla v(s)}{2} + r(s)\right\|}.
\end{equation}
It is shown in \cite{tao} that the quasipotential at $y$ with respect to $x^{\ast}$ coincides with $v(y) - v(x^{\ast})$.
This fact can be useful 
in the case where the quasipotential is hard  or impossible to find analytically, but an orthogonal decomposition of the form \eqref{ortdec} is easy to spot.
We will call the function $v(x)$ that satisfies the boundary condition $v(x^{\ast}) = 0$ and the Hamilton-Jacobi  \eqref{HJ} but is not the quasipotential 
(i.e., does not solve the minimization problem \eqref{Qpot}) a \emph{false quasipotential}.

We apply {\tt olim3D} to two 3D examples from \cite{tao}. 

{\bf Example 6.}
The vector field is given by:
\begin{equation}
\label{tao1}
b(x) = \left[\begin{array}{c} (1-x_3)^2\frac{x_1}{\sqrt{x_1^2+x_2^2}} -x_1-x_2 \\
(1-x_3)^2\frac{x_2}{\sqrt{x_1^2+x_2^2}} -x_2 + x_1\\
x_3 - x_3^3
\end{array}\right]
\end{equation}
The asymptotically stable equilibria are $(0,0,-1)$ and $(0,0,1)$. The transition state is the hyperbolic periodic orbit
$C = \{x~|~ x_1^2 + x_2^2 = 1, ~x_3 = 0\}$.
A false quasipotential is given by $v(x) = \tfrac{1}{2}(1-x_3^2)^2$. 
Obviously, $v(0,0,-1) = 0$, and $v(x) = \tfrac{1}{2}$ at any point such that $x_3 = 0$; in particular, on the hyperbolic periodic orbit $C$.
Hence, the true quasipotential $U$ is also $\tfrac{1}{2}$ at $C$.

For this example, we have chosen   the cube $[-2,2]^3$ to be the computational domain, the mesh size to be
 $513^3$, and the update factor $K=25$.
Note that $\xi(x)$ approaches $\infty$ in the neighborhood of $C$.

The vector field $b(x)$ is non-differentiable at $(0,0,-1)$. 
For initialization purposes, we have 
replaced the non-existent partial derivatives $\tfrac{db_1}{dx_3}$ and $\tfrac{db_2}{dx_3}$  in the Jacobian matrix
with zeros. The computations were performed with and without local factoring.

\begin{figure}[htbp]
\begin{center}
(a)\includegraphics[width=0.45\textwidth]{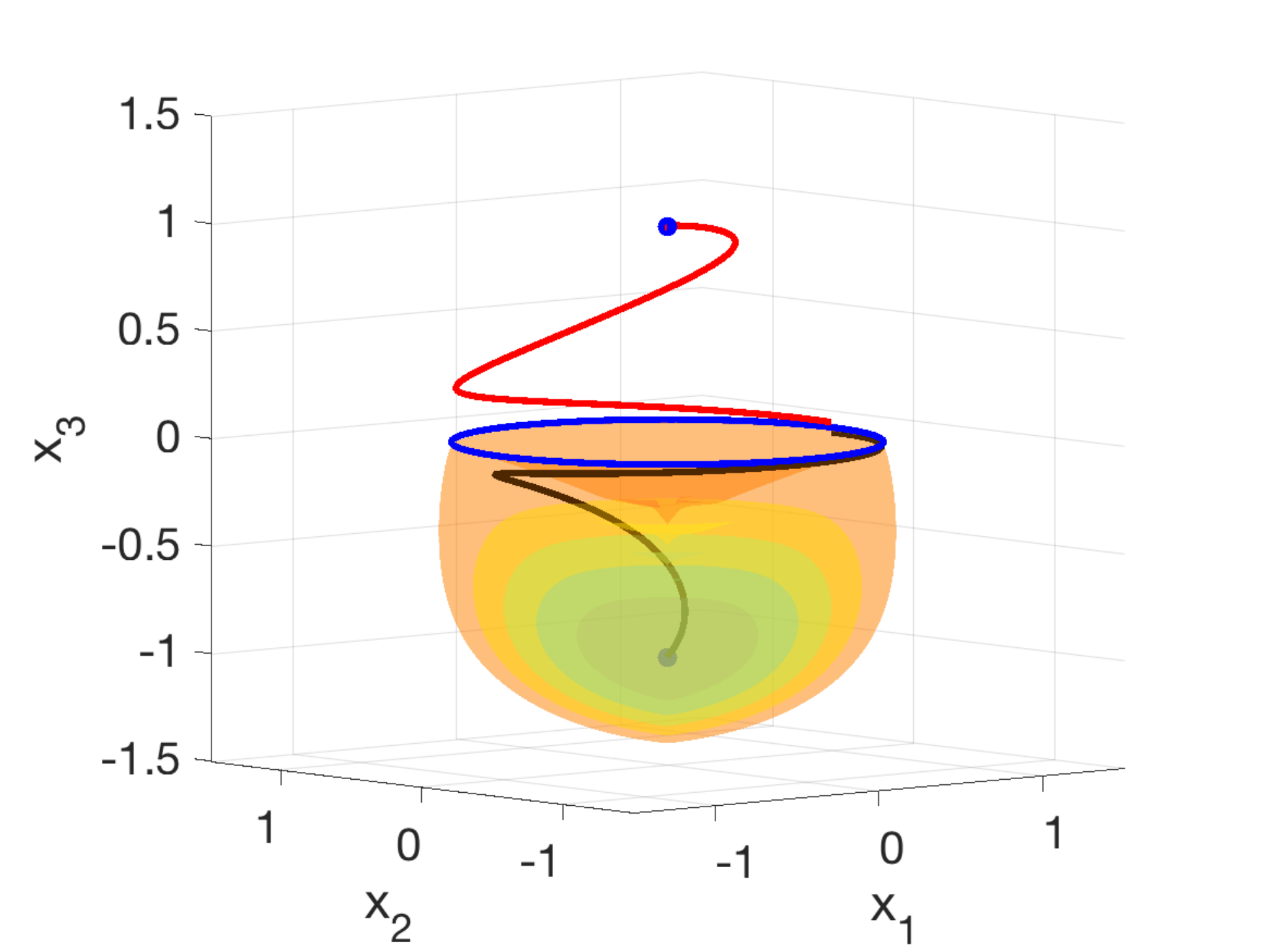}
(b)\includegraphics[width=0.45\textwidth]{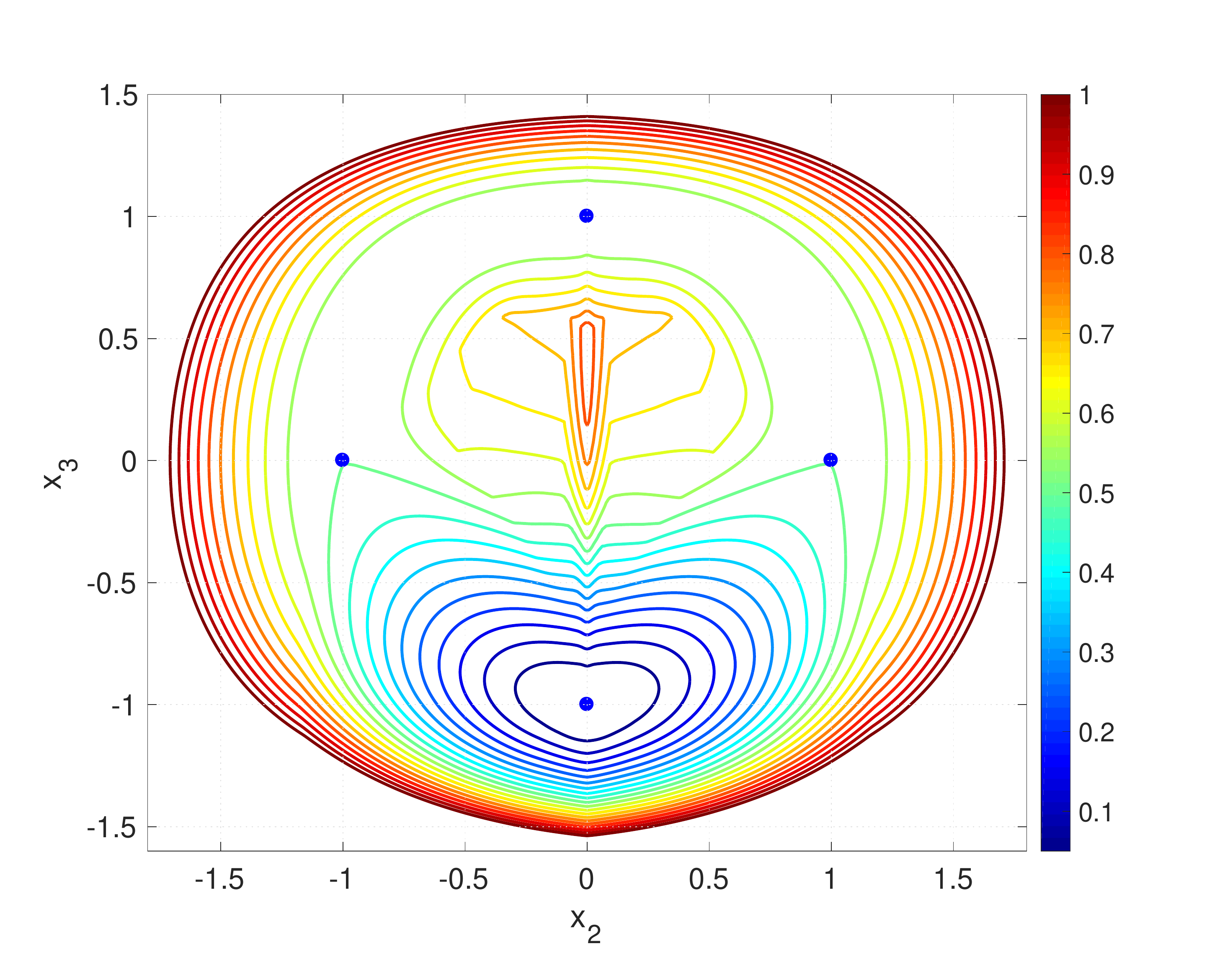}
(c)\includegraphics[width=0.45\textwidth]{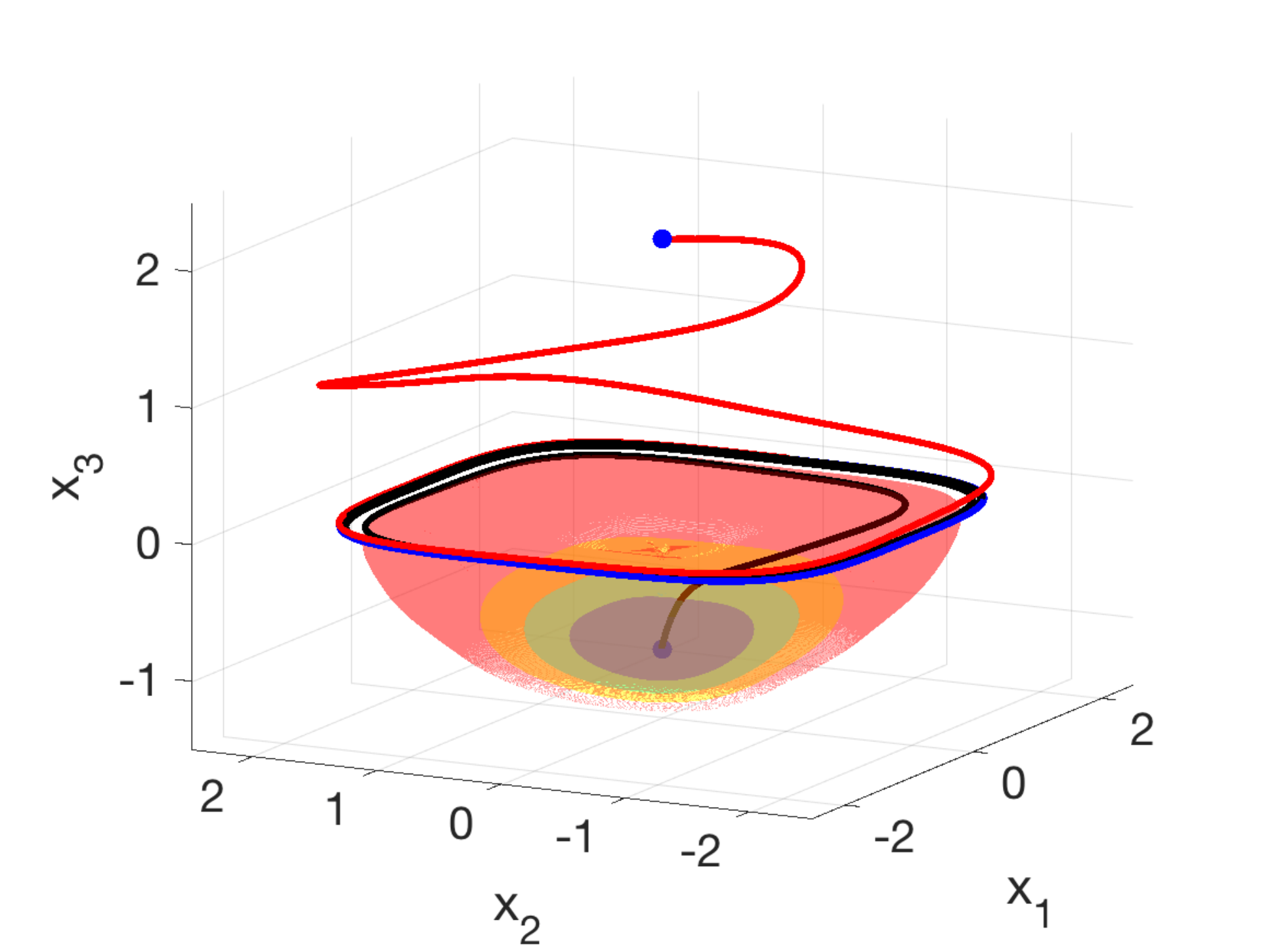}
(d)\includegraphics[width=0.45\textwidth]{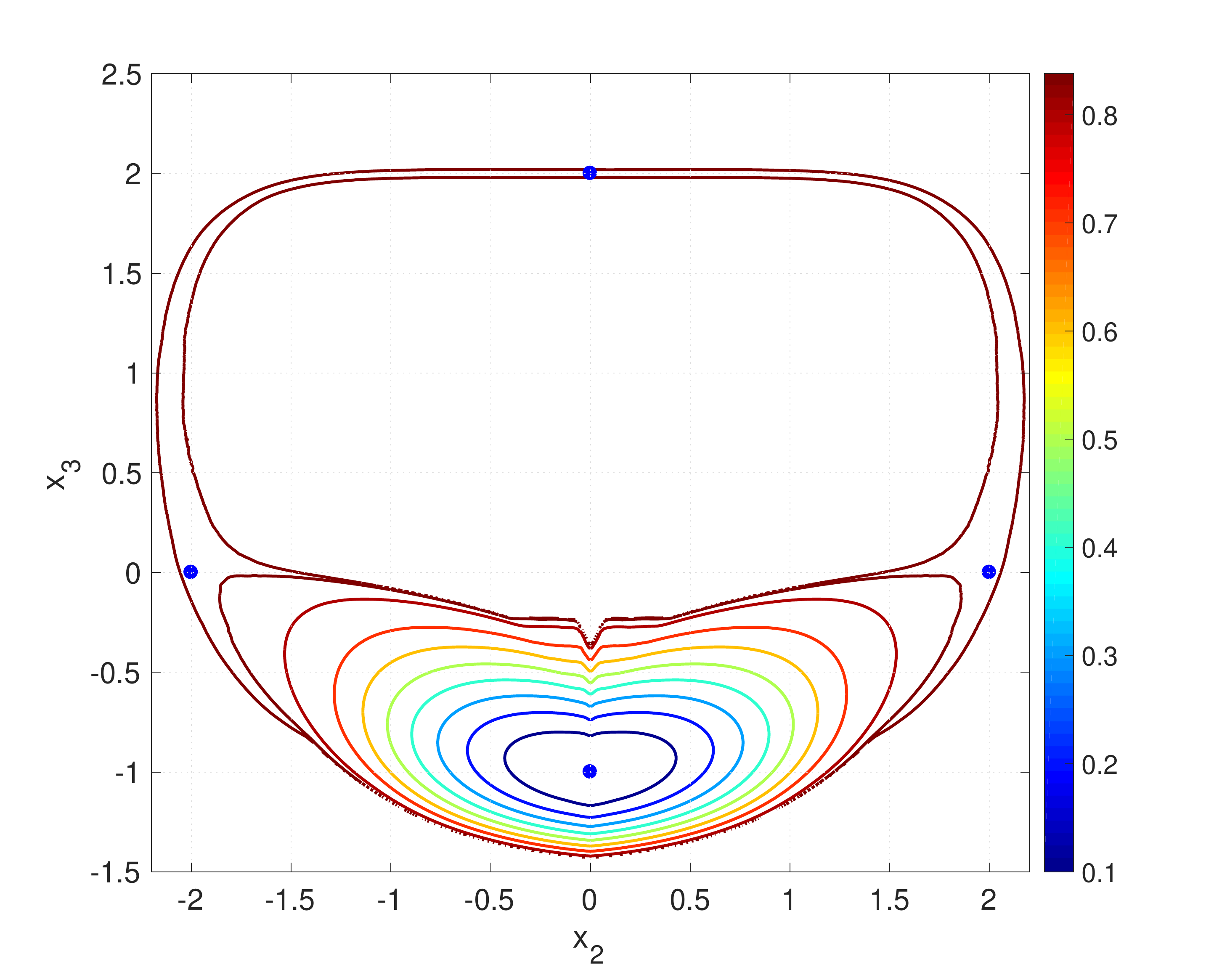}
\caption{(a--b): The level sets of the quasipotential and the MAP for Example 6. 
(a): The level sets correspond to the following values of $U$: 0.1, 0.2, 0.3, 0.4, and 0.5.
The blue curve is the hyperbolic periodic orbit $C =  \{(x~|~ x_1^2 + x_2^2 = 1, ~x_3 = 0\}$.
The black curve is a MAP from $(0,0,-1)$ to $C$ obtained by numerical integration. The red curve
is a trajectory starting near $C$ and going to $(0,0,1)$.  
(b): The intersection of level sets of the quasipotential with $yz$-plane.
The blue dots represent the equilibria and the intersection of $C$ with the $yz$-plane.
(c--d): The level sets of the quasipotential and the MAP for Example 7. 
(c): The level sets correspond to the following values of $U$: 0.2, 0.4, 0.6,  and $\tfrac{5}{6}$.
The blue curve is the hyperbolic periodic orbit $C =  \{(x~|~ x_1^4 + x_2^4 = 16, ~x_3 = 0\}$.
The black curve is a MAP from $(0,0,-1)$ to $C$ obtained by numerical integration. The red curve
is a trajectory starting near $C$ and going to $(0,0,2)$.  
(d): The intersection of level sets of the quasipotential with $yz$-plane.
The blue dots represent the equilibria and the intersection of $C$ with the $yz$-plane.
The two outermost contours correspond to $U = \tfrac{5}{6} \approx 0.833$ and $U = 0.838$ respectively.
}
\label{fig:tao}
\end{center}
\end{figure}

Level sets of the computed quasipotential corresponding to $U(x) = 0.1$, 0.2, 0.3, 0.4, and 0.5 and the MAP connecting the two equilibria  
obtained by numerical integration (see Section \ref{sec:background})
are shown in Figs. \ref{fig:tao} (a).
Level sets of the quasipotential intersected with the $yz$-plane 
are plotted in Fig. \ref{fig:tao} (b).
We found the value of the quasipotential at the transition state $C$ by linear interpolation from the computed quasipotential.
The maximal absolute error at $C$ is $1.23\times 10^{-4}$ for the computation without local factoring and $3.58\times 10^{-4}$ for the one with it.
The CPU time  was 9.9 hours.


{\bf Example 7.}
Here, unlike in Example 6, the  vector field is not rotationally symmetric:
\begin{equation}
\label{tao2}
b(x) = \left[\begin{array}{c} -(x_3+1)(x_3-2)\frac{x_1}{(x_1^4+x_2^4)^{1/4}} -x_1-x_2^3 \\
 -(x_3+1)(x_3-2)\frac{x_2}{(x_1^4+x_2^4)^{1/4}} +x_1^3 -  x_2\\
 -(x_3+1)(x_3-2)x_3
\end{array}\right]
\end{equation}
The asymptotically stable equilibria are $(0,0,-1)$ and $(0,0,2)$;  transition state is the hyperbolic periodic orbit
$C = \{x~|~ x_1^4 + x_2^4 = 16, ~x_3 = 0\}$.
A false quasipotential $v(x) = \tfrac{x_3^4}{2} - \tfrac{2x_3^3}{3} - 2x_3^2 + \tfrac{5}{6}$ is such that 
$v(0,0,-1) = 0$, and $v(x) = \tfrac{5}{6}$ at any point in the plane $x_3 = 0$.
Hence, the true quasipotential $U$ is also equal to $\tfrac{5}{6}$ at $C$.

We have chosen $[-2.2,2.2]^2\times[-1.5,2.5]$ computational domain and $661^2\times 601$ mesh, and set $K=25$.
The vector field $b(x)$ is non-differentiable at $(0,0,-1)$ like in Example 6, and the initialization has been performed similarly.

Level sets of the computed quasipotential corresponding to $U(x) = 0.2$, 0.4, 0.6,  and $\tfrac{5}{6}$ and the MAP connecting the two equilibria  
obtained by numerical integration 
are shown in Figs. \ref{fig:tao} (c).
Level sets of the quasipotential intersected with the $yz$-plane 
are plotted in Fig. \ref{fig:tao} (d).
Fig. \ref{fig:tao} (d) shows that the quasipotential changes very slowly in a quite large neighborhood of $C$.
Therefore, small inaccuracy in computing the quasipotential visibly shifts the level set $U = \tfrac{5}{6}$ 
in Fig. \ref{fig:tao} (c) that should have passed through the
blue curve representing $C$.

As in Example 6, the quasipotential at the transition state $C$  is found by linear interpolation.
Without local factoring, the maximal absolute error at $C$ is $1.32\cdot 10^{-3}$, while with local factoring, it decreases to $1.02\cdot 10^{-3}$.
The CPU time was 10.7 hours.

\subsection{An application to a genetic switch model}
\label{sec:gene}
Lv et al. \cite{lv} studied a 3D genetic switch model with two metastable states and positive feedback. 
They justified an application of Donsker-Varadhan type large deviation theory \cite{DV1,DV2} to this model, and used GMAM \cite{hey1,hey2}
to find MAPs and a 2D quasipotential landscape. The latter was done by taking the maximum of the quasipotential
with respect to one of the variables. 
In this work, we perturb the ODEs from \cite{lv} with small isotropic white noise and use them as a test model:
\begin{align}
d\mathsf{m}  & = \left(\frac{a_0\gamma_0 + ak_0\mathsf{d}}{\gamma_0 + k_0\mathsf{d}} - \gamma_m\mathsf{m}\right) + \sqrt{\epsilon}dw_1,  \notag \\
d\mathsf{n} & = \left(b\mathsf{m} - \gamma_n \mathsf{n} -2k_1\mathsf{n}^2 + 2\gamma_1\mathsf{d}\right) + \sqrt{\epsilon}dw_2 , \label{sdeMND}\\
d\mathsf{d} & = \left(k_1\mathsf{n}^2 - \gamma_1\mathsf{d}\right) + \sqrt{\epsilon}dw_3.\notag
 \end{align}
Here, $\mathsf{m}$, $\mathsf{n}$, and $\mathsf{d}$ are the numbers of the mRNA, protein, and dimer formed by the protein, respectively.
 The parameters were taken from Supplement S1 for \cite{lv}. 
 This system has two asymptotically stable equilibria $x_i:=[\mathsf{m}_i,\mathsf{n}_i,\mathsf{d}_i]$ and 
 $x_a=[\mathsf{m}_a,\mathsf{n}_a,\mathsf{d}_a]$ representing inactive and active states respectively,
 separated by a Morse index one saddle $x_s=[\mathsf{m}_s,\mathsf{n}_s,\mathsf{d}_s]$:
 \begin{align}
x_i &= [0.040206714231704188,1.6082685692681673,0.00025865277908958782], \notag \\
x_a &= [29.376860080598071,1175.0744032239231,138.07998531120592], \label{lveq} \\
x_s &= [10.5829, 423.3173, 17.9198].\notag
 \end{align}

\begin{figure}[htbp]
\begin{center}
(a)\includegraphics[width=0.45\textwidth]{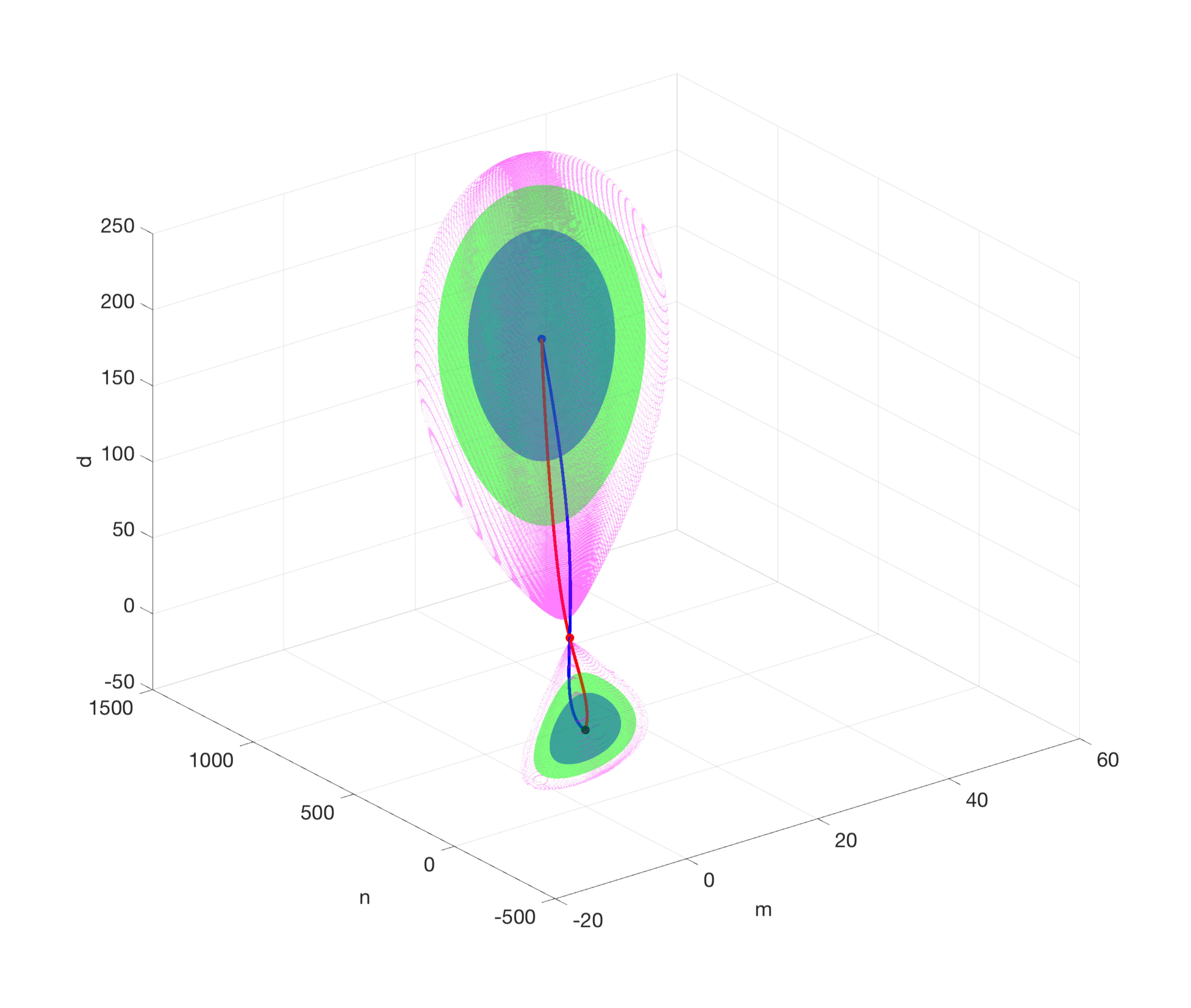}
(b)\includegraphics[width=0.45\textwidth]{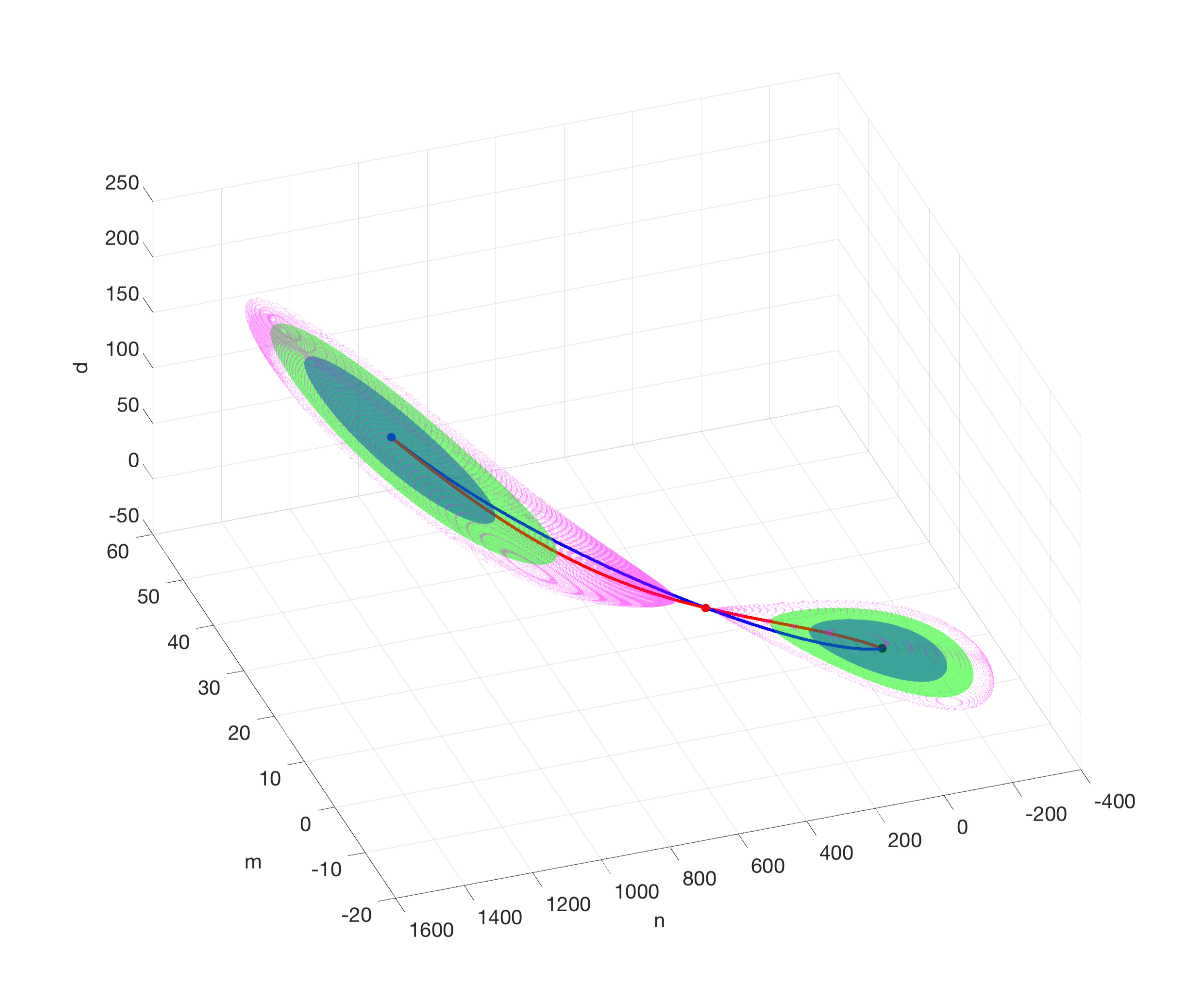}
\caption{(a--b) Two views of the level sets of the quasipotentials for a genetic switch model \cite{lv} (SDE \eqref{sdeMND})
 computed with respect to each equilibrium.
The red curves are the MAPs from the equilibria to the saddle. The blue curves are the trajectories from the saddle to the equilibria.
The equilibria with small  and large values of the variables $\mathsf{m}$, $\mathsf{n}$, and $\mathsf{d}$ represent the inactive and active
states, respectively. 
The level sets correspond to the values  $0.33\times U(x_s)$,  $0.66\times U(x_s)$, and  $0.99\times U(x_s)$ where $U(x_s)$ is the value
of the quasipotential at the saddle for each computation.
}
\label{fig:gene}
\end{center}
\end{figure}

We have computed the quasipotential with respect to each equilibrium.
The computational domains for the runs with the initial points at  $x_i$ and $x_a$ are the parallelepipeds
\begin{align*}
&[\mathsf{m}_i-18,\mathsf{m}_i-18]\times[\mathsf{n}_i-450,\mathsf{n}_i-450]\times[\mathsf{d}_i-36,\mathsf{d}_i-36],~~{\rm and}\\
&[\mathsf{m}_a-36,\mathsf{m}_a-36]\times[\mathsf{n}_a-792,\mathsf{n}_a-792]\times[\mathsf{d}_a-180,\mathsf{d}_a-180],
\end{align*}
respectively. They are set up so that the runs terminate as soon as the saddle $x_s$ is reached.
The mesh sizes were $65\times1601\times129$, $129\times3201\times257$, and  $257\times6401\times513$
for the computation of the quasipotential with respect to $x_i$, and $129\times2817\times641$ for the computation 
initialized at $x_a$. The level sets of the computed quasipotentials are shown in Fig. \ref{fig:gene}.
They correspond to the values  $0.33\times U(x_s)$,  $0.66\times U(x_s)$, and  $0.99\times U(x_s)$ where $U(x_s)$ is the value
of the quasipotential at the saddle for each computation. 

While the quasipotential for SDE \eqref{sdeMND} is not available analytically,
one can still estimate the computational error as follows.
One can find the MAP connecting the initial stable equilibrium
and the saddle  e.g. using GMAM \cite{hey1,hey2} and find the quasipotential at the saddle by numerical integration
of the geometric action. The resulting value of the quasipotential at the saddle can be compared with the one 
obtained by linear interpolation from the computation by {\tt olim3D}.
The values of the quasipotential with respect to $x_i$ and $x_a$ at the saddle $x_s$ 
found by integration along the MAPs are $325.5$ and $4164.5$, respectively.
We have performed computations using {\tt olim3D}  with and without local factoring
and using various values of $K$. The local factoring radius is $R_{\rm fac} = 100$.
Good values for $K$ are 25 and 20 for computations initialized at $x_i$ and $x_a$, respectively.
The values of the quasipotential at $x_s$ obtained by linear interpolation from the computations
by {\tt olim3D}  initialized at $x_i $ and $x_a$ are 329.1 and 4121.1 respectively. Comparing these values with the values above
gives an estimate of 0.011 for the maximal relative error.

Local factoring reduces the error. 
For example, for the computation initialized at $x_i$ on a  $257\times6401\times513$ mesh,
without local factoring, the quasipotential at $x_s$ is 329.7. 
Comparing the errors at $x_s$ which are $3.6$ and $4.2$ for computations with and without local factoring respectively,
we see that local factoring reduces the error by 16.7\%. In both cases, the CPU time is 4.94 hours. 

Our computation of the quasipotential in 3D shows that the level sets within the basins of attraction are oar-shaped, and the paddles are
 twisted with respect to each other.
This means that the stochastic dynamics  on the timescale when the genetic switch is possible
are virtually limited to a neighborhood of a 2D manifold. One can find this manifold from the output of {\tt olim3D}
and perform a more accurate 2D computation of the quasipotential. 
We plan to develop a technique for such dimensionally-reduced computations in future work.

\subsection{Why larger update lengths lead to smaller errors}
\label{sec:why}
\begin{figure}[htbp]
\begin{center}
(a)\includegraphics[width = 0.7\textwidth]{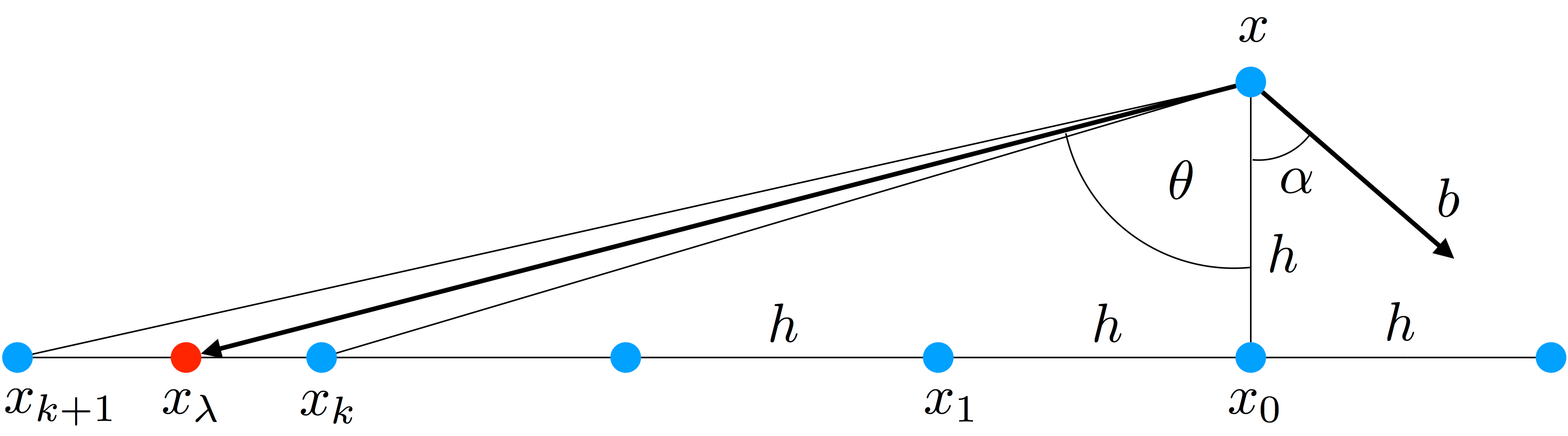}
(b)\includegraphics[width = 0.7\textwidth]{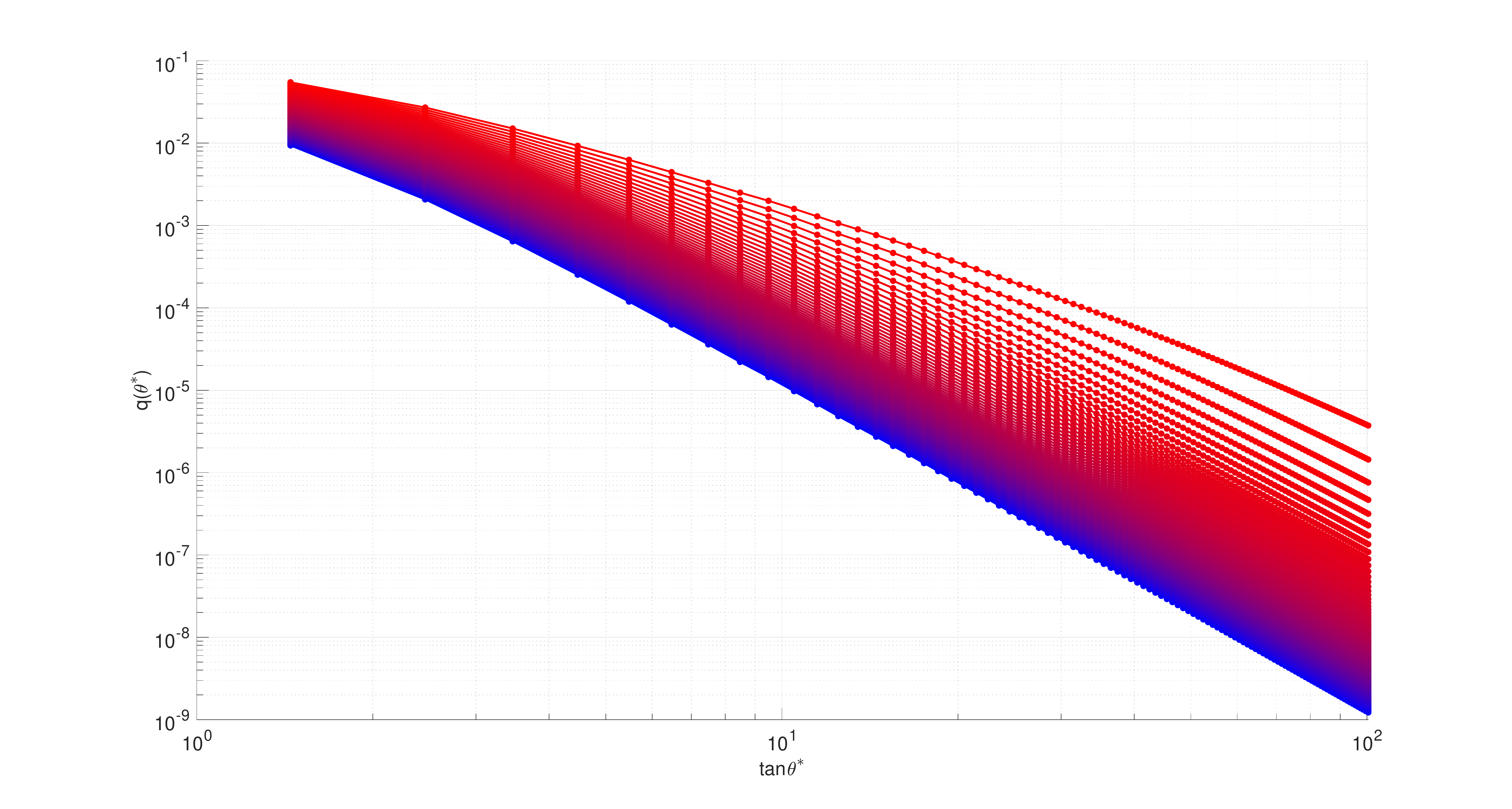}
\caption{An illustration for Section \ref{sec:why} that explains the phenomenon why larger update lengths may lead to smaller numerical errors.
(a): The set up of our 2D model. (b): Plots of $q(\theta)$ defined in  \eqref{B7} versus $\tan\theta$
for values of $\alpha$ from $\tfrac{\pi}{200}$ up to $\tfrac{\pi}{2} - \tfrac{\pi}{200}$ with step $\tfrac{\pi}{200}$. 
The bluest  curve corresponds to $\alpha = \tfrac{\pi}{200}$.
The reddest curve corresponds to $\alpha = \tfrac{\pi}{2} - \tfrac{\pi}{200}$.
}
\label{fig:AppB}
\end{center}
\end{figure}

 Figs. \ref{fig2:ex} (a) and (b) show the graphs of the normalized maximal absolute errors and normalized RMS errors as functions of $N$.
 As we have mentioned, the errors for Examples 3 and 5  with large $\xi$ 
 decay faster as $h$ decreases and surpass the errors for Examples 1, 2 and 4 where $\xi$ is relatively small.
 The least squares fits of these errors to $Ch^q$ where $h$ is the mesh step are shown in Table \ref{table:tests}.
 The convergence in Examples 3 and 5 is superlinear, while it is sublinear in Examples 1, 2 and 4.
 This phenomenon is explained as follows. 
 On average, larger $\xi$
 leads to larger  \emph{update lengths}; i.e., the distances $\|x-x_{{\lambda}^{\ast}}\|$ where $\lambda^{\ast}$ is the solution of the
triangle or simplex update minimization problems (Eqs. \eqref{2u} or \eqref{3u} respectively). If the mesh is fine enough, 
segments of MAPs of length $Kh$ are well approximated by straight line 
segments everywhere except for, perhaps, neighborhoods of equilibria.
The error due to the midpoint quadrature rule implemented in {\tt olim3D} decays as $O((Kh)^3)$, while the linear interpolation error
decays as $O(h^2)$. 
It turns out that the contribution of the interpolation error diminishes as the update length grows. 
This fact is responsible for both the visibly superlinear convergence rates and smaller numerical errors obtained 
on fine meshes in examples with large rotational components. 
We will explain this phenomenon for a 2D model shown in Fig. \ref{fig:AppB}(a).

Let the mesh points $\{x_0,~x_1, \ldots\}$ lying on a mesh line be {\sf Accepted Front}, and the point $x\in\mathcal{N}_1(x_0)$ 
lying on the next parallel mesh line be up for an update (see Fig. \ref{fig:AppB}(a)). Let $h$ be the mesh step.
We consider the case where the vector field $b$ is constant for simplicity. This assumption is justified as, provided that $b$ is continuous, the mesh 
can always be refined so that the variation of $b$ within the update radius $Kh$ is less than any prescribed small positive constant.
Let $b$ form an angle $\alpha$ with the ray $[x,x_0)$. The proposed one-point update value for $U(x)$ from the point $x_k$ is
\begin{equation}
\label{B1}
\mathsf{Q}_1(x_k,x) = U_k + \|b\|h\sqrt{k^2 + 1}\left(1 + \cos(\alpha + \arctan k)\right),
\end{equation}
where $U_k\equiv U(x_k)$. Here we have taken into account the fact that 
$$
\cos(\widehat{x_0x,b}) = \cos(\pi - (\alpha + \arctan k)) = \cos(\alpha + \arctan k).
$$
Suppose the numbers $k$, $U_k$, and $U_{k+1}$ are such  that the triangle update function 
\begin{equation} 
\label{B2}
f(\theta;\delta u): = U_k +h\|b\|\left( \delta u (\tan\theta - k) + \frac{1 + \cos(\alpha + \theta)}{\cos\theta}\right),~~
{\rm where}~~\delta u: = \tfrac{U_{k+1} - U_k}{\|b\|},
\end{equation}
 has a minimum in the interval $\arctan( k) < \theta < \arctan (k+1)$.  \eqref{B2} is triangle update  \eqref{2u} rewritten 
in terms of the angle $\theta$ defined as shown in Fig. \ref{fig:AppB}(a).  In terms of $\theta$, $\lambda = \tan(\theta) - k$.
After differentiating $f$, setting $f'$ to zero, and applying trigonometric formulas, we obtain:
\begin{equation}
\label{B3}
\frac{\partial f}{\partial \theta} = \|b\|h\frac{\delta u+ \sin\theta -\sin\alpha}{\cos^2\theta} = 0.
\end{equation} 
Hence, the minimizing angle $\theta^{\ast}$ must satisfy
\begin{equation}
\label{B4}
\theta^{\ast} = \arcsin\left(\sin\alpha - \delta u\right)\in [\arctan k,\arctan (k+1)].
\end{equation}

Now, we want to show that the difference between the triangle update, which is $f(\theta^{\ast};\delta u)$, and 
the minimal one-point update, which is $\min\{\mathsf{Q}_1(x_k,x),\mathsf{Q}_1(x_{k+1},x)\}$,
tends to zero as $k$ increases. We proceed as follows. 
For a fixed $\alpha\in\left(0,\tfrac{\pi}{2}\right)$ and each $k\in\mathbb{N}$ 
we find an interval of $\delta u$ such that the  function $f$ defined in  \eqref{B2} has a minimum at
$\theta^{\ast}\in[\arctan k,\arctan (k+1)]$:
\begin{equation}
\label{B5}
\sin\alpha -\frac{k}{\sqrt{k^2+1}}\le \delta u\le \sin\alpha -\frac{k+1}{\sqrt{(k+1)^2+1}}.
\end{equation}
Then we find the value of $\delta u$ in this interval such that the difference between the triangle update and the minimal value of the 
one-point update is maximal, i.e.
\begin{equation}
\label{B6}
\delta u^{\star} = \arg\max_{\delta u}\left[ \min\{|f(\theta^{\ast};\delta u) - \mathsf{Q}_1(x_k,x)|,|f(\theta^{\ast};\delta u) - \mathsf{Q}_1(x_{k+1},x)|\}\right].
\end{equation}
Since $\delta u^{\star}$ is uniquely defined for every $k$, we will write $\delta u^{\star}(k)$.
The minimizer of  $f(\theta;\delta u^{\star})$ is 
\begin{equation}
\label{B66}
\theta^{\star}: = \arcsin(\sin\alpha - \delta u^{\star}).
\end{equation}
The corresponding  $k^{\star} $ is equal to ${\rm floor}(\tan\theta^{\star})$.
Let us define
\begin{equation}
\label{B7}
q(\theta^{\star}) : = \frac{\min\{|f(\theta^{\star};\delta u^{\star}(k^{\star})) - \mathsf{Q}_1(x_k,x)|,|f(\theta^{\star};\delta u^{\star}(k^{\star})) 
- \mathsf{Q}_1(x_{k+1},x)|\}}{f(\theta^{\star};\delta u^{\star}(k^{\star}))}.
\end{equation}
In words, $q(\theta^{\star})$ is the upper bound of the 
quotient of the absolute value of the difference between the triangle update function and the minimal one-point update.
The graphs of $q(\theta^{\star})$ for $\alpha$ varying from $\tfrac{\pi}{200}$ to $\tfrac{\pi}{2}-\tfrac{\pi}{200}$ with step $\tfrac{\pi}{200}$
versus $\tan\theta^{\star}$ are shown   in Fig. \ref{fig:AppB}(b). These graphs show that
$q(\theta^{\star})$ decays not slower than  $(\tan\theta^{\star})^{-\gamma}$  as $\tan\theta^{\star}$ goes to infinity. The least squares fits for $50\le \tan\theta^{\star}\le 100$
 show that $\gamma> 3$.

Therefore, we conclude that the larger is the update length $l:= \|x_{\lambda} -x\|$, the closer are
the update values produced by a successful triangle update and the minimal one-point update.
The minimal one-point update, if the optimal angle $\theta$ is such that $\tan\theta$ is an integer, 
does not involve any interpolation error,
only errors that decay as $O(h^3)$ (see Section 4 in \cite{olim}). 
The update value produced by a successful triangle update  approaches the minimal one-point update value 
faster than $O((l/h)^{-3})$ as $(l/h)\rightarrow \infty$ . 
This sheds light  on the phenomenon of
numerical errors decaying faster in systems where update lengths tend to be larger.


\section{Conclusion}
\label{sec:conclusion}
We have presented a 3D quasipotential solver {\tt olim3D}. This is the first time  the quasipotential has been computed 
on a regular 3D mesh to the best of our knowledge. The C source code is available on M. Cameron's web page \cite{mariakc}.
Our method results from the further development of ordered line integral methods introduced in \cite{olim}.
They, in turn, inherit the general structure from the Dijkstra-like Hamilton-Jacobi solver ordered upwind method \cite{OUM2001,OUM2003}.

An important feature of {\tt olim3D} is that it computes the quasipotential on rather large meshes 
(e.g. $513^3$) within reasonable CPU times; i.e., within a few hours -- a straightforward promotion of
the 2D OLIM-Midpoint from \cite{olim} to 3D would take several days. Such a dramatic reduction of CPU 
time is the result of a number of technical innovations implemented in {\tt olim3D}: 
$(i)$ the new version of the hierarchical update strategy applied to \emph{all}
{\sf Considered } points,
$(ii)$ the use of the Karush-Kuhn-Tucker constrained optimization theory to  reject simplex updates that are going to be unsuccessful, and
$(iii)$ pruning of the number of admissible simplexes and a fast search for them.

We have conducted an extensive numerical study of the proposed solver. Our set of examples includes a genetic switch model \cite{lv},
systems with hyperbolic periodic orbits \cite{tao}, and a series of linear and nonlinear SDEs with different ratios
of magnitudes of their rotational and potential components.  
The practical results of this study are: 
$(i)$ a guideline for choosing the update parameter $K$,  $(ii)$ the conclusion that a local factoring
may or may not be helpful, depending on the problem, and $(iii)$ a way to estimate the numerical 
error using the integration along the MAP where the quasipotential is not given analytically.
 An interesting finding is that longer update lengths, occurring in systems with 
the ratio of the magnitudes of the rotational and  potential components is around 10, may lead to higher convergence rate than in systems
where this ratio is around 1. 

Our method is a new tool for the analysis and visualization of 3D dynamical systems perturbed by small white noise.
It can be used to determine the relative stability of attractors \cite{nolting}, and  
to find manifolds to which the dynamics of a system are virtually confined to.

The next stage of the development of quasipotential solvers
is motivated by the fact that a number of interesting real-life systems 
have a significant time-scale separation virtually bounding their high-dimensional stochastic dynamics to
a neighborhood a low-dimensional manifold \cite{aladip}. Our application to the genetic switch model \cite{lv}
exhibits this phenomenon: the level sets of the quasipotential are thin and stretched along a 2D manifold. 
Inspired by these facts, we will explore combining quasipotential solvers with techniques for learning manifolds from the dynamics -- e.g., as in  \cite{atlas} --
in our future work.


\section*{Acknowledgements}
We thank Prof. A. Vladimirsky for useful discussions.
This work was partially supported by M. Cameron's NSF Career grant DMS1554907.


%
%

 \appendix
\setcounter{equation}{0}
\renewcommand{\theequation}{A-\arabic{equation}}
\section*{Appendix A}
{\bf  Calculation of gradient \eqref{grad3} and Hessian \eqref{hess3}.}
%
The calculation of $\nabla f$ for $f(\lambda)$ given by  \eqref{fmin} involves the following ingredients: 
$\nabla  \|b_{m\lambda}\| $, $\nabla \|x-x_{\lambda}\|$, and $\nabla [b_{m\lambda}\cdot(x-x_{\lambda})]$.
Recalling that (Fig. \ref{fig:simplex})
\begin{align*}
b_{m\lambda} &= b(x_{m0}) + \lambda_1\left[ b(x_{m1})-b(x_{m0})\right] + \lambda_2\left[b(x_{m2}) - b(x_{m0})\right],\\
x-x_{\lambda} & = x - x_0 + \lambda_1(x_0-x_1) + \lambda_2(x_0-x_2),
\end{align*}
we compute:
\begin{align}
\nabla \|b_{m\lambda}\| &=\nabla \sqrt{(b_{m\lambda})_1^2 + (b_{m\lambda})_2^2 + (b_{m\lambda})_3^2} \notag \\
&=
\frac{1}{\|b_{m\lambda}\|}\left[\begin{array}{c}b_{m\lambda}\cdot(b_{m1}-b_{m0})\\b_{m\lambda}\cdot(b_{m2}-b_{m0})\end{array}\right]
=\frac{1}{\|b_{m\lambda}\|}B^\top b_{m\lambda}, \label{gb}
\end{align}
where $B$ is the $3\times 2$ matrix defined in  \eqref{Bmatr}.
A similar calculation gives:
\begin{equation}
\label{gx}
\nabla \|x-x_{\lambda}\| = \frac{1}{\|x-x_{\lambda}\|} X^\top (x-x_{\lambda}),
\end{equation}
where $X$ is $3\times 2$ matrix defined in  \eqref{Xmatr}.
The gradient of the dot product is
\begin{align}
\nabla [b_{m\lambda}\cdot (x-x_{\lambda})] &= \nabla\left[(b_{m\lambda})_1(x-x_{\lambda})_1 + (b_{m\lambda})_1(x-x_{\lambda})_2 + (b_{m\lambda})_1(x-x_{\lambda})_3\right]\notag \\
&=\left[\begin{array}{c}
(b_{m1}-b_{m0})\cdot (x-x_{\lambda}) + b_{m\lambda}\cdot (x_0-x_1)\\
(b_{m2}-b_{m0})\cdot (x-x_{\lambda}) + b_{m\lambda}\cdot (x_0-x_2)
\end{array}\right] = B^\top (x-x_{\lambda}) + X^\top b_{m\lambda}. \label{gdot}
\end{align}
Assembling the terms  from Eqs. \eqref{gb}--\eqref{gdot} and adding $\nabla U_{\lambda}$ we obtain  $\nabla f$ given by  \eqref{grad3}.

To compute the Hessian, we calculate $\frac{\partial^2 f}{\partial \lambda_1^2}$ and 
$\frac{\partial^2 f}{\partial \lambda_1\partial\lambda_2}$. The other entries of $H$ 
are deduced using symmetry. 
\begin{align}
\frac{\partial^2 f}{\partial \lambda_1^2} &= 2\frac{[(x_0-x_1)\cdot(x-x_{\lambda})][(b_{m1}-b_{m0})\cdot b_{m\lambda}]}{\|x-x_{\lambda}\| \|b_{m\lambda}\|} \notag \\
&- \frac{\|x-x_{\lambda}\|}{\|b_{m\lambda}\|^3}[(b_{m1}-b_{m0})\cdot b_{m\lambda}]^2 \notag \\
&+ \frac{ \|x-x_{\lambda}\|}{\|b_{m\lambda}\|}(b_{m1}-b_{m0})\cdot (b_{m1}-b_{m0})\notag \\
& - 
\frac{\|b_{m\lambda}\|}{\|x-x_{\lambda}\|^3}[(x_0-x_1)\cdot(x-x_{\lambda})]^2 \notag \\
&+ \frac{\|b_{m\lambda}\|}{ \|x-x_{\lambda}\|}(x_{0}-x_{1})\cdot (x_0-x_1)]\notag \\
&-[(b_{m1}-b_{m0})\cdot (x_0-x_1) +(b_{m1}-b_{m0})\cdot (x_0 - x_1)]. \label{d11}
\end{align}

\begin{align}
\frac{\partial^2 f}{\partial \lambda_1\partial \lambda_2} &= \frac{[(x_0-x_2)\cdot(x-x_{\lambda})][ (b_{m1}-b_{m0})\cdot b_{m\lambda}] +
[(x_0-x_1)\cdot(x-x_{\lambda})][ (b_{m2}-b_{m0})\cdot b_{m\lambda}] 
}{\|x-x_{\lambda}\| \|b_{m\lambda}\|} \notag \\
&-
\frac{\|x-x_{\lambda}\|}{\|b_{m\lambda}\|^3}[(b_{m1}-b_{m0})\cdot b_{m\lambda}][(b_{m2}-b_{m0})\cdot b_{m\lambda}] \notag \\
&+ \frac{ \|x-x_{\lambda}\|}{\|b_{m\lambda}\|}(b_{m1}-b_{m0})\cdot (b_{m2}-b_{m0}) \notag \\
&-
\frac{\|b_{m\lambda}\|}{\|x-x_{\lambda}\|^3}[(x_0-x_1)\cdot(x-x_{\lambda})][(x_0-x_2)\cdot(x-x_{\lambda})] \notag \\
&+ \frac{\|b_{m\lambda}\|}{ \|x-x_{\lambda}\|}(x_{0}-x_{1})\cdot (x_0-x_2) \notag \\
&-[(b_{m1}-b_{m0})\cdot (x_0-x_2) +(b_{m2}-b_{m0})\cdot (x_0 - x_1)]. \label{d12}
\end{align}
Putting Eqs. \eqref{d11} and \eqref{d12} together term-by-term and recognizing the matrices $B$ and $X$
we obtain  \eqref{hess3}. 



%
%

 \appendix
\setcounter{equation}{0}
\renewcommand{\theequation}{B-\arabic{equation}}
\renewcommand{\theexample}{B-\arabic{example}}
\section*{Appendix B}
{\bf On the Hessian of the simplex update function.}
Since it does not complicate our calculations, we consider a more  general situation where $x_{mi}$ are replaced with
$$
y_i = x_i + \theta(x_i - x),\quad \theta \in [0,1),\quad i=0,1,2.
$$
Let us denote the linear interpolant of $b$ in the triangle $(y_0,y_1,y_2)$ by $b_{\theta}$:
$$
b_{\theta}(\lambda) = b(y_0) + B_{\theta}\lambda,\quad{\rm where}\quad B_{\theta} : = [b(y_1)-b(y_0),b(y_2)-b(y_0)].
$$
It will replace $b_{m\lambda}$ in  \eqref{hess3}.
Taylor expanding $b_\theta$ near $\lambda = 0$ we obtain
\begin{equation}
\label{btaylor}
b_{\theta}(\lambda) = b_{\theta}(0) + JY\lambda + O(h^2),
\end{equation}
where $J$ is the $3\times 3$ Jacobian matrix of the vector field $b$, $Y:=[y_1-y_0,y_2-y_0]$ is a $3\times 2$ matrix, and $h$ is the size of the mesh step.
Note that the entries of $Y$ are of the order of $h$.  
We observe that $Y = -\theta X$ where $X$ is given by  \eqref{Xmatr}.
Hence the matrix $B_{\theta}$ can be approximated by
\begin{equation}
\label{Btheta}
B_{\theta} = -\theta JX + O(h^2).
\end{equation}
Plugging $B_{\theta}$ from  \eqref{Btheta} for $B$ in  \eqref{hess3}, using the notations
$$
\beta : = \frac{ \|b_{\theta}(\lambda)\|}{\|x-x_{\lambda}\|},\quad 
p:=\frac{x-x_{\lambda}}{\|x-x_{\lambda}\|},\quad 
q:=\frac{b_{\theta}(\lambda)}{\|b_{\theta}(\lambda)\|},
$$
and grouping terms, we obtain
\begin{align}
H_{\theta} & =-\theta (X^\top pq^\top JX + X^\top J^\top qp^\top X)  + \frac{\theta^2}{\beta} X^\top J^\top JX + \beta X^\top X \notag\\
&-\frac{\theta^2}{\beta}  X^\top J^\top q q^\top JX
-\beta X^\top pp^\top X +\theta\left(X^\top J^\top X + X^\top JX \right) +O(h^3) \notag \\
& = X^\top\left[ \theta\left[ (I - pq^\top)J + J^\top(I - qp^\top)\right]  + \beta (I - pp^\top) + \frac{\theta^2}{\beta}J^\top(I - qq^\top)J \right]X +O(h^3). \label{Htheta}
\end{align}
Typically, $\beta\gg 1$ since $\|x-x_{\lambda}\| \lesssim Kh$, while $\|b_\theta(\lambda)\| = O(1)$. 
Hence, if $K$ proportional  to $h^{-\alpha}$ for $\alpha \in(0,1)$ and $b_{\theta}(\lambda)\neq 0$, then $\beta\rightarrow \infty$ as $h\rightarrow 0$.
Therefore, as $h\rightarrow 0$, 
$$
H_{\theta}\approx \beta X^\top (I - pp^\top) X
$$
becomes positive definite for most simplexes, except for those where $\|b\|$ is small, 
i.e., of the order of $h$, near equilibria.

Let us demonstrate that the Hessian can have negative eigenvalues near equilibria on the following 2D example.
\begin{example}
\label{ex:negH}
In 2D, the Hessian is $1\times 1$, just the $H_{11}$ component of $H$ in  \eqref{hess3}.
Suppose $b$ is the gradient field 
$b(x) = -x$. Then $J = -I$ and the quasipotential is $U(x) = \|x\|^2$.
Let us consider a triangle update with
$$
x_0 =h\left[ \begin{array}{c}1\\0\end{array}\right],\quad 
x_1 =h\left[ \begin{array}{c}1\\-1\end{array}\right],\quad{\rm and}\quad
x =h\left[ \begin{array}{c}2\\-1\end{array}\right],
$$
as shown in Fig. \ref{fig:negH}(a). 
\begin{figure}[htbp]
\begin{center}
(a)\includegraphics[width = 0.45\textwidth]{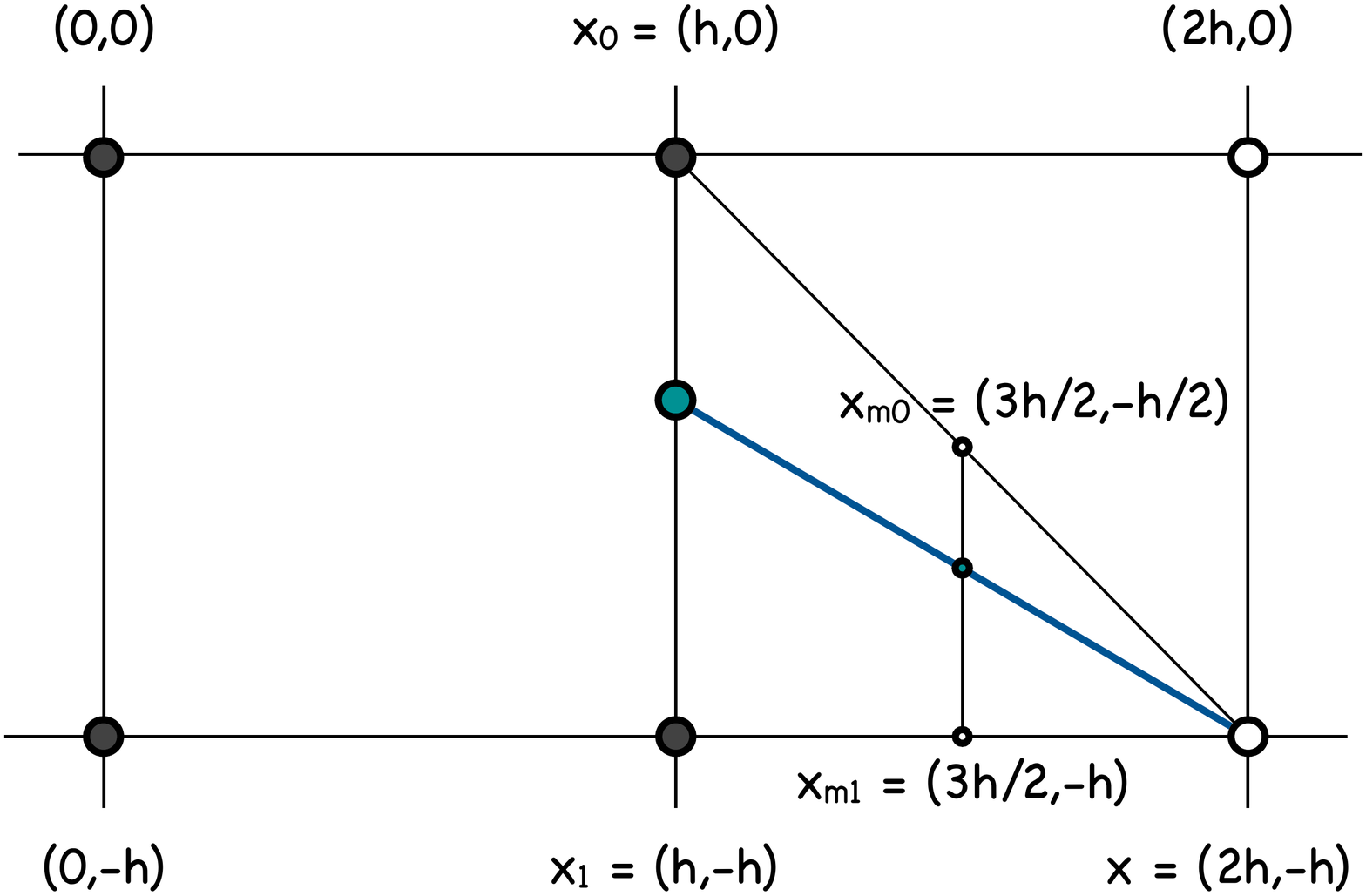}
(b)\includegraphics[width = 0.45\textwidth]{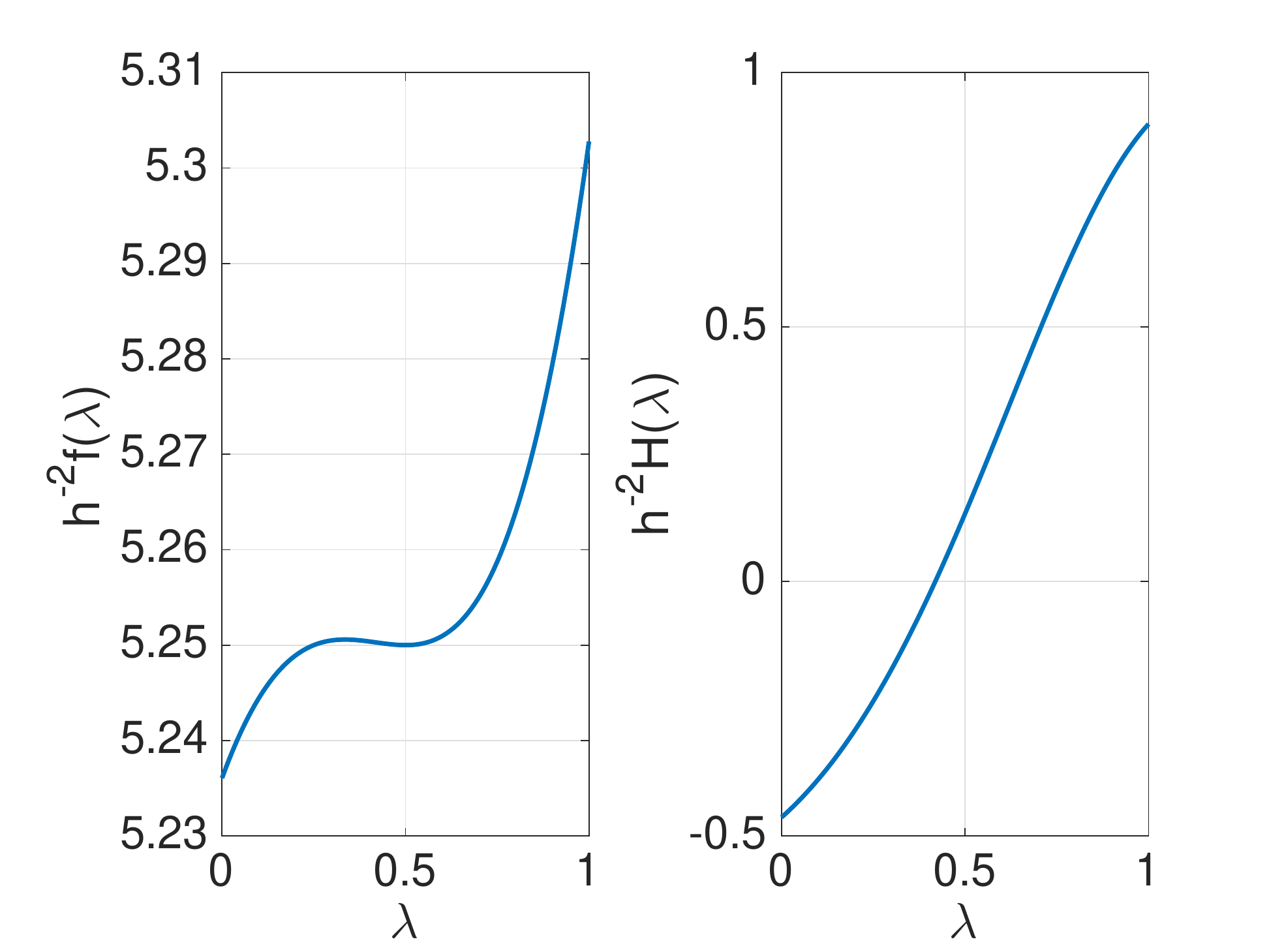}
\caption{An illustration for Example \ref{ex:negH}. (a) The triangle update.
(b) The graphs of $h^{-2}f(\lambda)$ and $h^{-2}H(\lambda)$.  
}
\label{fig:negH}
\end{center}
\end{figure}
The function $f(\lambda)$ defined in  \eqref{fmin} and its Hessian are
\begin{equation*}
f(\lambda) = h^2\left[1+\lambda + \frac{1}{2}\left(\left[1 + (1-\lambda)^2\right]^{1/2}\left[9+(\lambda+1)^2\right]^{1/2} + 2 - \lambda^2\right)\right],
\end{equation*}
\begin{align*}
H(\lambda) &= \frac{h^2}{2} \left[\frac{R_1}{R_2} +
\frac{R_2} {R_1} 
+\frac{2(\lambda^2-1)}{R_1 R_2 }
-\frac{(1+\lambda)^2 R_1}{R_2^3} 
-  \frac{(1 -\lambda)^2 R_2 }{R_1^3} - 2\right],\\
&{\rm where}~~ R_1: = \left[1 + (1-\lambda)^2\right]^{1/2},\quad R_2: = \left[9+(\lambda+1)^2\right]^{1/2}.
\end{align*}
The graphs of $h^{-2}f(\lambda)$ and $h^{-2}H(\lambda)$ are shown in Fig. \ref{fig:negH}(b).
The Hessian $H(\lambda)$ changes sign in the interval $\lambda\in[0,1]$. Hence it is not positive definite.
The global minimum of $f(\lambda)$ subject to $0\le \lambda\le 1$ is achieved at $\lambda=0$ and would be found using the one-point update;
$f(\lambda)$ also has a local minimum at $\lambda = 0.5$ corresponding to the true MAP arriving at $x$ from the origin. However, the linear interpolation
$$
u_{\lambda} = h^2(u_0 + \lambda(u_1 - u_0)) = h^2(1 + \lambda)
$$
exaggerates  $u$ on the segment $[x_0,x_1]$. As a result, the OLIM  picks the one-point update from $x_0$ 
because it turns out to be the smallest proposed update value.
The true value at $x$ is $U(x) = \|x\|^2 = 5h^2$. 
The one-point update value, $5.236h^2$, is closer to the true value than $f(0.5) = 5.250h^2$, as
obtained from the triangle update with the value of $\lambda$ corresponding to the true MAP.

\end{example}




\begin{thebibliography}{9}
%
\bibitem{Aurell} 
E. Aurell and K. Sneppen. Epigenetics as a First Exit Problem. Physical Review Letters. 88: 048101(1-4), (2002)

\bibitem{BS}
R. Bartels and G.~W. Stewart, 
Solution of the matrix equation AX+ XB = C,
 Comm A.C.M., {\bf 15}, 9, 820--826 (1972)

\bibitem{bouchet}
F. Bouchet and J. Reygner, {\it Generalisation of the Eyring-Kramers Transition Rate Formula to Irreversible Diffusion Processes}, 
Annales Henri Poincare, {17} (2016), 12, pp. 3499--3532
%
%
\bibitem{quasi}
M.~ K. Cameron, {\it Finding the Quasipotential for Nongradient SDEs},
Physica D: Nonlinear Phenomena, { 241} (2012), pp. 1532--1550
%
\bibitem{mariakc}
{\tt https://www.math.umd.edu/~mariakc/software-and-datasets.html}
%
%
%
%

\bibitem{CV}
A. Chacon and A. Vladimirsky,
Fast two-scale methods for Eikonal equations,
SIAM J. on Scientific Computing {\bf 34}, 2, A547--A578 (2012)

\bibitem{CF}
Z. Chen, M. Freidlin, Smoluchowski-Kramers approximation and exit problems, Stoch. Dyn. {\bf 5}, 4,  569--585 (2005),

\bibitem{chen}
Z. Chen, 
Asymptotic Problems related to Smoluchowski-Kramers approximation, Ph.D. Dissertation, UMD, 2006 
\href{https://drum.lib.umd.edu/bitstream/handle/1903/3791/umi-umd- 3634.pdf?sequence=1}
{https://drum.lib.umd.edu/bitstream/handle/1903/3791/umi-umd- 3634.pdf?sequence=1}

\bibitem{FitzHugh-Nagumo}
Z. Chen, J. Zhu, X.  Liu,  Crossing the quasi-threshold manifold of a noise-driven excitable system. Proc. R. Soc. A 473: (2017) 0058.

\bibitem{visc}
M.~G.Crandall, P.~L. Lions, 
Viscosity solutions of Hamilton-Jacobi-Bellman equations, 
Trans. Am. Math. Soc. {\bf 277}, 1--43 (1983)

\bibitem{atlas}
M. Crosskey and M. Maggioni,
ATLAS: a geometric approach to learning high-dimensional stochastic systems near manifolds,
SIAM J. Multiscale Model. Simul., {\bf 15}, 1, 110--156 (2017) 


\bibitem{olim}
D. Dahiya, and M. Cameron, Ordered Line Integral Methods for Computing the quasipotential, J. Sci. Comput. (2017), {\it to appear},  
\href{https://doi.org/10.1007/s10915-017-0590-9}{https://doi.org/10.1007/s10915-017-0590-9}

\bibitem{DC2}
D. Dahiya, and M. Cameron, 
An Ordered Line Integral Method for Computing the quasipotential in the case of Variable Anisotropic Diffusion,
Physica D (2018) {\it to appear}, https://doi.org/10.1016/j.physd.2018.07.002, arXiv:1806.05321


%

\bibitem{mam}
W. E, W. Ren, and E. Vanden-Eijnden, Minimum Action Method for the Study of Rare Events,
Comm. Pure Appl. Math., {\bf 57}, 0001--0020 (2004)


\bibitem{FLZ}
S. Fomel, S. Luo, and H. Zhao, {\it Fast sweeping method for the factored eikonal equation}, J. Comput. Phys., {\bf 228}, 17, 6440--6455 (2009)


\bibitem{FW} 
M.~I. Freidlin and A.~D. Wentzell, {\it Random Perturbations of Dynamical Systems}, 3rd
Ed, Springer-Verlag, Berlin Heidelberg, 2012.
%
%

\bibitem{hey1}
M. Heymann, E. Vanden-Eijnden, 
Pathways of maximum likelihood for rare events in non-equilibrium systems, application to nucleation in the presence of shear, 
Phys. Rev. Lett. {\bf 100}, 14, 140601  (2007)

\bibitem{hey2}
M. Heymann, E. Vanden-Eijnden, 
The geometric minimum action method: a least action principle on the space of curves, 
Comm. Pure Appl. Math. {\bf 61}, 8, 1052--1117  (2008)



\bibitem{ishii}
H. Ishii,  A simple direct proof of uniqueness for solutions of the Hamilton-Jacobi equations of eikonal type, 
Proc. Amer. Math. Soc. {\bf 100}, 2, 247--251 (1987)
%
%
\bibitem{LQ}
S. Luo and J. Qian, {\it Fast sweeping methods for factored anisotropic eikonal equations: multiplicative and additive factors}, 
J. Sci. Comput. {\bf 52}, 2, 360--382 (2012)


\bibitem{lv}
Cheng Lv, Xiaoguang Li, Fangting Li1, and Tiejun Li, 
Constructing the Energy Landscape for Genetic Switching System Driven by Intrinsic Noise,
\href{http://journals.plos.org/plosone/article?id=10.1371/journal.pone.0088167}{PLOS One, {\bf 9}, 2, e88167 (2014)}


\bibitem{NieNoise}
Q. Wang, W.~R. Holmes, J. Sosnik, T. Schilling, and Q. Nie, 
Cell Sorting and Noise-Induced Cell Plasticity Coordinate to Sharpen Boundaries between Gene Expression Domains, 
PLoS Comput. Biol. {\bf 13},1, e1005307. doi:10.1371/journal. pcbi.1005307

\bibitem{nocedal}
J. Nocedal and S.~J. Wright, {\it Numerical Optimization}, Second Edition, Springer, USA, 2006 

\bibitem{nolting}
B.~C. Nolting, K.~C. Abbot, Balls, cups, and quasipotentials: quantifying stability in stochastic systems,
Ecology 97(4), 850--864 (2016)

\bibitem{Rjournal}
B. Nolting, C.  Moore, C. Stieha, M. Cameron, K. Abbott,
QPot: an R package for stochastic differential
equation quasipotential analysis. R J. 8(2), 19--38 (2016)

\bibitem{NLSE}
G. Poppe and T. Schaefer,
Computation of minimum action paths of the stochastic nonlinear Schroedinger equation with dissipation,
arXiv:1804.10142

\bibitem{PC}
S. Potter and M. Cameron, Ordered Line Integral Methods for Solving the Eikonal Equation, {\it in preparation}

\bibitem{aladip}
M.~A. Rohdanz, W. Zheng, M. Maggioni, and C. Clementi,
Determination of reaction coordinates via locally scaled diffusion map,
The Journal of Chemical Physics 134, 124116 (2011); doi: 10.1063/1.3569857


\bibitem{QV}
D. Qi and A. Vladimirsky.
Corner cases, singularities, and dynamic factoring, arXiv:1801.04322v1




\bibitem{SethPNAS} 
J.~A. Sethian, A fast marching level set method for monotonically advancing
fronts, Proc. Natl. Acad. Sci. 93 (4) (1996) 1591--1595.

\bibitem{SethSIAM} 
J.~A. Sethian, Fast marching methods, SIAM Rev. 41 (2) (1999) 199--235.

\bibitem{SethBook}
J.~A. Sethian, Level Set Methods and Fast Marching Methods, Cambridge
University Press, 1999.

\bibitem{OUM2001}
J.~.A. Sethian, A. Vladimirsky, Ordered Upwind Methods for static Hamilton-Jacobi-Bellman equations, 
Proc. Natl. Acad. Sci. {\bf 98}, 11069--11074.  (2001)

\bibitem{OUM2003}
 J.~A. Sethian, A. Vladimirsky, Ordered Upwind Methods for static Hamilton-Jacobi-Bellman equations: theory and algorithms, 
 SIAM J. Numer. Anal.
{\bf 41}, 1, 325--363  (2003) 

\bibitem{sunzhou}
Y. Sun and X. Zhou, 
An Improved Adaptive Minimum Action Method for the Calculation of Transition Path in nongradient Systems,
Commun. Comput. Phys.{\bf 24}, 1, 44--68 (2018) doi: 10.4208/cicp.OA-2016-0230
%
\bibitem{tao}
M. Tao, {\it Hyperbolic periodic orbits in nongradient systems and small-noise-induced metastable transitions}, Physica D: Nonlinear Phenomena,
363 (2018) pp. 1-17

\bibitem{DV2}
H. Touchette, The large deviation approach to statistical mechanics. Phys. Rep. {\bf 478},  1--69 (2009)

\bibitem{DV1}
S. Varadhan, Large deviations and applications, SIAM, Philadelphia, 1984

\bibitem{savanna}
J.~D. Touboul, A.~C. Staver, and S.~A. Levin,
On the complex dynamics of savanna landscapes,
Proc. Natl. Acad. Sci. {\bf 115}, 7, E1336--E1345 (2018), https://doi.org/10.1073/pnas.1712356115

\bibitem{wilkinson}
J. Wilkinson, Two Algorithms Based on Successive Linear Interpolation, Computer Science, Stanford University, Technical Report CS-60 (1967)


\bibitem{zhou1}
X. Zhou, W. Ren, Weiqing, W. E, Adaptive minimum action method for the study of rare events,
J. Chem. Phys. 128, 104111 (2008)

\bibitem{zhou2}
X. Zhou, W. Ren, Weiqing, W. E, Study of noise-induced transitions in the Lorenz system using the minimum
action method. Commun. Math. Sci. 8(2), 341--355 (2010)


\end{thebibliography}
\end{document}